\documentclass[12pt]{amsart}
\usepackage{amscd,amssymb}
\usepackage[dvips]{graphics}
\input diagrams
\diagramstyle[PostScript=dvips]

\newcommand{\includefig}[1]{\raisebox{-3ex}{\resizebox{!}{7ex}{\includegraphics{pix/#1.eps}}}}

\newcommand{\tconkm}{\includefig{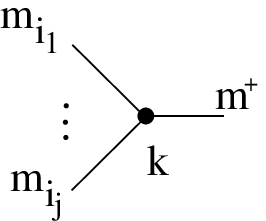}}
\newcommand{\tcongi}{\includefig{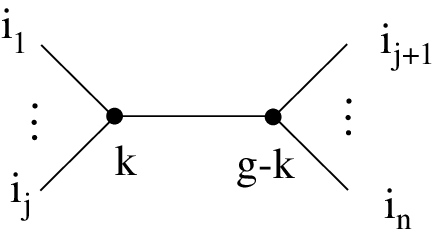}}

\newcommand{\tcongm}{\includefig{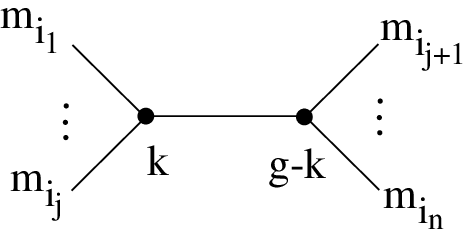}}
\newcommand{\tcongmpm}{\includefig{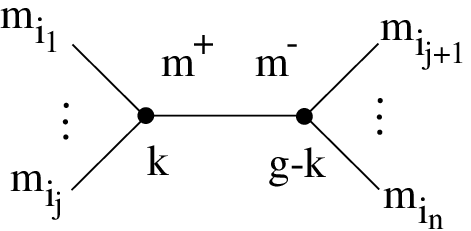}}
\newcommand{\ocongmpm}{\includefig{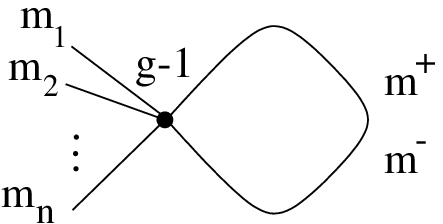}}
\newcommand{\tcongmkm}{\includefig{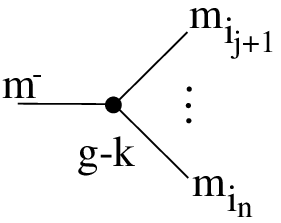}}

\newcommand{\ocongm}{\includefig{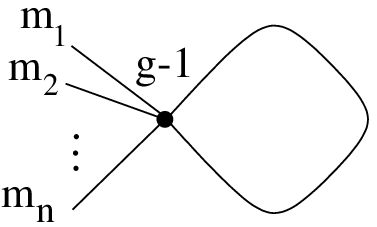}}

\newcommand{\ocongi}{\includefig{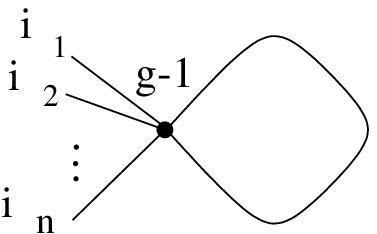}}
\newcommand{\ocongmcut}{\includefig{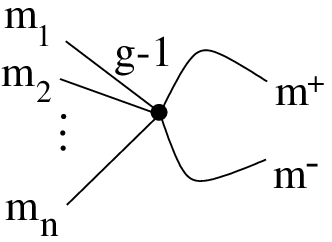}}

\newcommand{\scheme}[1]{\mathcal {#1}}
\newcommand{\mbar}{\M}

\newcommand{\mgnbar}{\overline{\scheme{M}}_{g,n}}
\newcommand{\mgnrmbar} {\mgnbar^{1/r,\bm}}
\newcommand{\mgnrbar} {\mgnbar^{1/r}}

\newcommand{\mgnv}{\mgnbar(V)}
\newcommand{\cgn}{\scheme{C}_{g,n}}
\newcommand{\cgnr}{\cgn^{1/r}}
\newcommand{\cgnrm}{\cgn^{1/r,\bm}}
\newcommand{\mgnrmbareven}{\mgnbar^{1/r,\bm,\text{even}}}
\newcommand{\mgnrmbarodd}{\mgnbar^{1/r,\bm,\text{odd}}}

\newcommand{\mrbar}{\M^{1/r}}
\newcommand{\pic}{\text{Pic\,}}        

\newcommand{\tensor}{\otimes}
\newcommand{\cross}{\times}

\newcommand{\irightarrow}{\rTo^{\sim}}

\newcommand{\ck}{{\mathcal K}}
\newcommand{\cl}{{\mathcal L}}
\newcommand{\ce}{{\mathcal E}}
\newcommand{\cf}{{\mathcal F}}
\newcommand{\co}{{\mathcal O}}
\newcommand{\cc}{{\mathcal C}}

\newcommand{\ch}{{\mathcal H}}
\newcommand{\cd}{{\mathcal D}}

\newcommand{\cm}{{\mathcal M}}

\newcommand{\sheafhom}{\ch\kern-.15em om}

\newcommand{\Aut}{\operatorname{Aut}} 
\newcommand{\bm}{\mathbf{m}}          %
\newcommand{\bmu}{\boldsymbol{\mu}}
\newcommand{\bnu}{\boldsymbol{\nu}}
\newcommand{\bs}{\mathbf{s}}
\newcommand{\bt}{\mathbf{t}}
\newcommand{\btau}{\boldsymbol{\tau}}

\newcommand{\bu}{\mathbf{u}}
\newcommand{\bx}{\mathbf{x}}          
\newcommand{\cft}{CohFT}              
\newcommand{\cfts}{CohFTs}            
\newcommand{\chh}{\hat{\ch}^{(r)}}          
\newcommand{\chr}{\ch^{(r)}}          
\newcommand{\cor}[1]{\langle\,{#1}\,\rangle}  
\newcommand{\ccor}[1]{\langle\langle\,{#1}\,\rangle\rangle}  
\newcommand{\cv}{c^{1/r}}            
\newcommand{\cvtwo}{c^{1/2}}           
\newcommand{\dil}{\cd}                
\newcommand{\ev}{\mathrm{ev}}         
\newcommand{\etat}{\widetilde{\eta}}  

\newcommand{\KdV}{\mathrm{KdV}}
\newcommand{\Lambdat}{\widetilde{\Lambda}} 
\newcommand{\M}{\overline{\MM}}       
\newcommand{\MM}{\scheme{M}}          

\newcommand{\nc}{{\mathbb{C}}}        
\newcommand{\nq}{{\mathbb{Q}}}        
\newcommand{\nz}{{\mathbb{Z}}}        
\newcommand{\pit}{\tilde{\pi}}     

\def\lift#1#2{
  \dimen0 = \unitlength
  \multiply\dimen0 by #1 \divide \dimen0 by 2
  \dimen1 = \dimen0
  \multiply \dimen1 by 7 \divide \dimen1 by 10
  \raise\dimen1
     \hbox{\hskip 0.3cm ${\vbox to \dimen0{}}$ \enspace #2}}

\newtheorem{thm}{Theorem}[section]
\newtheorem{lm}[thm]{Lemma}
\newtheorem{prop}[thm]{Proposition}
\newtheorem{crl}[thm]{Corollary}
\newtheorem{conj}[thm]{Conjecture}

\theoremstyle{definition}
\newtheorem{rem}[thm]{Remark}
\newtheorem{rems}[thm]{Remarks}

\newtheorem{df}[thm]{Definition}
\newtheorem{ex}[thm]{Example}
\theoremstyle{remark}

\newenvironment{ack}{\begin{sloppypar}\emph{Acknowledgments}:}{\end{sloppypar}}

\begin{document}
\addtocounter{section}{-1}
\title[Higher Spin Curves and Integrable Hierarchies]
{Moduli Spaces of Higher Spin Curves and Integrable Hierarchies}

\author
[T. J. Jarvis]{Tyler J. Jarvis}
\address
{Department of Mathematics, Brigham Young University, Provo, UT 84602, USA}
\email{jarvis@math.byu.edu}
\thanks{Research of the first author was partially supported by NSA
grant MDA904-99-1-0039}

\author
[T. Kimura]{Takashi Kimura}
\address
{Department of Mathematics, 111 Cummington Street, Boston
University, Boston, MA 02215, USA}
\email{kimura@math.bu.edu}
\thanks{Research of the second author was partially supported by NSF grant
  number DMS-9803427}

\author
[A. Vaintrob]{Arkady Vaintrob}
\address
{ Max-Planck-Institut f\"ur Mathematik, Vivatsgasse 7, 53111
Bonn, Germany \quad and \quad
Department of Mathematical Sciences, New
Mexico State University, Las Cruces, NM 88003, USA}
\email{vaintrob@math.nmsu.edu}

\date{\today}

\begin{abstract}
We prove the genus zero part of the generalized Witten conjecture,
relating moduli spaces of higher spin curves to Gelfand-Dickey
hierarchies. That is, we show that intersection numbers on the
moduli space of stable $r$-spin curves assemble into a generating
function which yields a solution of the semiclassical limit of the
$\KdV_r$ equations. We formulate axioms for a cohomology class on
this moduli space which allow one to construct a cohomological
field theory of rank $r-1$ in all genera. In genus zero it
produces a Frobenius manifold which is isomorphic to the Frobenius
manifold structure on the base of the versal deformation of the
singularity $A_{r-1}$. We prove analogs of the puncture, dilaton,
and topological recursion relations by drawing  an analogy with
the construction of Gromov-Witten invariants and  quantum
cohomology.
\end{abstract}
\maketitle

\section{Introduction}
\label{intro}

The moduli space $\mgnbar$ of stable curves of genus $g$ with $n$
marked points is a fascinating object. Mumford \cite{Mu}
introduced tautological cohomology classes associated to the
universal curve $\cgn \rTo \mgnbar$. Witten \cite{W3} conjectured
and Kontsevich \cite{Ko} proved that certain intersection numbers
of tautological cohomology classes on $\mgnbar$ have a generating
function which satisfies the equations of the Korteweg-de Vries
hierarchy (more precisely, that it is a $\tau$-function of the KdV
hierarchy satisfying some additional equations). This remarkable
result provided an unexpected link between the algebraic geometry
of these moduli spaces and integrable systems.

The spaces $\mgnbar$ can be generalized in two ways. The first way
is by choosing a smooth projective variety $V$ and considering the
moduli space $\M_{g,n}(V)$ of stable maps into $V$ from genus-$g$,
$n$-pointed, stable curves. When $V$ is a point, $\mgnv$ reduces
to $\mgnbar$.

The second way to generalize $\mgnbar$ is by considering the
moduli space $\mgnrbar$ of higher spin curves introduced
in~\cite{J,J2}.  Roughly speaking, a higher spin curve, or
$r$-spin curve, is an algebraic curve with an $r$-th root of its
(suitably twisted) canonical bundle. Forgetting the $r$-spin
structure reduces $\mgnrbar$ to $\mgnbar$. It is natural to ask if
Kontsevich's theorem admits a generalization to either of these
two cases.

The case of $\M_{g,n}(V)$ remains mysterious. It gives the
Gromov-Witten invariants of $V$ and their so-called gravitational
descendants, which assemble into a generating function whose
exponential is an analog of a $\tau$-function. In the case where
$V$ is a point, one recovers the $\tau$-function of the $\KdV$
hierarchy by Kontsevich's theorem. More generally, there is a
conjecture of Eguchi, Hori, and Xiong \cite{EHX} and of S.~Katz
which essentially states that this generating function is a
highest weight vector for a particular representation of the
Virasoro algebra. Presumably, there is some analog of an
integrable system which gives rise to this Virasoro algebra
action, should the conjecture hold.

On the other hand, the $\KdV$ hierarchy is just the first in a
series of integrable hierarchies $\KdV_r$, where   $r=2,3,\ldots$,
called the generalized $\KdV$, or Gelfand-Dickey hierarchies. In
the case of $r=2$, this is the usual KdV hierarchy. Each of these
hierarchies has a formal  solution, corresponding to the unique
$\tau$-function which satisfies an additional equation known as
the string (or puncture) equation. In \cite{W,W2}, Witten
formulated a generalization of his original conjecture, suggesting
that for each $r \ge 2$, there should exist moduli spaces  and
cohomology classes on  them  whose intersection  numbers assemble
into this $\tau$-function of the $\KdV_r$ hierarchy. The
corresponding moduli spaces of higher spin curves have recently
been constructed in~\cite{J,J2}. In this paper we present a
precise mathematical formulation of the generalized Witten
conjecture and prove it in several special cases including, in
particular, the case of genus zero.

Motivated by analogy with the construction of Gromov-Witten
invariants from the moduli space of stable maps, we introduce
axioms which must be satisfied by a cohomology class $\cv$ (called
the \emph{virtual class}) on the moduli space of $r$-spin curves
$\mgnrbar$ in order to obtain a cohomological field theory (\cft)
of rank $(r-1)$  in the sense of Kontsevich and Manin \cite{KM1}.
This virtual class on $\mgnrbar$ is an analog of the Gromov-Witten
classes of a variety $V$ (i.e., the pullbacks via the evaluation
maps of elements in  $H^\bullet(V)$). We realize this virtual
class in genus zero as the top Chern class of a tautological
bundle over $\M_{0,n}^{1/r}$  associated to the $r$-spin
structure. This yields a Frobenius manifold structure
\cite{Du,Hi,Ma} on the state space of the \cft\ which is
isomorphic to the Frobenius manifold associated to the versal
deformation of the $A_{r-1}$ singularity \cite{Du}. This is an
indication of the existence of a kind of ``mirror symmetry''
between the moduli space of $r$-spin curves and singularities.
According to Manin~\cite{Ma2}  ``isomorphisms of Frobenius
manifolds of different classes remain the most direct expression
of various mirror phenomena.'' Proving the generalized Witten
conjecture for all genera would provide further evidence of this
relationship.

As in the case of Gromov-Witten invariants, one can construct a
potential function from the integrals of the class $\cv$ on
different components of $\M_{g,n}^{1/r}$ to form the \emph{small
phase space} of the theory. The \emph{large phase space} is
constructed by introducing the tautological classes $\psi$
associated to canonical sections of the universal curve $\cgnr
\rTo \mgnrbar$, and can be regarded as a parameter space for a
family of \cfts. A \emph{very large phase space} (see \cite{KK,
FaPa, MaZ}), parametrizing an even larger space of \cfts, is
obtained by considering classes $\lambda$, associated to the Hodge
bundles,  and classes $\mu$, associated to the universal spin
structure.

     We show that the corresponding potential function satisfies
analogs of the puncture and dilaton equations and also a new
differential equation obtained from a universal relation
involving the class $\mu_1$.  These relations hold in all genera.
Topological recursion relations are also obtained from
presentations of these classes in terms of boundary classes  in
low genera.

Finally, using the new relation involving $\mu_1$, we show that
the genus zero part of the large phase  space potential
$\Phi_0(\bt)$ is completely determined by the geometry, and this
potential agrees with the generalized Witten conjecture in genus
zero.

Some of our constructions were foreshadowed by Witten, who
formulated his conjecture even before the relevant moduli spaces
and cohomology classes had been constructed, just as he had done
in the case of the topological sigma model and quantum cohomology.
We prove that his conjecture has a precise algebro-geometric
foundation, just as in the case of Gromov-Witten theory. Witten
also outlined a formal argument to justify his conjecture in genus
zero. Our work shows that the formulas that he ultimately obtained
for the large phase space potential function in genus zero are
indeed correct, provided that the geometric objects involved are
suitably interpreted. This is nontrivial even in genus zero
because the underlying moduli spaces are not schemes, but stacks.
We proceed further to prove relations between various tautological
classes associated to the $r$-spin structures and to derive
differential equations for the potential function associated to
them.

Notice also that one can introduce moduli spaces $\M_{g,n}^{1/r}(V)$ of
stable $r$-spin maps into a variety $V$,
where one combines the data of both the stable maps and the
$r$-spin structures. The analogous construction on these spaces
yields a Frobenius manifold which combines Gromov-Witten
invariants (and quantum cohomology) with the $\KdV_r$
hierarchies. Work in this direction is in progress \cite{JK}.

In the first section of this paper, we review the moduli space
$\mgnrbar$ of genus $g,$ stable $r$-spin curves, which was
introduced in \cite{J,J2}. We also discuss the stratification of
the boundary of $\mgnbar^{1/r}$. The boundary strata fall into two
distinct categories---the so-called Neveu-Schwarz and Ramond
types.

In the second section we introduce canonical morphisms,
tautological bundles, tautological cohomology classes, and
cohomology classes associated to the boundary strata of
$\mgnrbar$, and we derive a new relation involving the $\mu_1$
class.

In the third section, we define a cohomological field theory
(\cft) in the sense of    Kontsevich and Manin, its small phase
space potential function, and the associativity (WDVV) equation.
We then review the construction of Gromov-Witten invariants for
the moduli space of stable maps and define the large and very
large phase spaces in the Gromov-Witten  theory. Motivated by this
example, we explain how one may construct a \cft\ and the various
potential functions from analogous intersection numbers on
$\mgnrbar$, assuming that the virtual class $\cv$ exists.

In the fourth section we state axioms which $\cv$ must satisfy in
order to obtain a \cft.  We show that these axioms give a complete
\cft\ with a flat identity, and we construct the class $\cv$  in
genus zero, as well  as in the case $r=2$.

In the fifth section, we obtain analogs of the string and dilaton
equations for this $r$-spin \cft, and we find a new equation based
on the relation involving the $\mu_1$  class. We also prove the
analog of topological recursion relations in genus zero.

In the sixth section, we use the new relation for the class
$\mu_1$ to completely determine the genus-zero part of the large
phase space potential.

Finally, in the seventh section, we give a precise formulation of
the generalized Witten conjecture and prove that the genus-zero,
large phase space potential of the $r$-spin \cft\ yields a
solution to the semiclassical limit of the $\KdV_r$ hierarchy,
thereby proving the Witten conjecture in genus zero.  We conclude
with our own $W$-algebra conjecture, a $\KdV_r$-analog of a
refinement of the Virasoro conjecture \cite{EHX}.
\begin{ack}
  We would like to thank D.~Abramovich, M.~Adler, B.~Dubrovin,
  J.~Figueroa-O'Farrill,  E.~Getzler, A.~Kabanov, Yu.~Manin, P.~van~Moerbeke,
  A.~Polishchuk,  M.~Rosellen, J.~Stasheff, and Y.~Zhang for useful
  exchanges.  We   would also like to thank Heidi Jarvis for help with
  proofreading and typesetting this paper.   A.V.\ thanks the Max-Planck-Institut f\"ur
  Mathematik in Bonn for hospitality and financial support.
\end{ack}

\section{The Moduli Space of $r$-spin Curves}

In this section, we review the definition and some of the basic
properties of the moduli space $\mgnrbar$ of genus $g$,
$n$-pointed, stable $r$-spin curves.

\subsection{An overview of $\mgnrbar$}\label{overview}

As the definition of $\mgnrbar$ is rather involved, we motivate it
by starting with an intuitive approach to $r$-spin curves and
their moduli space.

A smooth $r$-spin curve is essentially just a curve with an $r$-th
root of the canonical bundle $\omega_X$ (suitably twisted). In
other words, it is a pair $(X,\cl)$ where $X$ is a smooth curve
and $\cl$ is a line bundle on $X$ such that $\cl^{\otimes r}$ is
isomorphic to the canonical bundle   $\omega_X$. Given a
collection of integers $\bm=(m_1, \dots, m_n)$, an $n$-pointed
smooth $r$-spin curve of type $\bm$ is a smooth $n$-pointed curve
$(X, p_1, \dots, p_n)$ with a line bundle $\cl$ on $X$, such that
$\cl^{\otimes r}$ is isomorphic to $\omega_X
(-\displaystyle\sum^n_{i=1} m_ip_i)$. For degree reasons such a
bundle exists only if $2g-2-\sum m_i$ is divisible by $r$.  When
this condition is met, there are $r^{2g}$ choices of $\cl$ on $X$.

If we want to compactify the space of smooth $r$-spin curves by
allowing the curve $X$ to degenerate to a stable curve, the above
definition of an  $r$-spin structure is insufficient. In
particular, there is often no line bundle $\cl$ such that
$\cl^{\otimes r}$ is isomorphic to $\omega_X (-\sum m_ip_i)$, even
when the degree condition is satisfied. One possible
solution---replacing line bundles by arbitrary rank-one,
torsion-free sheaves---permits too many potential candidates. The
correct structure required in this case amounts essentially to an
explicit choice of isomorphism (or homomorphism when $\cl$ is not
locally free)
        $b:\cl^{\otimes r} \rTo \omega_X (-\sum m_i p_i)$, with
some additional technical restrictions described in Definitions
\ref{def1} and \ref{def2}.

There are two very different types of behavior of this
torsion-free sheaf $\cl$ near a node $q\in X$. When it is still
locally free, the sheaf $\cl$ is said to be \emph{Ramond} at the
node $q$. If the sheaf $\cl$ is not locally free at $q$, it is
called \emph{Neveu-Schwarz}.

In the Ramond case, the homomorphism $b$ is still an isomorphism
(near the node $q$), but in the Neveu-Schwarz case it cannot be an
isomorphism. The local structure of the sheaf $\cl$ near a
Neveu-Schwarz node can be described as follows.

Near the node $q$, the curve $X$ has two coordinates $x$ and $y$,
such that $xy=0$; and the sheaf $\omega_X$ (or $\omega_X (-\sum
m_ip_i)$) is locally generated by $\frac{dx}{x} = -\frac{dy}{y}$.
Near $q$ the sheaf $\cl$ is generated by two elements $\ell_+$ and
$\ell_-$ supported on the $x$ and $y$ branches respectively (that
is, $x\ell_-=y\ell_+=0$).  The two generators may be chosen so
that the homomorphism $b:\cl^{\otimes r} \rTo \omega_X(-\sum m_i
p_i)$ takes $\ell^{\otimes r}_+$ to $x^{m_+
+1}(\frac{dx}{x})=x^{m_+}dx$, and so that $b$ takes $\ell^{\otimes
r}_{-}$ to $y^{m_-+1}(\frac{dy}{y})=y^{m_-}dy$, where
$(m_++1)+(m_-+1)=r$ is the order of vanishing of $b$ at the node
$q$.

   One more difficulty arises when $r$ is not prime---in this case
the moduli of stable curves with $r$-spin structure, as described
above, is not smooth.  The remedy is to include all $d$-spin
structures for every $d$ dividing $r$, satisfying some natural
compatibility conditions. This is described in Definition
\ref{def2}.

We now give the definition of $r$-spin curves.

\subsection{Higher spin curves}

\begin{df}\label{def:prestable}
A \emph{prestable curve} is a reduced, complete, algebraic curve with at
worst nodes as singularities.
\end{df}

\begin{df}\label{def1}
Let $(X, p_1, \dots, p_n)$ be a prestable, $n$-pointed, algebraic
curve,  $\ck$ be a rank-one, torsion-free sheaf on $X$, and $\bm =
(m_1,\ldots,m_n)$  be a collection of integers. A \emph{$d$-th
root of $\ck$ of type $\mathbf{m}$} is a pair $(\ce, b)$, where
$\ce$ is a rank-one, torsion-free sheaf, and $b$ is an
$\co_X$-module homomorphism $$
  b: \ce^{\tensor d}  \rTo \ck \otimes \co_X(-\sum m_ip_i)
$$
with the following properties:
\begin{itemize}
\item $d \cdot \deg \ce = \deg \ck-\sum m_i$
\item $b$ is an isomorphism on the locus of $X$ where $\ce$ is
locally free
\item for every point $p \in X$ where $\ce$ is not free, the
length of the cokernel of $b$ at $p$ is $d-1$.
\end{itemize}
\end{df}

The condition on the cokernel amounts essentially to the condition
that the order of vanishing of $b$ at a node should be $d$. For
any $d$-th root $(\ce,b)$ of type $\mathbf{m}$, and for any
$\mathbf{m}'$ congruent to $\mathbf{m} \pmod d$, we can construct
a unique $d$-th root $(\ce',b')$ of type $\mathbf{m}'$ simply by
taking $\ce'=\ce \otimes \co(1/d \sum (m_i-m_i')p_i)$.
Consequently, the moduli of curves with $d$-th roots of a bundle
$\ck$ of type $\mathbf{m}$ is canonically isomorphic to the moduli
of curves with $d$-th roots of type $\mathbf{m}'$. Therefore,
unless otherwise stated, we will always assume the type
$\mathbf{m}$ of a $d$-th root has the property that $0 \leq m_i<d$
for all $i$. Unfortunately, the moduli space of curves with $d$-th
roots of a fixed sheaf $\ck$ is not smooth when $d$ is not prime,
and so we must consider not just roots of a bundle, but rather
coherent nets of roots \cite{J}.  This additional structure
suffices to make  the moduli  space of curves with a coherent net
of roots smooth.

\begin{df}\label{def2}
Let $\ck$ be a rank-one, torsion-free sheaf on a prestable
$n$-pointed curve $(X, p_1, \ldots, p_n)$. A \emph{coherent net of
$r$-th roots of $\ck$ of type $\mathbf{m}=(m_1, \ldots, m_n)$} is
a pair $(\{\ce_d\}, \{c_{d, d'}\})$ of a set of sheaves and a set
of homomorphisms as follows.  The set of sheaves consists of a
rank-one, torsion-free sheaf $\ce_d$ on $X$ for every divisor $d$
of $r$; and the set of homomorphisms consists of an $\co_X$-module
homomorphism $$ c_{d,d'} : \ce^{\tensor d/d'}_{d} \rTo \ce_{d'} $$
for every pair of divisors $d',d$ of $r$,  such that $d'$ divides
$d$. These sheaves and homomorphisms must satisfy the following
conditions:
\begin{itemize}
\item $\ce_1=\ck$ and $c_{1,1}=\mathbf{1}$.
\item For each divisor $d$ of $r$ and each divisor $d'$ of $d$,
the  homomorphism $c_{d,d'}$  makes $(\ce_d, c_{d,d'})$ into a
$d/d'$-th root of $\ce_{d'}$ of type $\mathbf{m}'$, where
$\mathbf{m}'=(m'_1, \ldots, m_n')$ is the reduction of $\bm$
modulo $d/d'$ (i.e. $0\le m_i' < d/d'$ and $m_i \equiv m_i' \pmod
d/d'$).
\item The homomorphisms $\{c_{d,d'}\}$ are compatible.
That is, the diagram
$$
\begin{diagram}
(\ce^{\tensor d/d'}_{d})^{\tensor d'/d''} & \rTo^{(c_{d,d'})^{\tensor d'/d''}}
& \qquad\ce^{\tensor d'/d''}_{d'}\\
&  \rdTo^{c_{d,d''}}  & \dTo c_{d',d''} \\
& &  \ce_{d''}\\
\end{diagram}
$$
commutes for every $d''|d'|d|r$.
\end{itemize}
\end{df}

If $r$ is prime, then a coherent net of $r$-th roots is simply an
$r$-th root of $\ck$.  Even when $d$ is not prime, if the root
$\ce_d$ is locally free, then for every divisor $d'$ of $d$, the
sheaf $\ce_{d'}$ is uniquely determined, up to an automorphism of
$\ce_{d'}$. In particular, if $\bm'$ satisfies the conditions
$\bm' \equiv \bm \pmod {d'}$ and $0 \leq m'_i <d'$,  then the
sheaf $\ce_{d'}$ is isomorphic to $\ce_{d}^{\tensor d/d'}\otimes
         \co\left(\frac{1}{d'} \sum (m_i -m'_i)p_i\right)$.
\begin{df}
An \emph{$n$-pointed, $r$-spin curve of type $\bm = (m_1, \ldots,
m_n)$} is an $n$-pointed, prestable curve $(X, p_1, \ldots, p_n)$
with a coherent net of $r$-th roots of $\omega_X$ of  type
$\mathbf{m}$, where $\omega_X$ is the (canonical) dualizing sheaf
of $X$.  An $r$-spin curve is called \emph{smooth}  if $X$ is
smooth, and it is called \emph{stable} if $X$ is stable.
\end{df}

\begin{ex}
Smooth 2-spin curves of type $\mathbf{0}:=(0,0,\ldots,0)$
correspond to classical spin curves (a curve with a
theta-characteristic) $(X,\ce_2)$, \emph{with an explicit
isomorphism $\ce^{\tensor 2}_{2}  \irightarrow \omega$.}
\end{ex}

\begin{df}
An isomorphism of $r$-spin curves $$ (X, p_1, \ldots, p_n,   (\{
\ce_d\}, \{c_{d, d'}\})) \irightarrow (X', p'_1, \ldots, p'_n,
(\{\ce'_d\}, \{c'_{d,  d'}\})) $$ of the same type $\bm$ is an
isomorphism of pointed curves $$ \tau : (X, p_1,  \ldots, p_n)
\irightarrow (X', p'_1, \ldots, p'_n) $$ and a set of sheaf
isomorphisms $\{\beta_d : \tau^* \ce'_d \irightarrow \ce_d\},$
with $\beta_1$  being the canonical isomorphism $\tau^*
\omega_{X'}(-\sum_i m_i  {p'}_i)  \irightarrow \omega_X(-\sum
m_ip_i),$ and such that the homomorphisms $\beta_d$ are compatible
with all the maps $c_{d,d'}$ and $\tau^*c'_{d,d'}$.
\end{df}

Every $r$-spin structure on a smooth curve $X$ is determined, up
to isomorphism, by a choice of a line bundle $\ce_r$, such that
$\ce^{\tensor r}_r \cong \omega_X (-\sum m_i p_i)$.  In
particular, if $X$ has no automorphisms, then the set of
isomorphism classes of $r$-spin structures (if non-empty) of type
$\mathbf{m}$ on $X$ is a principal homogeneous space for the group
of $r$-torsion points of the Jacobian of $X$. Thus there are
$r^{2g}$ such isomorphism classes.

\begin{ex}\label{ex:modular}
If $g=1$ and $\mathbf{m}=\mathbf{0}$, then $\omega_X$ is
isomorphic to $ \co_X$, and a smooth $r$-spin curve is just an
elliptic curve $X$ with a line bundle $\ce_r$ corresponding to an
$r$-torsion point of $X$, together with an explicit isomorphism
$\ce^{\tensor r}_{r} \irightarrow \co_X$. In particular, the stack
of stable, one-pointed $r$-spin curves of genus one and type
$\mathbf{0}$ forms a gerbe over the disjoint union of modular
curves $\coprod_{d|r} X_1(d)$.
\end{ex}

\begin{df}
The stack of connected, stable, $n$-pointed, $r$-spin curves of
genus $g$ and type $\bm=(m_1,\ldots,m_n)$ is denoted by
$\mgnrmbar$.  The disjoint union
$\displaystyle\coprod_{\substack{\mathbf{m} \\ 0 \leq m_i <r}}
\mgnrmbar$ is denoted by $\mgnrbar$.
\end{df}

\begin{rem} \label{rem:restrict}
As mentioned above, no information is lost by restricting $\bm$ to
the range $0\le m_i \le r-1$, since when $\bm \equiv \bm' \pmod r$
every $r$-spin  curve of type $\bm$ naturally gives an $r$-spin
curve of type   $\bm'$  simply by $$ \ce_d \mapsto \ce_d \tensor
\co (\sum \frac{m_i-{m'}_i}{d} p_i). $$ Thus $\mgnrmbar$ is
canonically isomorphic to
$\overline{\scheme{M}}^{1/r,\mathbf{m}'}_{g,n}$.
\end{rem}

\subsection{Basic properties of the moduli
space}\label{onepointtwo}

\

In \cite{J} it was shown that $\mgnrbar$ is a smooth
Deligne-Mumford stack, finite over $\mgnbar$, with a projective,
coarse moduli space.  For $g>1$ the spaces $\mgnrmbar$ are
irreducible if $\gcd(r, m_1, \ldots, m_n)$ is odd, and they are
the disjoint union of two irreducible components if $\gcd(r, m_1,
\ldots, m_n)$ is even. When $r=2$ (and in fact, for all even $r$)
this is due to the well-known fact that even and odd theta
characteristics on a curve cannot be deformed into one another
\cite{mumford:theta-chars}. These two components will be denoted
$\mgnrmbareven$ and $\mgnrmbarodd$ respectively.

When the genus $g$ is zero  the moduli space
$\mbar^{1/r,\bm}_{0,n}$ is either empty (if $r$ does not divide $2
+ \sum m_i$), or is canonically isomorphic to $\mbar_{0,n}$. Note,
however, that this isomorphism is not an isomorphism of stacks,
since the automorphisms of elements of $\mgnrmbar$ vary
differently from the way that automorphisms of the underlying
curves vary.  We will discuss this further in Section~\ref{auto}.
In any case, $\mbar^{1/r,\bm}_{0,n}$ is always irreducible.

When the genus $g$ is one, the space
$\mbar^{1/r,\mathbf{m}}_{1,n}$ is the disjoint union of $d$
irreducible components, where $d$  is the number of divisors of
$\gcd(r,m_1, \ldots, m_n)$.  We will denote the irreducible (and
connected) component indexed by a divisor $e$ of $\gcd(r,m_1,
\ldots, m_n)$ by $\mbar^{1/r,\bm,(e)}_{1,n}$. When $\bm$ is zero,
as mentioned in Example~\ref{ex:modular}, the locus of smooth
$r$-spin curves in this component consists of $n$-pointed,
elliptic curves with a torsion point of exact order $e$.

Throughout this paper we will denote the forgetful morphism by
$p:\mgnrbar \rTo \mgnbar $, and the universal curve  by $\pi:\cgnr
\rTo \mgnrbar$.  As in the case of the moduli space of stable
curves, the universal curve possesses canonical sections
$\sigma_i:\mgnrbar \rTo \cgnr$ for $i\,=\,1,\ldots,n$.  Unlike the
case of stable curves, however, the universal curve
$\scheme{C}^{1/r, \bm}_{g,n} \rTo \mgnrmbar$ is not obtained by
considering $(n+1)$-pointed $r$-spin curves. The curve
$\scheme{C}^{1/r,   \bm}_{g,n}$ is birationally equivalent to
$\overline{\scheme{M}}^{1/r, (m_1, m_2, \dots, m_n, 0)}_{g,n+1}$,
but they are not isomorphic.

There is one other canonical morphism associated to these spaces;
namely, when $d$ divides $r$, the morphism $$ [r/d]:\mgnrmbar \rTo
\mbar^{1/d,\mathbf{m}'}_{g,n} \mathrm{\ and\ } [r/d]:\mgnrbar \rTo
\mbar^{1/d}_{g,n}, $$ which forgets all of the roots and
homomorphisms in the net of $r$-th roots except those associated
to divisors of $d$.  Here $\mathbf{m}'$ is congruent to
$\mathbf{m} \pmod d$ and $0 \leq m'_i <d$ for all $i \in
\{1,\ldots, n\}$. In the case that the underlying curve is smooth,
this is equivalent to replacing the line bundle $\ce_r$ by its
$r/d$-th tensor power (and then taking the tensor product with
$\co(1/d \sum(m_i-m'_i)p_i )$).

The two components $\M^{1/r,\bm,even}_{g,n}$ and $\M^{1/r, \bm,
odd}_{g,n}$ that arise in the case that $\gcd(r,m_1, \ldots, m_n)$
is even are just the preimages of the spaces of even and odd
theta-characteristics in $\mbar^{1/2,\mathbf{0}}_{g,n}$ under the
map $[r/2]:\mbar^{1/r,\mathbf{m}}_{g,n} \rTo
\mbar^{1/2,\mathbf{0}}_{g,n}$.

\subsection{Boundary behavior}\label{boundary}

\subsubsection{Neveu-Schwarz and Ramond nodes.}\label{sec:NSRnodes}

\

At any node $q$ of a prestable curve $X$, there  are two types of
local behavior of an $r$-th root $(\ce_r, b_r)$ of $\omega_X
(-\sum m_ip_i)$. It is either locally free at $q$, in which case
the homomorphism $b_r$ is an isomorphism near $q$, or it is
torsion-free, but not locally free at $q$. In the locally free
case we will say that the root $\ce_r$ is \emph{Ramond} at $q$,
and in the non-locally free case it will be called
\emph{Neveu-Schwarz}.

If the $r$-th root sheaf $\ce_r$ of an $r$-spin structure on $X$
is Ramond at every node of $X$, then the whole net of roots is
completely determined (up to isomorphism) by the root $(\ce_r,
c_{r,1})$  as follows: $$\ce_d = \ce^{\tensor r/d}_r
           \otimes \co\left(\frac{1}{d} \sum (m_i -m'_i)p_i\right),$$
and $$c_{d,1} = c_{r,1} \otimes I,$$ where $I$ is the identity
homomorphism $$\co (\frac{1}{d}\sum (m_i-m'_i)p_i) \rTo \co
(\frac{1}{d} \sum (m_i-m'_i)p_i). $$

\begin{rem}\label{rem:NS}
The Neveu-Schwarz case is more complicated, in some sense, than
the Ramond case, but this is compensated for by the fact that the
cohomology classes defined by boundary strata with Neveu-Schwarz
nodes factor in a nice combinatorial way. Moreover, the hope of
constructing a cohomological field theory from $\mgnrbar$ is based
on the expectation that one can construct a canonical cohomology
class which vanishes on the strata where the $r$-spin structure
has an $r$-th root sheaf which is Ramond at some node. (This would
follow from Axiom~4 in Section~\ref{axioms}.)
\end{rem}

\subsubsection{Local structure at Neveu-Schwarz nodes.}\label{sec:NSnode}

\

Recall from Section~\ref{overview} (see also \cite{J,J2,W}) that
near a Neveu-Schwarz node $q$, an $r$-th root $(\ce_r, b_r)$ of
$\omega_X(-\sum m_ip_i)$ is uniquely determined by an $r$-th root
$(\tilde{\ce_r}, \tilde{b_r})$ of the bundle $\omega_{\tilde{X}}
(-\sum  m_ip_i - m^+ q^+-m^-q^-)$ on the normalization $\nu:
\tilde{X} \rTo X $ of $X$ at the node $q$.  Here $q^+$ and $q^-$
are the inverse images of $q$ under $\nu$, and $m^+$ and $m^-$ are
non-negative integers\footnote{Note that the integers $u$ and $v$
of the papers \cite {J} and \cite{J2} are $m^+ + 1$ and $m^-+1$,
respectively.} which sum to $r-2$. If $x$ and $y$ are local
parameters of $X$ near the node $q$ satisfying the equation $x
y=0$, then the sheaf $\ce$ is generated locally by the sections
$(x^{m_+}dx)^{1/r}$ and $(y^{m_-}dy)^{1/r}$; and $\tilde{\ce}_r$
is generated by $(x^{m_+}dx)^{1/r}$ on the $x$ branch of
$\tilde{X}$, and it is generated by $(y^{m_-}dy)^{1/r}$ on the $y$
branch of $\tilde{X}$. The points $q^+$ and $q^-$ are given by
$\{x=0\}$ and $\{y=0\}$, respectively, on $\tilde{X}$. The sheaf
$\tilde{\ce_r}$ is simply $\nu^* \ce_r$ modulo torsion; and $\nu_*
\tilde{\ce_r}$ is $\ce_r$, with $b_r$ induced from $\tilde{b_r}$
by adjointness.  We will call the integers $m^+$ and $m^-$ the
\emph{order of the $r$-spin structure} at the node, along the $x$
or $y$ branch, respectively.\footnote{The order $m^+$ and $m^-$ of
the $r$-spin structure along the $x$ or $y$ branch of a node is
not to be confused with the order of vanishing of the structure
maps. Indeed, if the $r$-th root bundle $\ce_r$ is Neveu-Schwarz
at a node, the order of vanishing of the map $c_{1,r}$ at that
node is exactly $(m^+ + 1)+(m^-+1)  = r$.}

In the case that $m^++1$ and $m^- +1$ are relatively prime,   one
can show (see \cite{J}) that $\tilde{\ce}$ and $\tilde{b}$
uniquely determine the entire net. However,  if $\gcd(m^++1,
m^-+1) =d$, then although $\ce_r$ still completely determines the
Neveu-Schwarz roots, $d$ divides $r$, and the $d$-th root $(\ce_d,
c_{d,1})$ of the net is locally free (Ramond), as are all roots
$(\ce_{d'}, c_{d',1})$ for every $d'$ dividing $d$. In particular,
although generators $(x^{m_+}dx)^{1/r}$ and $(y^{m_-}dy)^{1/r}$ in
$\ce_r$ determine $(x^{m_+}dx)^{1/d}$ and $(y^{m_-}dy)^{1/d}$, we
must identify $(\frac{dx}{x})^{1/d}$ with $(-\frac{dy}{y})^{1/d}$.
However,  this identification is only determined by
$(\ce_r,c_{r,1})$ up to a non-canonical choice of a $d$-th root of
unity. If the normalization $\tilde{X}$ at $q$ has two connected
components, then the $d$-th root $(\ce_d, c_{d,1})$ is determined
up to (non-canonical) isomorphism by $(\ce_r, b_r)$, but if
$\tilde{X}$ is connected, then $(\ce_d, c_{d,1})$ is not
determined by $(\ce_r, b_r)$, since an additional choice of a
$d$-th root of unity is required to construct $\ce_d$ from
$\ce^{\otimes r/d}_r$ (see Section~\ref{glue}).

\subsubsection{More detailed study of the Ramond case}

 \label{ramond}

\

Let $(\ce,b)$  be an $r$-th  root of $\omega_X(-\sum m_ip_i)$
which is Ramond at a node $q$ of $X$. The restriction of $\ce$ to
$q$ gives an exact sequence $$ 0 \rTo \mathfrak{m}_q \otimes \ce
\rTo \ce \rTo \ce|_q \rTo 0, $$ where $\mathfrak{m}_q$ is the
maximal ideal of the point $q$.  The sheaf $\ce':= \ce \otimes
\mathfrak{m}_q$ is a rank-one, torsion-free sheaf of degree
$(2g-2-\sum m_i)/r-1$ on $X$, and pulling $\ce'$ back to the
normalization  $\nu: \tilde{X} \rTo X$ of $X$ at $q$ gives, modulo
torsion, a rank-one, torsion-free sheaf $\ce'' := \nu^*
\ce'/\text{torsion}$, such that $\nu_* \ce''$ is equal to $\ce'$.

If $x$ and $y$ are local coordinates on $\tilde{X}$ near $q^+$ and
$q^-$ respectively, then $\ce''$ is locally generated by
$x(\frac{dx}{x})^{1/r}$ (respectively $y (\frac{dy}{y})^{1/r}$).
Therefore, the homomorphism $$ b'': {\mbox{$\ce''$}}^{\otimes r}
\rTo \nu^* \omega_X (-\sum m_i p_i) = \omega_{\tilde{X}} (-\sum
m_i p_i + q^+ + q^-), $$ induced by $b: \ce^{\tensor r} \rTo
\omega_{\tilde{X}}( -\sum m_i p_i)$, factors through $$
\omega_{\tilde{X}}( -\sum m_i p_i - (r-1)q^+ -(r-1)q^-) \rTo
\omega_{\tilde{X}}( -\sum m_i p_i +q^+ +q^-). $$

Thus the $r$-th root $(\ce_r, c_{r,1})$ can be Ramond at the node
if and only if $m^+ =r-1$ and $m^-=r-1$ satisfy the degree
conditions $$ \deg_{X^{(j)}} \omega_{\tilde{X}} -\sum m_i-m^+ -m^-
\equiv 0 \pmod r $$ on every connected component $X^{(j)}$ of
$\tilde{X}$.  In the case that $(\ce_r, c_{r,1})$ is Ramond, we
will define the order of the $r$-spin structure at the node to be
$m^+=m^-=r-1$ along both branches of the underlying curve.

Similarly, $$(\nu^* \ce)^{\otimes r} \rTo^{\nu^*b} \nu^*
\omega_X(-\sum m_ip_i) = \omega_{\tilde{X}}(-\sum m_ip_i+q^++q^-)
$$ corresponds to the choice $m^+=m^-=-1$.  In this special case
there is a residue map that canonically identifies $\ce_r|_{q^+}$
and $\ce_r|_{q^-}$ with $\mathbb{C}$.
\begin{prop}\label{residue}
If $(\ce_r,b)$ is an $r$-th  root of $\omega_X(-\sum m_ip_i)$ with
$m_i= -1$ for some $i$, then there is an isomorphism
\begin{equation}
  \label{eq:residue}
R_{p_i}:\ce_r|_{p_i} \irightarrow \mathbb{C}
\end{equation}
which is canonical up to a choice of an $r$-th root of unity.
\end{prop}

An immediate consequence is the following corollary:
\begin{crl}\label{trivial}
If $\sigma_i$ is the $i$-th  section of the universal curve $$
\pi:\cgnrm \rTo \mgnrmbar $$ with $\mathbf{m}=(m_1,\ldots, m_n)$
and $m_i=-1$ for some $i$, then the pullback $\sigma^*_i(\ce_r)$
of the universal $r$-th root $\ce_r$ is a trivial line bundle,
isomorphic (canonically up to a choice of an $r$-th root of unity)
to $\co_{\mgnrmbar}$.
\end{crl}

\begin{proof}{(of the proposition)}
Let $z$ be a local parameter on $X$ near $p$, so that the sheaf
$\ce_r$ is locally generated by an element $(\frac{dz}{z})^{1/r}$
which is well defined up to an $r$-th root of unity. Define the
map $R_p :\ce_r|_p \irightarrow \mathbb{C}$ by
$R_p((a_0+a_1z+\ldots)(\frac{dz}{z})^{1/r}):=a_0$.  It is easy to
check that this definition is independent of local parameter, and
hence defines a canonical isomorphism.
\end{proof}

\subsection{Graphs}

Much of the information about the structure of the boundary of
$\mgnrmbar$ can be encoded  in terms of decorated graphs.

Recall that the \emph{(dual) graph} of an $n$-pointed prestable
curve $(X, p_1,\ldots,p_n)$ consists of the following elements:
\begin{itemize}
\item \emph{Vertices},
corresponding to the irreducible components of $X$: a vertex $v$
is labeled with a non-negative integer $g(v)$, the (geometric)
genus of the component;
\item \emph{Edges},
 corresponding to the nodes of the curve:  an
edge connects two vertices (possibly even the same vertex, in
which case the edge is called a {\em loop}) if and only if the
corresponding node lies on the associated irreducible components;
\item \emph{Tails}, corresponding to the marked points $p_i \in X$,
                             $i=1,\ldots, n$:
a tail labeled by the integer $i$ is attached at the vertex
associated to the component of $X$ that contains $p_i$.
\end{itemize}

\begin{df}
A \emph{half-edge} of a graph $\Gamma$ is either a tail or one of the
two ends of a ``real'' edge of $\Gamma$.
We denote by $V(\Gamma)$ the set of vertices of  $\Gamma$ and by
$n(v)$ the number of half-edges of $\Gamma$ at the vertex $v$.
\end{df}

The following definition describes a class of graphs that are dual graphs
of stable pointed curves.

\begin{df}

Let $\Gamma$ be a graph. The number $$ g(\Gamma) =
\mathrm{dim}H^1(\Gamma) + \sum_{v\in V(\Gamma)}g(v) $$ is called
the \emph{genus} of a graph $\Gamma$.

A graph $\Gamma$ (not necessarily connected) is called
\emph{stable} if $2g(v)-2 + n(v) > 0$ for every  $v\in V(\Gamma)$
(in particular, it satisfies $2g(\Gamma)-2 +n > 0$, where $n$ is
the number of tails of $\Gamma$).
\end{df}

To describe strata of the moduli space of $r$-spin curves,  we
decorate the graphs with additional data coming from the $r$-spin
structure. In particular, the type $\mathbf{m} = (m_1,\ldots ,
m_n)$ gives a marking to each of the tails.
\begin{df}
Fix an integer $r\geq 2$. A \emph{decorated stable graph} is a
stable graph                  with a marking  of each half-edge by
a non-negative integer $m<r$, such that for each edge $e$ the
marks $m^+$ and $m^-$ of the two half-edges of $e$ satisfy $$ m^+
+ m^- \equiv r-2 \pmod r. $$
\end{df}
As   mentioned in Section~\ref{boundary}, decorated stable graphs
with $n$ tails and genus $g$ correspond to boundary strata in
$\mgnbar^{1/r}$.
\medskip

\begin{df}
Given a stable $r$-spin curve $\mathfrak{X}$ of type $ \bm=(m_1,
\dots, m_n)$, the \emph{decorated dual graph} of $\mathfrak{X}$ is
the dual graph $\Gamma$ of the underlying curve $X$, with the
following additional markings.  The $i$-th tail is marked by
$m_i$, and each half-edge associated to a node of $X$ is marked by
the order ($m^+$ or $m^-$) of the $r$-spin structure along the
branch of the node associated to that half-edge.
\end{df}

\begin{df}
Let $\Gamma$ be a connected stable graph (or a decorated stable
graph) with $n$ tails and of genus $g$ . We denote by
$\mbar_\Gamma$ (or by $\mrbar_\Gamma$) the closure in $\mgnbar$
(or in $\mgnbar^{1/r}$) of the moduli space of stable curves (or
$r$-spin curves) whose dual graph is $\Gamma$. If
$\Gamma=\displaystyle{\coprod_{i\in I}}\Gamma_i$ is the disjoint
union of connected subgraphs $\Gamma_i$ then we denote by
$\M_{\Gamma}$ the product $\displaystyle{\prod_{i\in
I}}\,\M_{\Gamma_i}$, and similarly $\M_{\Gamma}^{1/r} =
\displaystyle{\prod_{i\in I}}\,\M^{1/r}_{\Gamma_i}$.
\end{df}

\subsection{Automorphisms of $r$-Spin Curves}\label{auto}

As mentioned in Section \ref{onepointtwo}, even in the genus zero
case, where there is a unique $r$-spin structure of a given type
$\bm=(m_1, \dots, m_n)$ for each genus-zero curve (provided  $\bm$
satisfies the degree requirement $\sum m_i \equiv 2 \pmod r$), the
automorphisms of the $r$-spin structure ensure that
$\M^{1/r,\bm}_{0,n}$ is not isomorphic, as a stack, to $\M_{0,n}$.
The automorphisms of $r$-spin structures will play an important
role later in this paper, particularly in the determination of the
degrees of morphisms and the properties of various cohomology
classes under restriction and pullback. Consequently, we need to
understand the group of automorphisms of an $r$-spin curve.

First, we introduce some notation. Let $\mathfrak{X}=(X,
p_1,\dots,p_n, (\{\ce_d\},\{ c_{d,d'}\}))$ be an $r$-spin curve,
and let $\Gamma$ be its decorated dual graph.  Let $V$ be the set
of vertices of $\Gamma$ and $E_{nl}$ be the set of edges
\emph{which do not start and end at the same vertex} (i.e.,
non-loops). Furthermore, for each $v \in V$, denote by $X_v$ the
irreducible component of $X$ associated to vertex $v$; and denote
by $F_v$ the set of all half-edges attached to $v$ in $\Gamma$.
For each $f$ in $F_v$ let $p_f$ be the point of $X_v$ associated
to $f$, and let $m_f$ be the marking of $f$. Finally, for each
vertex $v$, if $i_v: X_v \rTo X$ denotes the obvious inclusion,
then the $r$-spin structure $(\{\ce_d\},\{ c_{d,d'}\})$ pulls back
to a collection $(\{ i_v^* \ce_d\}, \{ i_v^* c^v_{d,d'}\}))$ of
sheaves and morphisms on $X_v$. However, these sheaves are not
necessarily torsion free. Taking the quotient of each sheaf $i_v^*
\ce_d$ by its torsion submodule gives a rank-one, torsion-free
sheaf, which we denote by $\ce^v_d$.  It is easy to see that the
homomorphisms $\{ c^v_{d,d'}\}$, induced on the $\{\ce^v_d\}$ from
the homomorphisms $i_v^* \ce_d$, make $\mathfrak{X}_v = (X_v,
p_{f^v_1}, \dots, p_{f^v_{k_v}}, (\{\ce^v_d \}, \{ c^v_{d,d'}\}))$
into an $r$-spin curve.  We will call $\mathfrak{X}_v$ the
\emph{restriction of $\mathfrak X$ to the curve $X_v$}. Any $e \in
E_{nl}$ consists of two half-edges $f^+_e$ and $f^-_e$; and we
denote by $d_e$ the integer $$d_e := \gcd (m_{f^+_e}+1,
m_{f^-_e}+1)=\gcd(m_{f^+_e}+1, r).$$

\begin{prop}
If the underlying pointed curves $(X_v, p_{f^v_1}, \dots,
p_{f^v_{k_v}})$ of the $r$-spin curves $\mathfrak{X}_v$ have no
non-trivial automorphisms (this is true for a generic curve with
$g+n>2$), then
\begin{enumerate}
\item for each $v \in V$ the automorphism group $\Aut(\mathfrak{X}_v)$ of
$\mathfrak{X}_v$ is isomorphic to $\mathbf{\mu}_r$, the group of
$r$-th roots of unity; and
\item given any orientation of the edges of the dual graph $\Gamma$ of
$\mathfrak{X}$, the automorphism group $\Aut(\mathfrak{X})$ of
$\mathfrak{X}$ fits into the following exact sequence
\begin{equation}\label{eq:auto}
1 \rTo \Aut(\mathfrak{X}) \rTo \prod_{v \in V} \Aut
(\mathfrak{X}_v) \rTo^{\partial} \prod_{e \in
E_{nl}}\mathbf{\mu}_{r/d_e}.
\end{equation}
\end{enumerate}
Here the map $\partial$ is defined as follows.  The orientation of
each edge $e$ determines the choice of half-edges $f^+_e$ and
$f^-_e$ and of corresponding vertices $v^+_e$ and $v^-_e$ (the
unique vertex of $\Gamma$ attached to $f^+_e$ or $f^-_e$,
respectively).  The map $\partial$ maps the element $\prod \zeta_v
\in \displaystyle \prod_{v \in V} \mathbf{\mu}_r$ to the element
$\prod(\zeta^{d_e}_{v^+_e} \zeta^{-d_e}_{v^-_e})$.  Note that,
although the map $\partial$ depends upon the given orientation of
the edges of $\Gamma$, the kernel of $\partial$ is independent of
orientation.
\end{prop}

\begin{proof}
If $(X_v, p_{f^v_1}, \dots, p_{f^v_e})$ has no automorphisms, then
each term $\ce^v_d$ in the $r$-spin structure $(\{\ce_d^v\},\{
c_{d,d'}^v\})$ is locally free (since $X_v$ is smooth) of rank one
and has automorphism group $\mathbb{C}^*=H^0(X_v,\co^*_{X_v})$,
which acts on $\ce_d$ by multiplication.  However,  an
automorphism of the $r$-spin structure must also be compatible
with the structure maps $\{c_{d,d'}\}$.  In particular,
compatibility with the isomorphism $c_{r,1}:(\ce_r^v)^{\otimes r}
\irightarrow \omega_{X_v}(-\sum_{f \in F_v} m_f p_f)$ shows that
any automorphism $\sigma_r$ of $\ce_r$ must satisfy $(\sigma_r)^r
= 1$. Moreover, compatibility with $c_{d,d'}$ shows   that
$\sigma_{d'}=(\sigma_d)^{d/d'}$ for every $d'$ dividing $d$ and
$d$ dividing $r$. Thus $\sigma_r$ corresponds to some $ \zeta \in
\mathbf{\mu}_r$, and for every $d$ dividing $r$, the automorphism
$\sigma_d$ is just $ \zeta^{r/d}$. This proves the first part of
the proposition.

For the second part, it is easy to see that any automorphism of
the whole $r$-spin curve $\mathfrak{X}$ induces, by restriction,
an automorphism of $\mathfrak{X}_v$ for each $\mathfrak{X}_v$, and
the map $$\Aut (\mathfrak{X}) \rTo \prod_{v \in V} \Aut
(\mathfrak{X}_v)$$ is injective.

Moreover, for any edge $e \in E_{nl}$ corresponding to half edges
$f^+$ attached to vertex $v^+$ and $f^-$ attached to vertex $v^-$,
an automorphism $\sigma$ of $\mathfrak{X}$ will induce
automorphisms $\sigma_{v^+} = \zeta_+ \in \mathbf{\mu}_r$ and
$\sigma_{v^-}= \zeta_- \in \mathbf{\mu}_r$; and these
automorphisms must agree whenever $\ce_r$ is Ramond (locally free)
at the node $p_e$ corresponding to edge $e$.  Similarly, if
$\ce_d$ is Ramond at $p_e$, then $(\zeta_+)^{r/d}$ must equal
$(\zeta_-)^{r/d}$.  However,  if the sheaf $\ce_d$ is
Neveu-Schwarz at $p_e$, then $\sigma_{v^+}$ and $\sigma_{v^-}$ act
on distinct vector spaces (the sheaf $\ce^{v^+}_d|_{p_{f^+_e}}$ is
not the same as $\ce^{v_-}_d|_{p_{f^-_e}}$), and so $\ce_d$
imposes no compatibility condition on $\sigma_{v^+}$ and
$\sigma_{v^-}$.

Since $\ce_d$ is Ramond at the node $p_e$ precisely when $d$
divides both $m_{f^+}+1$ and $m_{f^-}+1$, we have that the
condition imposed by compatibility for $\sigma_{v^+}$ and
$\sigma_{v^-}$ is precisely the equality $(\sigma_{v^+})^{r/d}
=(\sigma_{v^-})^{r/d}$ for $$d=d_e := \gcd(m_{f^+}+1,
m_{f^-}+1)=\gcd(m_{f^+}+1, r).$$ It is clear that any $\prod
\sigma_v \in \displaystyle{\prod} \Aut (\mathfrak{X}_v)$ which
meets this compatibility condition at each edge $e \in E_{nl}$
defines a global automorphism $\sigma \in \Aut (\mathfrak{X})$.
\end{proof}

Of course, since we only care about the kernel of $\partial$, the
right-most term in  exact sequence~(\ref{eq:auto}) might as well
include all the edges, including loops.  This is because for any
loop $e$, the map $\partial$ will always map every element of
$\displaystyle \prod_{v \in V} \Aut(\mathfrak{X}_v)$ to $1 \in
\mathbf{\mu}_{r/d_e}$, since the vertices $v^+$ and $v^-$ (and
hence also $\sigma_{v^+}$ and $\sigma_{v^-}$) are the same.

For $r$-spin curves of genus $1$ and $2$ with a generic involution
(i.e., in $\M^{1/r}_{1,1}$ and $\M^{1/r}_{2,0}$), we have the
following description of the group of automorphisms.  Note that
when $g=n=1$, degree requirements force $m_1$ to be zero.  Also
when $g=2$ and $n=0$, degree requirements force $r$ to be $2$.

\begin{prop}
\

\begin{enumerate}
\item If the underlying curve $(X,p)$ of a smooth $r$-spin curve
$\mathfrak{X} \in \mbar^{1/r,\mathbf{0},(j)}_{1,1}$ of index $j$
has  no automorphisms other than the elliptic involution, then the
automorphism group of $\mathfrak{X}$ is $$\Aut \mathfrak{X}=
\begin{cases} \mu_r \cross \mathbb{Z}/2 & \text{if }  j =1 \text{
or } j=2 \\ \mu_r & \text{if }  j>2.
\end{cases}$$
\item If the underlying curve $X$ of a smooth $2$-spin curve
$\mathfrak{X} \in \MM^{1/2}_{2,0}$ has no automorphisms other than
the hyperelliptic involution, then the automorphism group of
$\mathfrak{X}$ is $$\Aut \mathfrak{X} = \mu_r \cross
\mathbb{Z}/2.$$
\end{enumerate}
\end{prop}

\begin{proof}
In the case of $g=n=1$, smooth $r$-spin curves of type $0$ and
index $j$ correspond to the torsion points of $X$ of exact order
$j$. It is well known that the involution $i: X \rTo X$ acts
without fixed points on the points of exact order $j$, unless $j$
is $1$ or $2$, in which case the involution fixes all $2$-torsion
points (including the identity, corresponding to the trivial
bundle $\co_X=\omega_X$).

It is easy to check (e.g., by explicitly writing out the
coordinates) that for sheaves $\ce_r$ of index $2$ (or $1$),
corresponding to $2$-torsion, there is a canonical choice of
isomorphism $\tau: i^*\ce_r \rTo \ce_r$ such that
\begin{itemize}
\item $i^* \tau \circ \tau$ is trivial,
\item any other isomorphism $i^* \ce_r \rTo \ce_r$ differs
from $\tau$ by an element $\zeta \in \mu_r$, and
\item  $\tau$ commutes with all elements of $\mu_r$.
\end{itemize}
Thus $\Aut \mathfrak{X}$ has order $2r$ and is Abelian; and if $r$
is even, then every automorphism has order dividing $r$. Thus the
     proposition follows in genus $1$.

The proof in genus $2$ is similar, but simpler, since every
$2$-spin structure is fixed by the involution.
\end{proof}

\subsection{Gluing}
\label{glue}

\

In the case of moduli spaces of ordinary stable curves, if a graph
$\Gamma$ is obtained from another graph $\tilde{\Gamma}$ (not
necessarily connected) by gluing together two tails of
$\tilde{\Gamma}$ (thus producing a new edge), there is a natural
\emph{gluing morphism}
\begin{equation}
\label{eq:gluing}
  \rho: \mbar_{\tilde{\Gamma}} \rTo \mbar_{\Gamma} \rInto \mgnbar.
\end{equation}

It corresponds to gluing together the punctures on a curve $X \in
\mbar_{\tilde{\Gamma}}$ associated to the two tails, and thus the
curve $\rho(X) \in \mbar_{\Gamma}$ will have an additional node
and two  fewer punctures than $X$. A similar gluing operation
sometimes exists for $r$-spin curves, but even then we often need
to include extra data.

\subsubsection{Two irreducible components}
\label{twocomp} Consider the case of an $r$-spin curve with a
single Neveu-Schwarz node.  Assume the normalization of the
underlying curve $X$ at the node has two connected components:
$X^{+}$ of genus $k$, and $X^{-}$ of genus $g-k$. The decorated
dual graph of the $r$-spin curve looks like this $$ \Gamma
=\tcongmpm. $$ The $r$-th root bundle $\ce_r$ factors as $\ce_r =
\ce^{+}_r \oplus \ce^{-}_r$, with $\ce^{+}_r$ an $r$-th root of
$\omega_{X^{+}} (-\displaystyle\sum^j_{l=1} m_{i_l}
p_{i_l}-m^+q^+)$  on $X^{+}$ and $\ce^{-}$ an $r$-th root of
$\omega_{X^{-}} (-\displaystyle\sum^n_{l=j+1} m_{i_l} p_{i_l} -
m^-q^-)$ on $X^{-}$. Here $m^+$ and $m^-$ can easily be calculated
(for degree reasons) as the unique non-negative integers summing
to $r-2$ such that the $r$-th roots of the corresponding (twisted)
canonical bundles exist; namely,
\begin{eqnarray}
\label{eq:mplus} m^+ &\equiv& 2k-2-\sum^j_{l=1} m_{i_l} \pmod r
\qquad \mathrm{ and } \\ \label{eq:mminus}
  m^- &\equiv& 2(g-k)-2 -\sum^n_{l=j+1} m_{i_l} \pmod r.
\end{eqnarray}
Of course, since the degree of the original bundle was divisible
by $r$, each of the two relations implies the other.  The
Neveu-Schwarz case occurs exactly when the solutions $m^+$ and
$m^-$ to the congruences~(\ref{eq:mplus}) and~(\ref{eq:mminus})
lie in $\{0,1,\ldots,r-2\}$ (whereas the Ramond case occurs when
the solution is $m^+ \equiv r-1$, $m^-\equiv r-1$)\cite{J2}. If
$\gcd(m^++1,m^-+1)$ is one, then these data completely determine
the $r$-spin structure on $X$, and in a canonical way, so that in
this case we have a well-defined gluing morphism
\begin{equation}
\label{eq:gluetr}
\rho: \mrbar_{\Gamma_1\sqcup\Gamma_2} =\mrbar_{\Gamma_1} \times
\mrbar_{\Gamma_2} \rTo \mrbar_{\Gamma}  \rInto
\mbar^{1/r,\mathbf{m}}_{g,n}.
\end{equation}
If, however, $d :=\gcd(m^++1,m^-+1)$ is greater than one, there is
no canonical morphism of stacks, as there is in (\ref{eq:gluetr}).
 To construct a gluing morphism would require an isomorphism
   $\phi:\ce_d^+|_{q^+} \irightarrow \ce_d^-|_{q^-}$
that makes $\ce_d$ Ramond at the node. In particular, $\phi$ must
be compatible with the isomorphism $(\ce^+_d)^{\otimes d}
\irightarrow \omega (-\sum m_i p_i)$ and $(\ce^-_d)^{\otimes d}
\irightarrow \omega(-\sum m_i p_i)$.  We call such an isomorphism
$\phi$ a \emph{gluing datum}.

\begin{df}
Given a prestable $(n+2)$-pointed curve $(X, p_1, \dots, p_{n+2})$
(not necessarily connected) and a $d$-th root $(\ce,b)$ of
$\omega_X(-\sum m_i p_i)$, such that $m_{n+1}=m_{n+2}=-1$, denote
by $\overline{X}$ the curve obtained from $X$ by identifying the
point $p_{n+1}$ with $p_{n+2}$.  A \emph{gluing datum} for
$(\ce,b)$ is an isomorphism $\phi:\ce|_{p_{n+1}} \irightarrow
\ce|_{p_{n+1}}$ which is compatible with the $d$-th root maps
$\ce|^{\otimes d}_{p_{n+1}} \rTo^{b}
\omega_{\overline{X}}(-\sum^n_{i=1} m_i p_i) \rTo^{b^{-1}}
\ce|^{\otimes d}_{p_{n+2}}.$
\end{df}

Returning to the case of a curve $X$ with a single node and two
irreducible components, since  the normalized curve has two
connected components, the root $(\ce_d, c_{d,1})$ is determined
\emph{up to non-canonical isomorphism} by $(\ce_r,c_{r,1})$.
Still, choosing one gluing datum $\phi:\ce_d^+|_{q^+} \irightarrow
\ce_d^-|_{q^-}$ does not give a morphism of stacks because an
automorphism of the $r$-spin structure on $X^+$ or on $X^-$
changes the gluing datum $\phi$, and thus induces a different (but
isomorphic) $r$-spin structure on the curve $X$.

\subsubsection{Irreducible curve with one node}\label{node}

\

In the case of an $r$-spin curve $X$ with a single component and
one node we have the dual graph $$ \Gamma_{irr}=\ocongmpm. $$ This
determines an $r$-spin structure on the normalized curve
$\tilde{X}$, with the dual graph $$
\tilde{\Gamma}_{irr}=\ocongmcut. $$ If $d=\gcd(m^+ +1,m^-
+1)=\gcd(m^++1, r)=1,$ then all roots $\ce_d$ in the $r$-spin
structure are of  the Neveu-Schwarz type and we can define the
gluing morphism
\begin{equation}
\rho:\mrbar_{\tilde{\Gamma}_{irr}} \rTo
\mrbar_{\Gamma_{irr}} \rInto \mbar^{1/r,\mathbf{m}}_{g,n}
\end{equation}
in the obvious way.

However, if $d=\gcd(m^+ +1,r)$ is greater than one, then, as
mentioned above, an additional gluing datum is required to
construct $(\ce_d, c_{d,1})$ from $(\ce^{\otimes r/d}_r,
c_{r,1})$. In this case, set $u^+=\frac{m^++1}{d}$ and
$u^-=\frac{m^-+1}{d}$ and define $$
  \ce'_d:=\ce^{\otimes r/d}_r \otimes \co (u^+p^++u^-p^-),
$$ where $p^+$ and $p^-$ are the inverse images under
normalization of  the node.  This shows that $\ce'_d$ is a $d$-th
root of $\omega_{\tilde{X}}(-\sum m_i p_i + p^+ +p^-)$; that is,
$m^+$ and $m^-$ are both replaced with $-1$. To construct $\ce_d$
from $\ce'_d $ we need to choose an isomorphism $$
  \phi:\ce'_d |_{p^+} \rTo \ce'_d|_{p^-}
$$ compatible with the isomorphisms $$ \ce^{'\otimes d}_d |_{p^+}
\irightarrow\omega_{\tilde{X}} (p^++p^-)\irightarrow \ce^{'\otimes
d}_d |_{p^-}, $$ and there are exactly $d$ such isomorphisms.

Unlike in the case of the tree, an automorphism of the $r$-spin
structure on the normalized curve induces the same automorphism on
both sides of the gluing datum $\phi$, and thus it preserves
$\phi$. Consequently, we expect $d=\gcd(m^+ +1,r)$ gluing
morphisms $$ \rho_{\phi}:\mbar^{1/r}_{\tilde{\Gamma}_{irr}} \rTo
\M^{1/r}_{\Gamma_{irr}}\rInto \M^{1/r, \mathbf{m}}_{g,n}, $$
indexed by the set of  different choices of $\phi$. However,  in
order to construct such a  morphism, one needs to be able to
define the gluing data in families.  That is, if $\pi :
\cc^{1/r}_{\tilde{\Gamma}_{irr}} \rTo
\M^{1/r}_{\tilde{\Gamma}_{irr}}$ is the universal curve over
$\M^{1/r}_{\tilde{\Gamma}_{irr}}$, and if $D^+$ and $D^-$ are the
loci in of the two sections $\sigma^+$ and $\sigma^-$ of $\pi$ to
be glued, then we need to define an isomorphism $$
\phi:\ce'_d|_{D^+} =\sigma^{+*}(\ce'_d)\irightarrow
\sigma^{-*}(\ce'_d)=\ce'_d|_{D^-}.$$  Such an isomorphism exists
precisely because $\ce'_d$ corresponds to the graph
$\tilde{\Gamma}_{irr}$ with $m^+=m^-=-1$. By Corollary
\ref{trivial}, both $\ce'_d|_{D^+}$ and $\ce'_d|_{D^-}$ are
trivial bundles on $\M^{1/r,(\bm,-1,-1)}_{g-1,n+2}$, and hence are
isomorphic.  Thus if $\mathfrak{g}$ is the set of isomorphisms
$\ce'_d|_{D^+} \irightarrow \ce'_d|_{D^-}$ of the sheaf
$\ce'_d=\ce^{\otimes r/d}_{r} \otimes \co(u^+D^++u^-D^-)$, we have
a morphism $$ \rho_{\Gamma_{irr}}:\M^{1/r}_{\tilde{\Gamma}_{irr}}
\cross \mathfrak{g} \rTo \M^{1/r}_{\Gamma_{irr}} \rInto
\M^{1/r,\bm}_{g,n}. $$

\section{Tautological cohomology classes}

Unless otherwise stated, all cohomology groups in the paper are considered
with coefficients in $\nq$.

\subsection{Definitions}

\

There are many natural cohomology classes in $H^\bullet(\mgnrmbar,\nq)$;
these include the classes induced by pullback from $\mgnbar$, as
well as classes induced by replacing the canonical (relative
dualizing) sheaf $\omega_{\pi}$ with $\ce_r$ in the usual
constructions of tautological cohomology classes on $\mgnbar$.

In particular, we have the $i$-th Chern classes $\lambda_i$ of the
Hodge bundle $\pi_* \omega_\pi$.  However, $\pi_* \ce_r$ is not
especially well behaved. Instead, we prefer to use the K-theoretic
pushforward $\pi_!\ce_r$ (also called $R \pi_* \ce$).\footnote{The
notation $\pi_!$ is from algebraic topology and is not to be
confused with the sheaf-theoretic \emph{direct image with compact
supports}, which we never use in this paper.} Recall that for any
coherent sheaf $\cf$ on the universal curve $\pi:\cgnrm \rTo
\mgnrmbar$, the element $\pi_! \cf$ of $K_0(\mgnrmbar)$ is the
difference $\pi_* \cf -R^1 \pi_* \cf$ ($\pi$ has relative
dimension $1$).  Here $R^1 \pi_* \cf$ is the sheaf whose fiber
over a point $p$ of the base $\mgnrmbar$  is the vector space
$H^1(\pi^{-1}(p),\cf|_{\pi^{-1}(p)})$. Serre duality shows that
$H^1(X,\omega)$ is canonically isomorphic to $\nc$, and so
$R^1\pi_* \omega_\pi$ is a trivial line bundle. Hence $\pi_*
\omega_\pi =\pi_! \omega_\pi+\co$ in $K_0(\mgnrmbar)$.  Therefore,
we have an equality of Chern polynomials $$ c_t \pi_* \omega_{\pi}
=c_t\pi_!\omega_{\pi} = 1+\lambda_1 t + \lambda_2 t^2 + \dots. $$
Tautological classes $\nu_i$ are defined as components of the
Chern character of the Hodge bundle
\begin{equation}
\label{eq:nu} ch_t \pi_* \omega_\pi = 1+ch_t \pi_! \omega_\pi
=g+\nu_1 t + \nu_2 t^3 +  \nu_3 t^5 + \dots.
\end{equation}
(The even components of $ch_t \pi_* \omega_\pi$ vanish by
Mumford's theorem~\cite{Mu}.) Similarly, we define classes $\mu_i$
as components of the Chern character of $\pi_! \ce_r$:
\begin{equation}
\label{eq:mu}
ch_t \pi_! \ce_r = D+ \mu_1t + \mu_2 t^2 + \dots .
\end{equation}
Here $-D$ is the Euler characteristic
$\chi(\ce_r|_{\cc_{\overline{s}}})$ of $\ce_r$ on any geometric
fiber $\cc_{\overline{s}}$ of $\pi$,  and by Riemann-Roch
\begin{equation}
\label{eq:degree}
D=\frac{1}{r}\left((r-2)(g-1)+\sum_i m_i\right).
\end{equation}

Serre duality shows that for any $\cf$ we have
$\pi_!(\sheafhom(\cf,\omega_{\pi}))=\pi_! \cf$, so  $\pi_! \co =
\pi_! \omega_\pi$.  More importantly, for purposes of comparison
with Witten's calculations of \cite{W}, in the special case that
$R^1 \pi_* \sheafhom (\ce_r, \omega_\pi) = 0$ (or equivalently,
$\pi_* \ce_r=0$), the bundle $\mathcal{V} := \pi_{*}
\sheafhom(\ce_r,\omega_\pi)$ of \cite{W} corresponds to $\pi_!
\ce_r.$

In addition to the Hodge-like classes, there are those induced by
the canonical sections $\sigma_i$ of $\pi: \cgnrm \rTo \mgnrmbar$.
These are classes $$ \psi_i := c_1(\sigma^*_i (\omega_{\pi}))
\quad \mathrm{and} \quad \tilde{\psi}_i := c_1(\sigma^*_i (\ce_r))
$$ (and also class $\tilde{\psi}^{(d)}_i$ for each divisor $d$ of
$r$). When working in $\pic \mgnrmbar$, we will abuse notation and
use $\psi_i$ to indicate the line bundle
$\sigma^*_i(\omega_{\pi})$, and $\tilde{\psi}_i$ the line bundle
$\sigma^*_i(\ce_r)$. Finally, there are the boundary classes. In
particular, if $A \sqcup B$ is a partition of $\{1, \ldots, n\}$
into two subsets, we denote by $\alpha_{k;A}$ the class of the
divisor associated to $r$-spin curves with the dual graph $\Gamma$
of the form $\tcongm$, with $\{i_1, \ldots, i_j\}=A$ and
$\{i_{j+1}, \ldots, i_n\}=B$.  Of course there is an obvious
equality: $$\alpha_{k;A} = \alpha_{g-k;B}.$$ Since the graph
$\Gamma$ is a tree, there is a unique choice of $m^+$ and $m^-$,
given the original type $\bm$ and the partition $A \sqcup B$.

If $g$ is greater than $1$, and if $r$ and all of the $m_i$ are
even, then the moduli space has two components, $\mgnrmbareven$
and $\mgnrmbarodd$. If $ 2 \le k \le g-2,$ then $\alpha_{k;A}$ is
the sum of four divisors---two on each irreducible component of
the moduli space.  In particular, there are two divisors in $\pic
\mgnrmbareven$ with dual graph $\tcongm$. The first is the locus
where both vertices of the graph (irreducible components of the
underlying curve) are endowed with an even $r$-spin structure; and
the second is where both vertices are endowed with an odd $r$-spin
structure.  Similarly, in $\pic \mgnrmbarodd$, the two divisors
correspond to the two ways of endowing the vertices with $r$-spin
structures of differing parities.

In the case of $k=0$ and $g>1$, the divisor $\alpha_{0;k}$ is the
sum of only two divisors, corresponding to the parity of the
$r$-spin structure on the other vertex (of genus $g$).  If $k=1$
then $\alpha_{1;A}$ is the sum of (potentially many) divisors
corresponding to the choices of index for the $r$-spin structure
on the vertex of genus $1$, (as well as the choices of index or
parity for the remaining vertex).

Finally, denote by $\tilde{\delta}_{irr,m^+}$ the divisor
associated to the graph $\Gamma= \ocongm$ with the $r$-spin
structure inducing $\ocongmcut$ on the normalization; and denote
the divisor corresponding to the Ramond root by
$\tilde{\delta}_{irr,r-1}$.\footnote{Beware that the index $m$ of
$\tilde{\delta}_{irr,m}$ differs from the index $j$ of $\gamma$ in
\cite{J3} by one.} The divisor $\tilde{\delta}_{irr,m^+}$ is not
necessarily irreducible, since different choices of gluing will
induce distinct (and disjoint) divisors, all in
$\tilde{\delta}_{irr,m^+}$. Again, there is an obvious equality:
$$\tilde{\delta}_{irr,m} = \tilde{\delta}_{irr,r-2-m}.$$

\subsection{Basic properties of the tautological classes}

\

The following Proposition describes relations between various
elements in the Picard group of $\mgnrmbar$. It is a
straightforward generalization of the corresponding result for the
case $\bm = \mathbf{0}$ proved in \cite{J3}.

\begin{prop}
\

\begin{itemize}
\item
The forgetful map $p: \mgnrmbar \rTo \mgnbar $ induces an
injection $$ p^*: \pic \mgnbar \otimes \mathbb{Q} \rTo \pic
\mgnrmbar \otimes \mathbb{Q}. $$
\item
Let $\delta_{k;A}$ denote the pullback to $\pic \mgnrmbar$ of the
class in $\pic \mgnbar$ associated to the union of all strata in
$\mgnbar$  with the dual graph $\Gamma =\tcongi$, with $A=\{i_1,
\ldots, i_j\}$.  The pullback $\delta_{k;A}$ is related to
$\alpha_{k;A}$ as follows.  $$ \delta_{k;A} =
\frac{r}{\gcd(m^++1,r)} \alpha_{k;A}, $$ where $m^+$ is determined
by $k, \bm$, and $A$, as in Section~\ref{twocomp}.
\item
Let $\delta_{irr}$ be the pullback to $\pic \mgnrmbar$ of the
divisor of all curves in $\mgnbar$ with dual graph
$\Gamma_\mathrm{irr} = \ocongi$.  The pullback $\delta_{irr}$ can
be expressed in terms of the $\tilde{\delta}_{irr,m}$ as follows.
$$ \delta_{irr}= \sum_{r-2/2 \leq m < r} \frac{r}{\gcd(r,m+1)}
\tilde{\delta}_{irr,m}. $$
\end{itemize}
\end{prop}

\medskip

The fact that, for an $r$-th  root $(\ce, b)$ of $\omega(-\sum
m_ip_i)$, the map $b$ is almost an isomorphism means that
$\tilde{\psi}_i$ and $\psi_i$ are closely related.

\begin{prop}\label{bundles}
The line bundles $\sigma^*_i\omega $ and $\sigma^*_i(\ce_d)$ on
the stack $\mgnrmbar$ are related by $$ r \sigma^*_i(\ce_r) \cong
(m_i +1) \sigma^*_i (\omega) $$
   and
$$
  d \sigma^*_i(\ce_d)\cong (m'_i+1) \sigma^*_i(\omega),
$$ where $m'_i$ is the smallest non-negative integer congruent to
$m_i \pmod d$.  Therefore, in $\pic \mgnrmbar \otimes \mathbb{Q}$
we have $$ \tilde{\psi}_i= \frac{m_i+1}{r} \psi_i $$
   and
$$
\tilde{\psi}^{(d)}_i = \left(\frac{m'_i +1}{m_i+1}\right) \frac{r}{d}
\tilde{\psi}_i.
$$
\end{prop}

Before proving the proposition, we recall the following well-known
fact.

\begin{lm}\label{psi}
Let $\pi: \cgn \rTo \mgnbar$ be the universal $n$-pointed curve,
and $\sigma_i$ the $i$-th  tautological section of $\pi$.  If
$D_i$ is the divisor of $\cgn$ associated to $\sigma_i$, and if
$\omega=\omega_{\pi}$ is the canonical (relative dualizing)
sheaf, then
$$
\sigma^*_i (\co (-D_j)) \cong
\begin{cases}
\sigma^*_i (\omega|_{D_i}) = \psi_i & \mathrm{if} \ i=j \\
  \co & \mathrm{if} \ i \ne j.
\end{cases}
$$
\end{lm}
\begin{proof}{(of Lemma \ref{psi})}
When $i \neq j$ the bundle $\sigma^*_i (\co (-D_j))$ is trivial
because the sections are disjoint.  In the case $i=j$ the result
follows from the fact that taking residues gives an  isomorphism
between $\omega_{\pi} (D_i)|_{D_i}$ and $\co_{D_i}$.
\end{proof}

\begin{proof}{(of Proposition \ref{bundles})}
 The map $c_{d,1}: \ce^{\otimes d}_d  \rTo \omega(-\sum m'_jp_j)$ pulls
back, via $\sigma^*_i$, to give $$ \sigma^*_i c_{d,1}: (\sigma^*_i
\ce)^{\otimes r} \rTo \sigma^*_i \omega(-\sum m'_jp_j). $$
  Since $\operatorname{im}(\sigma_i)$
is disjoint from the nodes of $X$, $\sigma^*_ic_{d,1}$ is an
isomorphism, even on the boundary strata where $\ce$ fails to be
locally free.

Consequently, we have
$$\sigma^*_i(c_{d,1}):\sigma^*_i(\ce^{\otimes d}_d) \irightarrow
\sigma^*_i (\omega \otimes \co(-\sum_j m'_j
D_j))=\psi_i+m'_i\psi_i=(m'_i+1)\psi_i.$$
\end{proof}

\subsection{Non-trivial relations involving the class $\mu_1$}
\

\begin{prop} \label{prop:murelation}
Define the boundary divisors $   \epsilon $
 and $\delta$ by
 \begin{eqnarray*}
\epsilon=&& \sum_{k,A} \frac{(m^++1)(m^-+1)}{u_{k;A}}\alpha_{k;A}
\\ && + \sum_{\frac{r}{2}-1 \leq m < r -
1}\frac{(m+1)(r-m-1)}{v_m} \tilde{\delta}_{irr,m}
\end{eqnarray*}
and
$$
\delta= \sum \delta_{k;A} = \sum_{\frac{r}{2}-1 \leq m < r}
\frac{r}{v_m} \tilde{\delta}_{irr,m} + \sum_{k,A} \frac{r}{u_{k;A}}
\alpha_{k;A},
$$
where $v_m = \gcd(m+1,r)$, $u_{k;A} =\gcd (m^++1, r)$, and $m^+$ is
determined by $k,A$, and $\bm$
via relation~(\ref{eq:mplus}).

  Then the
following relation holds in $\pic \mgnrmbar $: $$ r \epsilon =
(2r^2-12r+12)\lambda_1-2r^2 \mu_1 +(r-1) \delta+\sum_{1\leq i\leq
n} m_i(r-2-m_i)\psi_i. $$
\end{prop}

The proof of the Proposition is almost identical to its
counterpart in \cite[Theorem 4.3.3]{J3} except that $\ce_r$ is not
an $r$-th  root of $\omega$, but rather of $\omega(-\sum m_ip_i)$.
The only extra information necessary to prove the proposition is
the content of the following two lemmas and the fact that the
divisor we have called $\epsilon$ is the product $<\tilde{\ce_r,
\mathfrak{E}}>$ in \cite{J3}.

\begin{lm}\label{deligne}
Let $\pi:\cc_{g,n}\rTo \mgnbar$ be the universal curve and $$ <,>
:\pic \cc_{g,n} \cross \pic \cc_{g,n} \rTo \pic \mgnbar $$ be
Deligne's bilinear product defined by $$ <\cl,\cm>:= \det(\pi_!
(\cl \otimes \cm) - \pi_! \cl -\pi_! \cm +\pi_! \co). $$ If $D$ is
the image of a section $\sigma: \mgnbar \rTo \cgn$, then for any
line bundle $\cl$ on $\cgn$ the product  $<\cl,\co_{\cc}(D)>$ is
equal to the restriction of $\cl$ to $D$; that is, it is just
$\sigma^*\cl$.
\end{lm}

This is proved in \cite[Prop. 6.1.3]{De}, but it also follows from
the fact that when the base is a smooth curve $B$, then
$\deg\!<\cl,\cm>$ is exactly the usual intersection number $(\cl.
\cm)$.  Since line bundles on $\mgnbar$ are completely determined
by their degree on smooth curves in $\mgnbar$
\cite{arbarello-cornalba:mg-pic}, and since in the case of a
smooth, one-parameter base the degrees agree, the lemma is true in
general.

\begin{lm} We have
$$
<\co(D_i), \co(D_j)>=
\begin{cases} - \sigma^*_i \omega_{\pi} = -\psi_i &\mathrm{if\ } i=j\\
0 & \mathrm{if\ } i \neq j.\end{cases}
$$

\end{lm}

This follows immediately from Lemmas \ref{psi} and \ref{deligne}.

\begin{crl}
If $g=0$ then classes $\lambda_1$ and $\tilde{\delta}_{irr,m}$
vanish (for all $m$),  and
\begin{eqnarray*}
2r^2\mu_1&=&(r-1)\delta +\sum_i m_i(r-2-m_i)\psi_i -r \epsilon \\
&=& \sum_{A \subseteq [n]}
\frac{r}{u_{0;A}}[r-1-(m^++1)(m^-+1)]\alpha_{0;A}
\\ && + \sum_{1 \leq i \leq n} m_i (r-2-m_i) \psi_i.
\end{eqnarray*}
\end{crl}

\section{Cohomological field theory}
\label{cft} In this section we begin with a review of the notion
of a cohomological field theory (\cft) in the sense of Kontsevich
and Manin \cite{KM1}. This is an object which formalizes the
expected factorization properties of the theory of topological
gravity coupled to topological matter. The  Gromov-Witten
invariants associated to a smooth, projective  variety $V$
correspond to the physical situation where the matter sector
arises from the topological sigma model \cite{KM1,RT}. The
analogous intersection numbers associated to the moduli space of
$r$-spin curves have their physical origins in a different choice
of the matter sector. Our goal in this section is to give a
precise formulation of these notions in terms of the  moduli
spaces described above.

\subsection{Axioms of \cft}
\begin{df}
A \emph{(complete)  cohomological field theory (\cft) of rank $d$} (denoted
by $(\ch,\eta,\Lambda)$ or just $(\ch,\eta)$) is a $d$-dimensional vector
space $\ch$ with a metric $\eta$ and a collection $\Lambda\,:=\,
\{\,\Lambda_{g,n}\,\}$ of $n$-linear $H^\bullet(\mgnbar)$-valued forms on
$\ch$
\begin{equation}\label{eq:lambda}
\Lambda_{g,n} \in \,H^\bullet(\mgnbar)\,\otimes \ch^{*\tensor n} =
\operatorname{Hom}(\ch^{\tensor n},H^\bullet(\mgnbar))
\end{equation}
defined for stable pairs $(g,n)$ and satisfying the following axioms
{\bf C1--C3} (where $\{e_0,\ldots,e_{d-1}\}$ is a fixed basis of
$\ch$, $\eta^{\mu\nu}$ is the inverse of the matrix of the metric
$\eta$ in this basis, and the summation convention has been
used).
\begin{enumerate}
\item[\bf C1.] The element
$\Lambda_{g,n}$ is invariant under the action of the symmetric
  group $S_n$.
\item[\bf C2.]
Let
\begin{equation}
\label{eq:gluetree}
\rho_{\mathrm{tree}}:
\M_{\Gamma_1 \sqcup \Gamma_2}=\M_{\Gamma_1} \times \M_{\Gamma_2}
\rTo \M_{\Gamma_\mathrm{tree}}  \rInto \M_{g,n}
\end{equation}
be the gluing morphism~(\ref{eq:gluing}) corresponding to the stable graph
$$
\Gamma_\mathrm{tree}=\tcongi
$$
and the two graphs $\Gamma_1$ and $\Gamma_2$ obtained by cutting
the edge of $\Gamma_\mathrm{tree}$. Then the forms
$\Lambda_{g,n}$ satisfy the composition property
\begin{eqnarray}
\label{eq:cfttree}
\rho_{\mathrm{tree}}^* \Lambda_{g,n}(\gamma_1,\gamma_2\,\ldots,\gamma_n)&=&
\\
\Lambda_{k,j}(\gamma_{i_1},\ldots,\gamma_{i_j},e_\mu)
\eta^{\mu\nu} &\otimes&
\Lambda_{g-k,n-j}(e_\nu\,\gamma_{i_{j+1}},\ldots,\gamma_{i_n})
\nonumber
\end{eqnarray}
for all $\gamma_i \in \ch$.

\item[\bf C3.]
Let
\begin{equation}
\label{eq:glueloop}
\rho_{\mathrm{loop}}:
\M_{\tilde{\Gamma}}=\M_{g-1,n+2} \rTo \M_{\Gamma_\mathrm{loop}} \rInto
\M_{g,n}
\end{equation}
be the gluing morphism~(\ref{eq:gluing}) corresponding to the stable graph
\begin{equation}
  \label{eq:gtree}
\Gamma_\mathrm{loop}=\ocongi
\end{equation}
and the graph $\tilde{\Gamma}$ obtained by cutting the loop of
$\Gamma_\mathrm{loop}$.
Then
\begin{equation}
\label{eq:cftloop}
\rho_{\mathrm{loop}}^*\,\Lambda_{g,n}(\gamma_1,\gamma_2,\ldots,\gamma_n)\,=\,
\Lambda_{g-1,n+2}\,(\gamma_1,\gamma_2,\ldots,\gamma_n, e_\mu,
e_\nu)\,\eta^{\mu\nu}.
\end{equation}
\medskip
The pair $(\ch,\eta)$ is called the \emph{state space} of the \cft.
\medskip

An element $e_0\in \ch$ is called \emph{a flat identity} of the
\cft\ if, in addition, the following equations hold.
\item[\bf C4a.]
For all $\gamma_i$ in $\ch$ we have
\begin{equation}
  \label{eq:identity}
\Lambda_{g,n+1}(\gamma_1,\ldots, \gamma_n, e_0)\, =
\,\pi^*\Lambda_{g,n}(\gamma_1,\ldots, \gamma_n),
\end{equation}
where $\pi:\M_{g,n+1} \to \M_{g,n}$ is the universal curve on $\mgnbar$ and
\item[\bf C4b.]
\begin{equation}
\label{eq:identity2}
\int_{\M_{0,3}}\,\Lambda_{0,3}(\gamma_1,\gamma_2,e_0) =
\eta(\gamma_1,\gamma_2).
\end{equation}
\end{enumerate}
A \cft\ with flat identity is denoted by $(\ch,\eta,\Lambda,e_0).$
\medskip
A \emph{genus $\tilde{g}$ \cft} on the state space $(\ch,\eta)$ is the
collection of forms
$\{\,\Lambda_{g,n}\,\}_{g\,\leq\,\tilde{g}}$ that satisfy only those of the
equations~(\ref{eq:cfttree}),~(\ref{eq:cftloop}), ~(\ref{eq:identity}),
and~(\ref{eq:identity2}), where $g\le \tilde{g}$.
\end{df}

\begin{rems}
\

{\bf 1.} In general, the state space $\ch$ of \cft\ is
$\mathbb{Z}_2$-graded, but here, for simplicity, we are assuming
that $\ch$ contains only even elements, since this is the only
case that will arise in this paper.

{\bf 2.} The definition of a \cft\ given above has an equivalent
dual description in terms of homology. Consider the maps
$H_\bullet(\mgnbar)\, \to\, T^n \ch^*$ given by
$[c]\,\mapsto\,\int_{[c]}\Lambda_{g,n}$. These maps are called the
\emph{($n$-point) correlators of the \cft}. A structure of a
(complete) \cft\ on $(\ch,\eta)$ is equivalent to the requirement
that these correlators endow $(\ch, \eta)$ with the structure of
an algebra over the modular operad $\{\,H_\bullet(\mgnbar)\,\}$ in
the sense  of Getzler and Kapranov \cite{GK}.

{\bf 3.}
Clearly, the definition of a cohomological field theory extends
from $\nc$ to more general ground rings $\ck$.
\end{rems}

Let $\Gamma$ be a stable graph, then there is a canonical composition
map $\rho_\Gamma$
\begin{equation}
\label{eq:totglue}
\rho_\Gamma:\prod_{v\in V(\Gamma)}\,\M_{g(v),n(v)}\,\to\,
\M_\Gamma\,\to\,\mgnbar
\end{equation}
where $V(\Gamma)$ denotes the set of vertices of $\Gamma$.
Since the map $\rho_\Gamma$ can be constructed from
gluing morphisms~(\ref{eq:gluetree}) and~(\ref{eq:glueloop}),
the forms $\Lambda_{g,n}$ satisfy a restriction property
\begin{equation} \label{restriction}
\rho_\Gamma^*\Lambda_{g,n}\,=\, \rho_\Gamma^{-1}(\,\bigotimes_{v\in V(\Gamma)}
\Lambda_{g(v),n(v)}\,)
\end{equation}
where \[ \rho^{-1}_\Gamma\,:\,\bigotimes_{v\in V(\Gamma)}
T^{n(v)}\,\ch^*
\,\to\, T^n\,\ch^* \]
contracts the factors $T^n\ch^*$ by means of the inverse of the metric $\eta$
and successive application of equations~(\ref{eq:cfttree})
and~(\ref{eq:cftloop}) .
There is a parameter which can be introduced into the definition of a
\cft. This parameter can be regarded as a coupling constant in the theory.
\begin{lm}
\label{lm:CouplingConstant} Let $(\ch, \etat,\Lambdat,e_0)$ be a
\cft\ with flat identity $e_0$ and let $\lambda$ be a nonzero
parameter. If we define $\Lambda\,=\,\{\,\Lambda_{g,n}\,\}$, where
\[
\Lambda_{g,n}\,:=\,\lambda^{2g-2}\,\Lambdat_{g,n}
\]
and
\[
\eta\,:=\,\lambda^{-2}\,\etat ,
\]
then $(\ch, \eta, \Lambda, e_0)$ is a \cft\ with flat identity.
\end{lm}
The proof is obvious.
\begin{df}
The {\em small phase space potential function} of the \cft\
$(\ch,\etat,\Lambdat)$ is a formal series $\Phi\,\in\,\nc[[\ch]]$ given by
\begin{equation}
   \label{eq:small}
\Phi(\bx)\, := \,\sum_{g\,=\,0}^\infty\,\Phi_g(\bx),
\end{equation}
where
\[
\Phi_g(\bx)\,:=\,\sum_n\,\frac{1}{n!}\,\int_{\mgnbar}\,
\left<\,\Lambdat_{g,n}\,,\,\bx^{\otimes n}\, \right>.
\]
Here $\left<\,\cdots\,\right>$ denotes evaluation, the sum
over $n$ is understood to be over the stable range, and
$\bx\,=\,\sum_{\alpha}\,x^\alpha\,e_\alpha$, where
$\{\,e_\alpha\,\}$ is a basis of $\ch$.
\end{df}
\begin{rem}
The small phase space potential function of the \cft\
$(\ch,\eta,\Lambda,e_0)$ associated to $(\ch,\etat,\Lambdat,e_0)$ as in Lemma
\ref{lm:CouplingConstant} may be regarded as an element in
$\lambda^{-2}\nc[[\ch,\lambda^{2}]]$.
\end{rem}
All of the information of a genus zero \cft\  is encoded in this potential.
\begin{thm} \cite{KM1,Ma}
\label{thm:wdvvg}
An element $\Phi_0$ in $\nc[[\ch]]$ is the potential of a rank
$d$, genus zero \cft\ $(\ch,\eta)$ if and only if it contains
only terms which are of cubic and higher order in the coordinates
$x^0,\ldots,x^{d-1}$ (corresponding to a basis
$\{\,e_0,\,\ldots\,e_{d-1}\,\}$ of $\ch$) and it satisfies the
{\em associativity, or WDVV} (Witten-Dijkgraaf-Verlinde$^2$)
equation
\[
\label{eq:WDVV}
\partial_{a} \partial_{b} \partial_{e}\Phi_0\, \eta^{ef}\, \partial_{f}
\partial_{c} \partial_{d}\Phi_0\, = \, \partial_{b} \partial_{c}
\partial_{e} \Phi_0\, \eta^{ef} \,\partial_{f} \partial_{a}
\partial_{d}\Phi_0,
\]
where $\eta^{ab}$  is the inverse matrix of the matrix of
$\eta$ in the basis $\{e_a\}$,  $\partial_a$ is derivative with
respect  to $x^a$, and the summation convention has been used.

Conversely,  a genus zero \cft\ structure on $(\ch,\eta)$ is
uniquely determined by its potential $\Phi_0$, which must satisfy
the WDVV equation.
\end{thm}

A genus zero \cft\ with flat identity is essentially equivalent to
endowing the state space $(\ch,\eta)$ with the structure of a
formal Frobenius manifold \cite{Du,Hi,Ma}. The theorem follows
from the  work of Keel \cite{Ke}, who proved that
$H^\bullet(\M_{0,n})$ is generated by boundary classes and that
all relations between boundary divisors arise from lifting the
basic codimension one relation on $\M_{0,4}$. Higher genus \cfts\
are not completely determined by their potentials because
$H^\bullet(\mgnbar)$ is not generated by its boundary classes.
Nonetheless, Getzler \cite{Ge} proved that the potential function
of a genus one \cft\ must satisfy an analogous (but much more
complicated) equation arising from a codimension-two relation in
$H^\bullet(\M_{1,4})$.

\subsection{Gromov-Witten invariants and their potentials}
\
Our construction of \cfts\ from the moduli space of stable
$r$-spin curves is guided by analogy with the moduli space of
stable maps and Gromov-Witten invariants. Let us briefly review
this construction. Let $V$ be a smooth projective variety,
$\ch\,=\,H^\bullet(V,\nc)$, and $\eta$ the Poincar\'e pairing. Let
$\mgnv$ be the moduli stack of stable maps into $V$ of genus $g$
with $n$ marked points. The \emph{Gromov-Witten invariants of $V$}
are multilinear maps $\ch^{\otimes n}\,\to\,\nc$ given by
\begin{equation}
\cor{\tau_0(\gamma_1)\,\cdots\,\tau_0(\gamma_n)}_g\,=\,
\lambda^{2g-2}
\int_{[\mgnv]^{\mathrm{vir}}}\,\ev_1^*\,\gamma_1\cup\cdots\cup\ev_n^*\,
\gamma_n,
\end{equation}
where $[\mgnv]^{\mathrm{vir}}$ is the virtual fundamental class of the
moduli stack $\mgnv$ and $\lambda$ is a formal parameter.
The corresponding \emph{small phase space potential}
$\Phi(\bx)$
is defined  by~(\ref{eq:small})
where the genus $g$ part is given by
\[
\Phi_g(\bx)\,=\,\,\sum_n\,\lambda^{2g-2}\int_{[\mgnv]^{\mathrm{vir}}}\,
\ev_1^*\,\bx\cup\cdots\cup\ev_n^*\,\bx,
\]
$\bx\,=\,\sum_a\,x^a\,e_a$, and $\{\,e_0,\ldots,e_{r}\,\}$ is a
basis for $\ch$ such that $e_0$ is the identity element. If $V$ is
a convex variety, then $\M_{0,n}(V)$ is a smooth stack and its
virtual fundamental class coincides with its topological
fundamental class. In this situation, \cite{FP} shows that
$(\ch,\eta)$ forms a genus zero \cft\ with potential $\Phi_0$.
This result can be generalized to higher genera and to more
general varieties, as well.
\begin{rem}
In the usual definition of Gromov-Witten invariants, there is no factor of
$\lambda^{2g-2}$ in the definition of the correlators but this factor is
inserted into the potential function by hand. We have chosen our conventions
so that this factor appears instead in the correlator but is not explicitly
inserted into the potential function.
\end{rem}
The \emph{gravitational descendants} are defined by twisting the
Gromov-Witten classes with the tautological $\psi$ classes as
follows:
\begin{equation}
  \label{eq:gravit}
\left<
  \tau_{a_1}(\gamma_1)\cdots\tau_{a_n}(\gamma_n)\right>_g
  :=\lambda^{2g-2}\,\int_{[\mgnv]^\mathrm{vir}}\!\ev_1^*\gamma_1\cup
  \psi_1^{a_1}\cup\cdots \cup  \ev_n^*\gamma_n \cup \psi_n^{a_n}
\end{equation}
for all $a_1,\,\ldots\,,a_n\,=\,0,1,2,\ldots$ and $\gamma_1,
\,\ldots,\,\gamma_n$ in
$H^\bullet(V)$. This gives rise to the \emph{large phase space potential}
$\Phi(\bt)\,\in\,\lambda^{-2}\nc[[\bt,\lambda^2]]$ where
$\bt\,=\,(\bt_0,\bt_1,\ldots)$ and $\bt_n\,=\,(t_n^0,\ldots,t_n^r)$,
which is defined by
\begin{equation}
\label{eq:LPSPotential}
\Phi(\bt)\,:=\,\sum_g\,\Phi_g(\bt)
\end{equation}
and
\begin{equation}
\label{eq:LPSPotentialg}
\Phi_g(\bt)\,:=\,\sum\,\left<\,
  \tau_{a_1}(e_{\alpha_1})\,\cdots\,\tau_{a_n}(e_{\alpha_n})\,
  \right>_g\,t_{a_1}^{\alpha_1}\,\cdots\,t_{a_n}^{\alpha_n}\,\frac{1}{n!}.
\end{equation}
Setting $t_n^\alpha\,=\,0$ for $n\,\geq\,1$ and $x^\alpha\,=\,t_0^\alpha$
reduces the large phase space potential $\Phi(\bt)$ to the small phase space
potential $\Phi(\bx)$.

When $V$ is a point, Kontsevich's theorem gives that $$
Z(\bt)\,:=\,\exp(\Phi(\bt)) $$ is a $\tau$-function of the $\KdV$
hierarchy. In addition, Kontsevich showed that $Z(\mathbf{t})$ is
a highest weight vector for the Virasoro Lie algebra, a condition
which allows one to completely solve for these intersection
numbers. The existence of a similar Virasoro algebra action has
been conjectured by Eguchi, Hori, and Xiong \cite{EHX} in the case
where $V$ is not a point. Evidence for this conjecture is mounting
\cite{GePa,FaPa2}. A \emph{very large phase space} has recently
been introduced in \cite{KK,MaZ} for the case where $V$ is a point
and for more general varieties in \cite{FaPa} by including
variables corresponding to the Hodge classes $\nu_i$ as well.
These additional variables parametrize an even larger family of
\cfts\ than just the large phase space coordinates \cite{KK}. We
will shortly introduce, in addition, variables associated to the
$r$-spin structure (see~(\ref{eq:vlarge})).

\subsection{$r$-spin \cft}
\
We perform an analogous construction of a very large phase space
where the role of the moduli space of stable maps $\M_{g,n}(V)$
is played by the moduli space of stable $r$-spin curves
$\M_{g,n}^{1/r}$. Unlike the moduli space of general stable maps,
the moduli space of stable $r$-spin curves is a smooth stack.
Intersection theory is therefore simpler in this case than for
the case of stable maps. However, the difficulty lies instead in
the construction of the analogs of the Gromov-Witten classes.

In the next section, we will introduce axioms which a collection
of cohomology classes (called a \emph{virtual class})
$\cv_{g,n}(\bm)$ in $H^\bullet(\mgnrmbar)$ must satisfy in order
to insure that the following result holds.

\begin{thm}
\label{th:rcft}
Let $(\chr,\etat)$ be a vector space of dimension $r-1$ with a
basis  $\{\,e_0,\,\ldots,\, e_{r-2}\,\,\}$ and metric $\etat$
given by
\begin{equation}
\etat(e_\mu,e_\nu)\,:=\,\etat_{\mu\nu}\,:=\,\frac{1}{r}\,\delta_{\mu+\nu,r-2}.
\end{equation}
Let $\cv_{g,n}(\bm)$ be a virtual class in $H^\bullet(\mgnrmbar)$
satisfying Axioms 1 through 5 from the  next section, and let
$p\,:\,\mgnrmbar\,\to\,\mgnbar$ be the map which forgets the
$r$-spin structure. Let $\Lambdat\,:=\,\{\,\Lambdat_{g,n}\,\}$ be defined by
\begin{equation}
\Lambdat_{g,n}^{(\bs,\bu)}(e_{m_1},\ldots,e_{m_n})\,:=\,p_*\,(\cv_{g,n}(\bm)
\,\exp(\bs\,\cdot\,\bmu\,+\,\bu\,\cdot\,\bnu)),
\end{equation}
where these forms have values in the ring $\nc[[\bs,\bu]]$, then
$(\ch,\etat,\Lambdat^{(\bs,\bu)},e_0)$ is a \cft\  satisfying
Axiom C4a. Furthermore, if
\begin{equation}
\Lambdat_{g,n}\,:=\,\Lambdat_{g,n}^{(\mathbf{0},\mathbf{0})},
\end{equation}
then $(\ch,\etat,\Lambdat,e_0)$ is a \cft\ with flat identity. The latter
will be called {\em the $r$-spin \cft}.
Restricting the $r$-spin \cft\ to genus zero shows that $(\chr,\eta)$ is
endowed with the structure of a Frobenius manifold.
\end{thm}
This theorem is proved in Section \ref{sec:verify}.
\begin{crl}
\label{crl:rCFT}
Let $(\ch,\eta,\Lambda,e_0)$ be constructed from $(\ch,\etat,\Lambdat,e_0)$
above by setting
\begin{equation}
\eta(e_\mu,e_\nu)\,:=\,\eta_{\mu\nu}\,:=\,\frac{1}{r\lambda^2}\,
\delta_{\mu+\nu,r-2},
\end{equation}
\begin{equation}
\Lambda_{g,n}^{(\bs,\bu)}(e_{m_1},\ldots,e_{m_n})\,:=\,\lambda^{2g-2}\,
\Lambdat_{g,n}^{(\bs,\bu)}(e_{m_1},\ldots,e_{m_n}),
\end{equation}
and
\begin{equation}
\Lambda_{g,n}\,:=\,\Lambda_{g,n}^{(\mathbf{0},\mathbf{0})},
\end{equation}
then $(\ch,\eta,\Lambda^{(\bs,\bu)},e_0)$ is a \cft\ satisfying Axiom C4a
and  $(\ch, \eta, \Lambda, e_0)$ is a \cft\ with flat identity.
\end{crl}
\begin{proof}
This is a direct consequence of Lemma \ref{lm:CouplingConstant} and the
previous theorem.
\end{proof}
The classes $\cv_{g,n}(\bm)$ are analogs of the Gromov-Witten classes in this
theory. The analogs of the gravitational descendants~(\ref{eq:gravit}) are
given by
\begin{equation}
\label{eq:gravspin}
\cor{\tau_{a_1}(e_{m_1})\,\cdots\,\tau_{a_n}(e_{m_n})\,}_g\,:=\,
\lambda^{2 g-2}\int_{\mgnrmbar}\,
\psi_1^{a_1}\,\cdots\,\psi_n^{a_n}\,\cv_{g,n}(\bm),
\end{equation}
and the large phase space potential function is defined by the
equations (\ref{eq:LPSPotential}) and (\ref{eq:LPSPotentialg}).
The small phase space potential function is defined by
restricting~(\ref{eq:gravspin}) to correlators with $a_i=0$.
\emph{We will see that the case of
$\lambda\,=\,\frac{1}{\sqrt{r}}$
  corresponds to the generalized Witten conjecture. This corresponds to the
  metric
\begin{equation}
\label{eq:metric}
\eta(e_{m_1},e_{m_2})\,=\,\delta_{m_1+m_2,r-2}
\end{equation}
and the forms
\begin{equation}
\label{eq:rlambda}
\Lambda_{g,n}(e_{m_1},\ldots,e_{m_n})\,:=\,r^{1-g}\,p_*\,\cv_{g,n}(\bm).
\end{equation}
}
\begin{rem}
The astute reader will notice in the sequel that, strictly speaking,
the state space of this $r$-spin \cft~(\ref{eq:rlambda}) should be instead,
$(\chh,\hat{\eta})$ where $\chh$ is an $r$-dimensional vector
space with basis $\{e_0,\ldots,\- e_{r-1}\}$ and a metric
given by
$$
\hat{\eta}(e_a,e_b)=
\begin{cases}
  1 &\mathrm{if} \ a + b \equiv \,(r-2)\,\mod\,r \\
  0 & \mathrm{otherwise}.
\end{cases}
$$
However, it follows from the axioms for $\cv$ in
Section~\ref{axioms} that the obvious orthogonal decomposition
$\chh\,=\,\chr\, \oplus\,\ch',$ where $\ch'$ is the trivial
one-dimensional \cft\  with basis $\{\,e_{r-1}\,\}$, is a direct
sum of
\cfts.  For this reason, we can (and will) restrict
ourselves to the state space $(\chr,\eta).$
\end{rem}

\section{Virtual classes}\label{vc}
To endow the pair $(\ch^{(r)}, \eta)$ from Theorem~\ref{th:rcft}
with the structure of a \cft\ by equation~(\ref{eq:rlambda}), we
must define cohomology classes $\cv_{g,n}(\bm)$. We will call this
collection of classes an \emph{$r$-spin virtual class}. It should
satisfy the axioms described below. \emph{Throughout this section,
we will restrict ourselves to the case where the coupling constant
$\lambda$ is $\frac{1}{\sqrt{r}}$, unless otherwise stated.} This
is done purely for convenience as analogous results hold for
general $\lambda$ as well.
\subsection{Axioms for the virtual class}\label{axioms}
\
\begin{df}
An \emph{$r$-spin virtual class} is an assignment of a cohomology
class
\begin{equation}
\label{eq:cvirt} \cv_\Gamma \in H^{2 D} (\mbar^{1/r}_\Gamma, \nq)
\end{equation}
 to every genus $g$, stable, decorated graph $\Gamma$ with
 $n$-tails.
Here, if the tails of $\Gamma$ are marked with the $n$-tuple
$\bm=(m_1,\ldots,m_n)$, then the dimension $D$ is
\begin{equation}
\label{eq:deg} D = \frac{1}{r}\bigl( (r-2)(g-\alpha)+\sum_{i=1}^n
m_i \bigr),
\end{equation}
and $\alpha$ is the number of connected components of $\Gamma$. In
the special case where $\Gamma$ has one vertex and no edges, we
denote $\cv_\Gamma$ by $\cv_{g,n}(\bm)$.   These classes must
satisfy the axioms below.
\begin{description}
\item[\bf Axiom 1a (Connected Graphs)]
Let $\Gamma$ be a connected, genus $g$, stable, decorated graph
with $n$ tails. Let $E(\Gamma)$ denote the set of edges of
$\Gamma$. For each edge $e$ of $\Gamma$, let
$l_e\,:=\,\gcd(m_e^++1,r)$, where $m_e^+$ is an integer decorating
a half-edge of $e$. The classes $c^{1/r}_{\Gamma}$ and
$c^{1/r}_{g,n}(\bm)$ are related by
\begin{equation}
\label{eq:cvgamma} \cv_{\Gamma}= \left(
\prod_{e\,\in\,E(\Gamma)}\,\frac{r}{l_e}\right) \, \tilde{i}^*
\,\cv_{g,n}(\bm) \in H^{2 D} (\mbar^{1/r}_{\Gamma}),
\end{equation}
where $\tilde{i}\, :\,\M_\Gamma^{1/r}\,\rInto\M_{g,n}^{1/r,\bm}$ is the
canonical
inclusion map.
\item[\bf Axiom 1b (Disconnected Graphs)]
Let $\Gamma$ be a stable, decorated graph which is the disjoint
union of connected graphs $\Gamma^{(d)}$, then the classes
$\cv_{\Gamma}$ and $\cv_{\Gamma^{(d)}}$ are related by
\[
\cv_\Gamma\,=\,\bigotimes_d\,\cv_{\Gamma^{(d)}}\,\in\,H^\bullet(\M^{1/r}_\Gamma).
\]
\smallskip

\item[\bf Axiom 2  (Convexity)]
Consider the universal $r$-spin structure $(\{\ce_d\},\{
c_{d,d'}\})$ on the universal curve $\pi:\cc^{1/r,\bm}_{g,n}\rTo
\mbar^{1/r,\bm}_{g,n}$. For each irreducible (and connected)
component of $\M^{1/r,\mathbf{m}}_{g,n}$ (denoted here by
$\M^{1/r,\mathbf{m},(d)}_{g,n}$  for some index $d$), if
$\pi_*\ce_r=0$ on $\M^{1/r,\mathbf{m},(d)}_{g,n}$, then
$c^{1/r}_{g,n}(\mathbf{m})$ restricted to
$\M^{1/r,\mathbf{m},(d)}_{g,n}$ is $c_D (-R^1\pi_*\mathcal{E}_r)$,
the top Chern class of the bundle with fiber $H^1(X,\ce_r)^*$ at
$[(X,p_1,\dots,p_n,(\{\ce_d\},\{c_{d,d'}\}))], \in
\M^{1/r}_{g,n}$.
\smallskip
\item[\bf Axiom 3 (Cutting edges)]
Given any genus $g$ decorated stable graph $\Gamma$ with $n$ tails
marked with $\bm$, we have a diagram
\begin{equation}
 \begin{diagram} \label{eq:Cutting}
 & &
 \M_{\tilde{\Gamma}}
\times_{\M_{\Gamma}}
 \M^{1/r}_{\Gamma} & \rTo^{\tilde{\mu}} &
 \M^{1/r}_{\Gamma} &
 \rInto^{\tilde{i}} & \M^{1/r}_{g,n} \\
 & \ldTo_{p_1} \\
 \M^{1/r}_{\tilde{\Gamma}} & & \dTo^p & & \dTo^p & &
 \dTo^p \\
 & \rdTo_{p_2}\\
  & & \M_{\tilde{\Gamma}} & \rTo^{\mu} &
 \M_{\Gamma} &
 \rInto^i & \M_{g,n}.
 \end{diagram}
\end{equation}
where $\M_{\tilde{\Gamma}}$ is the stack of stable
 curves with
 graph $\tilde{\Gamma}$, the graph obtained by cutting all edges of
 $\Gamma$, and $\M^{1/r}_{\tilde{\Gamma}}$ is the stack of
 stable $r$-spin curves with graph $\tilde{\Gamma}$ (still marked
 with $m^{\pm}$ on each half edge).
 $p_1$ is the following morphism:
 The fiber product consists of triples of an $r$-spin curve $(X/T,
 \{\mathcal{E}_d, c_{d,d'}\})$, a stable curve $\tilde{X}/T$, and a
 morphism $\nu : \tilde{X} \rTo X$, making $\tilde{X}$ into the
 normalization of  $X$.  Also, the dual graphs of $X$ and $\tilde{X}$
 are $\Gamma$   and $\tilde{\Gamma}$, respectively.  The associated $r$-spin
 curve in $\M^{1/r}_{\tilde{\Gamma}}$ is simply
 $(\tilde{X}/T,   \nu^*\{\mathcal{E}_d,c_{d,d'}\})$.
We require that
\[
p_{1*} \tilde{\mu}^* c^{1/r}_{\Gamma} =  r^{|E(\Gamma)|}
c^{1/r}_{\tilde{\Gamma}} ,\] where $E(\Gamma)$ is the set of edges
of $\Gamma$ that are cut in
 $\tilde{\Gamma}$.
\item[\bf Axiom 4 (Vanishing)]
\label{vanish}
  If $\Gamma$ contains a tail marked with $m_i=r-1$, then
$\cv_{\Gamma}=0.$
\item[\bf Axiom 5 (Forgetting tails)]
\label{forget} Let $\widehat{\Gamma}$ be a stable graph whose
$i$-th tail is marked by $m_i=0$, $\Gamma$ be the stable graph
obtained by removing the $i$-th tail, and
$$\pi:\mbar^{1/r}_{\widehat{\Gamma}} \rTo\mbar^{1/r}_{\Gamma} $$
be the forgetful morphism.  The classes $\cv_{\widehat{\Gamma}}$
and $\pi^*\cv_{\Gamma}$ are related by
$$c^{1/r}_{\widehat{\Gamma}}=\pi^*c^{1/r}_{\Gamma}.$$
\end{description}
\end{df}

\begin{rems} \label{rem:virt}
\
\begin{enumerate}
\item The factor of $\prod_{e}\,(\frac{r}{l_e})$ in Axiom 1
arises from the fact that the right hand square in equation
(\ref{eq:Cutting}) is not quite Cartesian. Rather, because of ramification of
$p$ over $\M_\Gamma$, we have that for any cohomology class $c$ on
$\M_{g,n}^{1/r}$,
\begin{equation}
i^* p_*\,c\,=\, \left(\prod_e \frac{r}{l_e}\right)\,p_*
\tilde{i}^*\,c .
\end{equation}
\item The astute reader will notice in Axiom 3 that, unlike the case of
  a tree, if $\Gamma$ contains a
  loop and if  the $r$-spin structure is Ramond at the corresponding node,
   the dimensions $D_{\Gamma}$ and $D_{\tilde{\Gamma}}$ of the
  virtual classes $\cv_{\Gamma}$ and $\cv_{\tilde{\Gamma}}$ are
  different.  Thus Axiom 3 actually requires the vanishing of both
  $p_{1*} \tilde{\mu}^* \cv_{\Gamma}$ and $\cv_{\tilde{\Gamma}}$ in
  this case.  Of course, for the Ramond case, Axiom 4 already requires
  the vanishing of $\cv_{\tilde{\Gamma}}$, since the cut half-edges
  are both marked with $r-1$.  Thus for any graph (tree or
  otherwise),
  in the Ramond case
  Axiom 3 amounts essentially to the requirement that $p_{1*} \tilde{\mu}^*
  \cv_{\Gamma}$ vanish.
\item
Although the vanishing of $H^0$ or $\pi_*$ is often called
concavity, the Serre dual $\sheafhom (\ce_r,\omega)$,
corresponding to Witten's sheaf $\mathcal{V}$ in \cite{W}, is
convex $(H^1$ vanishes) exactly when $\ce_r$ is concave. Moreover,
$\pi_! \mathcal{V}=\pi_! \ce_r$. Therefore, we use the term convex
to describe  the case when $\pi_* \ce_r =0$.

\item
One might think that the class $c_D(\pi_!\ce_r)$ would be a good
candidate for a virtual class, since it coincides with $\cv$ in
the convex case. However, this is not the case (see Section~\ref{sec:M11}).
\item Witten has described \cite[\S 1.3]{W} an analytic
construction of a class that he calls the ``top Chern class,'' but
it is not  clear  that this class satisfies the above axioms.
Witten's index-like construction is reminiscent of the analytic
construction of a virtual fundamental class of the moduli space of
stable maps in the theory of Gromov-Witten invariants. Ideally,
one should be able to construct $\cv$ by methods similar to those
used in algebraic constructions of the fundamental class.

\item Although, as explained in
Remark~\ref{rem:restrict}, the restriction $0\leq m_i<r$ does not
change the moduli space $\mgnrbar$, it does give a different
choice of $\cv$.  Indeed, replacing $m_i$ by $m_i+r$ changes the
dimension of $\cv$ by $1$ and corresponds (up to a multiplicative
constant) to the first descendant of the classes associated to
$m_i$. Thus on a given moduli space $\mgnrbar$ there are
potentially several (but still only finitely many, for dimensional
reasons) choices of $\cv$ and the corresponding $\cft$. However,
without the restriction $0 \leq m_i < r$, the corresponding metric
$\eta$ is not necessarily invertible, and several other unusual
considerations also arise. These issues will be treated in a
forthcoming paper \cite{JK}. In the remainder of this  paper, we
will assume that $0 \leq m_i <r$ except where explicitly stated
otherwise.

\end{enumerate}
\end{rems}

\subsection{Verification of the \cft\ axioms for $(\ch^{(r)},\eta)$}
\label{sec:verify}
\

In this section we give a proof of Theorem~\ref{th:rcft}, first
for the case where $\bs\,=\,\bu\,=\,\mathbf{0}$, and then in
general.

\subsubsection{The case $\bs\,=\,\bu\,=\,\mathbf{0}$}

\

Let $\cv_{g,n}$ be a cohomology class on $\mgnrbar$ satisfying the
axioms of Section~\ref{axioms}. We will show that the collection
of classes $\{\Lambda_{g,n}\}$ given by~(\ref{eq:rlambda}) satisfies the
\cft\ axioms   C1--C4 with state space $(\ch^{(r)},\eta)$. Axiom C1 clearly
holds by the definition of $\{\Lambda_{g,n}\}$.

Let $\rho=\rho_{tree}$ be the gluing morphism~(\ref{eq:gluetree}).
Condition~(\ref{eq:cfttree}) of Axiom C2 is equivalent to
\begin{align}
\label{eq:axiomtwo}
\rho^*\,p_*\cv_{g,n}(\bm) r^{1-g}\,=&\,\\
\sum_{a,b=0}^{r-2} r^{1-k}\, p_*\cv_{k,j+1}(m_{i_1},\ldots,m_{i_j},a) \tensor
r^{1-(g-k)}\, p_* &\, \cv_{g-k,n-j+1}(b,m_{i_{j+1}},\ldots,m_{i_nj})\eta^{ab}
\nonumber
\end{align}
for all $ 0 \le m_i \le r-1 $ and $\bm\,=\,(m_1,\ldots,m_n)$.
Consider the decorated stable graph $$ \Gamma =\tcongmpm ,$$ and
denote the graph obtained from cutting its edge by
$\tilde{\Gamma}\,=\,\Gamma_{1}\,\sqcup\,\Gamma_{2}$ where $$
\Gamma_{1}\,=\,\tconkm \quad \mathrm{and} \quad
\Gamma_{2}\,=\,\tcongmkm. $$ Since the spaces
$\M^{1/r}_{\Gamma_{i}}$ are non-empty for $0\le m^\pm \le r-1$
only when $m^\pm$ are determined by the
conditions~(\ref{eq:mplus}) and~(\ref{eq:mminus}), the sum in the
right-hand side of equation~(\ref{eq:axiomtwo})  has only one
non-vanishing term. By the definition of the metric
$\eta$~(\ref{eq:metric}), Axiom~C2 reduces to the following:
\begin{equation}
  \label{eq:gammapm}
\rho^* p_*\cv_{g,n}(\bm) r^{1-g} =
\begin{cases}
(p\times p)_*(\cv_{\Gamma_{1}} r^{1-k} \tensor \cv_{\Gamma_{2}}
r^{1-(g-k)} ), & \mathrm{\ if\ } 0\le m^+ \le r-2\\ 0 &\mathrm{\
if\ } m^+=r-1,
\end{cases}
\end{equation}
with $m^-=r-2-m^+$.
In other words, we must show that
\begin{equation}
\rho^* p_*c^{1/r}_{g,n}(\mathbf{m}) =r(p_* c^{1/r}_{\Gamma_1}
\otimes p_* c^{1/r}_{\Gamma_2})=r \cdot p_{2*}
(c^{1/r}_{\tilde{\Gamma}}),
\end{equation}
where we have the diagram
$$
\begin{diagram}
&& \M_{\tilde{\Gamma}} \times_{\M_{\Gamma}} \M^{1/r}_{\Gamma} &
\rTo^{\tilde{\chi}}& \M^{1/r}_{\Gamma} & \rTo^{\tilde{i}}&
\M^{1/r}_{g,n}\\ & \ldTo^{p_1}\\ \M^{1/r}_{\tilde{\Gamma}}& &
\dTo^p & & \dTo^p & & \dTo^p\\ & \rdTo_{p_2} \\ &&
\M_{\tilde{\Gamma}}& \rTo_{\chi} & \M_{\Gamma} & \rTo_i &
\M_{g,n},
\end{diagram} $$
and the map $\rho$ is $ i \circ \chi$. However, if we let $l$ be
$\gcd(m^++1,r)$, then
\begin{eqnarray*}
\rho^* p_* c^{1/r}_{g,n}(\mathbf{m})&=& \chi^* i^* p_*
c^{1/r}_{g,n} \\ & = & \chi^* p_* \tilde{i}^*(r/l) c^{1/r}_{g,n}\\
(\text{by the def. of } c^{1/r}_{\Gamma}) & = &
\chi^*p_*c^{1/r}_{\Gamma}\\ &=& p_* \tilde{\chi}^*
c^{1/r}_{\Gamma}\\ &=& p_{2*} p_{1*} \tilde{\chi}^*
c^{1/r}_{\Gamma}\\ (\text{by Axiom 3}) & =
&(p_{2*}c^{1/r}_{\tilde{\Gamma}}) r.
\\
\end{eqnarray*}
This gives~(\ref{eq:axiomtwo}); therefore, Axiom~C2 is verified.
\medskip
The statement~(\ref{eq:cftloop}) of  Axiom~C3 is equivalent
by~(\ref{eq:rlambda}) and~(\ref{eq:metric})  to
\begin{equation}
\label{eq:axiomthree}
  \rho^* p_*\cv_{g,n}(\bm) \,=\,
  r\,\sum_{m^+=0}^{r-2}\,p_*\,\cv_{g-1,n+2}(m_1,\ldots,m_n,m^+,m^-)~,
\end{equation}
where
$\rho=\rho_{\mathrm{loop}}$ is the gluing morphism~(\ref{eq:glueloop})
and
\begin{equation}\label{eq:mneg}
m^-=\begin{cases}
r-2-m^+ & \mathrm{if} \ 0\le m^+ \le r-2\\
r-1 &  \mathrm{if} \ m^+=r-1 .
\end{cases}
\end{equation}

Let
\begin{equation}
\label{eq:loopgr}
\tilde{\Gamma}_{m^+}=\ocongmcut \quad \mathrm{and} \quad
\Gamma_{m^+}= \ocongmpm
\end{equation}
   be decorated stable   graphs. Let $\tilde{\Gamma}$ and
$\Gamma=\Gamma_\mathrm{loop}$ be the corresponding underlying
(undecorated) graphs, respectively.
We have the commuting diagram
\begin{diagram}
& & \bigsqcup_{ (r-2)/2 \leq m^+ < r} F_{\Gamma, m^+} & \rTo^{\tilde{\chi}} &
\bigsqcup_{ (r-2)/2\leq m^+ <r}
\M^{1/r}_{\Gamma, m^+} & \rTo^{\tilde{i}}&
\M^{1/r}_{g,n}\\
& \ldTo_{p_1} \\
\bigsqcup_{ (r-2)/2 \leq m^+<r}\M^{1/r}_{\tilde{\Gamma}} & & \dTo^p & & \dTo^p & & \dTo^p \\
& \rdTo_{p_2}\\ & &
\M_{\tilde{\Gamma}} &
\rTo_{\chi} &
\M_{\Gamma} & \rTo_i & \M_{g,n}.\\
\end{diagram}
where $F_{\Gamma, m^+} := \M_{\tilde{\Gamma}} \times_{\M_{\Gamma}}
\M^{1/r}_{\Gamma,m^+}$. Therefore, if $l_e$ is defined to be
$\gcd(m^++1,r)$, then
\begin{eqnarray*}
\rho^* p_* c^{1/r}_{g,n} &=&\chi^* p_* \sum_{m^+} (\tilde{i}^*
c_{g,n}) r/l_e\\ &=& \chi^* p_* \sum_{m^+} c^{1/r}_{\Gamma,
m^+}\\ &=& p_* \tilde{\chi}^* \sum_{m^+} c^{1/r}_{\Gamma,m^+}\\ &=&
p_{2*}p_{1*} \tilde{\chi}^* \sum_{m^+} c^{1/r}_{\Gamma, m^+}\\ &=&
  r \cdot p_{2*} \sum_{m^+}
c^{1/r}_{\tilde{\Gamma}, m^+}.
\end{eqnarray*}
This proves axiom~C3.
To prove Axiom~C4a, consider the Cartesian
square
$$
\begin{diagram}
\mbar^{1/r,(m_1,\ldots,m_n,0)}_{g,n+1} & \rTo^{\tilde{\pi}} &
\mbar_{g,n}^{1/r,\mathbf{m}} \\ \dTo_p & & \dTo_p\\ \mbar_{g,n+1}
& \rTo^\pi & \mbar_{g,n}.
\end{diagram}
$$ By~(\ref{eq:lambda}) and Axiom~5 (forgetting tails) we have the
required equation~(\ref{eq:identity})
\begin{eqnarray*}
\pi^*\Lambda_{g,n}(e_{m_1},\ldots,e_{m_n})
&=&r^{1-g}\,\pi^*\,p_*\cv_{g,n}(\mathbf{m})
\,=\,r^{1-g}\,\,p_*\tilde{\pi}^*\cv_{g,n}(\mathbf{m})\\ &=&
r^{1-g}\,\cv_{g,n+1}(m_1,\ldots,m_n,0)
=\Lambda_{g,n+1}(e_{m_1},\ldots,e_{m_n},e_0) .
\end{eqnarray*}
Finally, a direct calculation yields Axiom 4b (see Proposition
\ref{prop:potential}).

\subsubsection{The general case} \label{sec:generalCohFT}

\

The proofs of Axioms C1 and C4a remain the same. We only need to
prove axioms C2 and C3.  Before doing so, we will need a lemma on
regular imbeddings and base change.

 \begin{lm} \label{lm:imbeddings}
 Let $i$ be a regular imbedding, and let
 \begin{diagram}
 Y' & \rTo_{\tilde{i}} & X' \\
 \dTo^{\pi'} & & \dTo^{\pi}\\
 Y & \rTo^i & X
 \end{diagram}
  be a Cartesian square with $\pi$ and $\pi'$ proper and
flat, and $\mathcal{E}$ a coherent sheaf on $X'$,  flat over $X$.
 If $X$ and $Y$ both carry an ample invertible sheaf, then the
 Chern character commutes with base change,
that is $i^*ch\pi_! \ce = ch \pi'_! \tilde{i}^* \ce$.
 \end{lm}

\begin{proof} First, we claim that because $i$ is a
regular imbedding it has finite Tor dimension; that is, there is
an integer $N$ such that for every coherent $\co_X$-module $\cf$,
the $\co_Y$-modules $\mathcal{T}or^{\co_X}_j(\co_Y,\cf)$ vanish
for $j>N$.  This can be seen as follows.  We may assume that $X$
is $\operatorname{Spec} A$ and $Y$ is $\operatorname{Spec} A/(x)$
for some regular element $x$ in a ring $A$. This gives the free
resolution $0 \rTo A \rTo^{(x)} A \rTo A/(x) \rTo 0$ of $\co_Y$,
and shows $i$ has finite Tor dimension.

Since $\pi$ is flat, the sheaves $\co_{X'}$ and $\co_Y$ are Tor
independent over $X$; that is, $Tor^{\co_X}_j(\co_{X'},\co_Y)=0$
for all $j>0$.

Let $Li^*$ and $L\tilde{i}^*$ be the left derived functors of
$i^*$ and $\tilde{i}^*$, respectively. Proposition 5.13 of
\cite{Srinivas} states that if $i$ has finite Tor dimension, and
if $\co_{X'}$ and $\co_Y$ are Tor-independent over $X$, then $Li^*
\pi_! =\pi'_! L\tilde{i}^*$.

However,  since $\ce$ is flat over $X$, we have $L\tilde{i}^*\ce:=
\sum^\infty_{j=0} (-1)^j Tor^X_j (\ce,\co_Y)= \ce\otimes_{\co_X}
\co_Y =\tilde{i}^*\ce_j$; and since $Li^*$ commutes with the Chern
character, the lemma holds.
\end{proof}

Now we prove that Axioms C2 and C3 hold.

First consider the Cartesian square
\begin{diagram}
\mathcal{C}_{\Gamma} & \rTo^{\hat{i}} & \mathcal{C}\\ \dTo^{\pi} &
& \dTo^{\pi}\\ \M^{1/r}_{\Gamma} & \rTo_{\tilde{i}} &
\M^{1/r}_{g,n},
\end{diagram}
where $\pi$ is the universal curve.  Let $\mathcal{E}_r$ be the
$r$-th root from the universal $r$-spin structure on
$\mathcal{C}$. The morphism  $\pi$ is projective, and
$\M^{1/r}_{g,n}$ carries an ample invertible line bundle
\cite[3.1.1]{J},  so Lemma \ref{lm:imbeddings} gives
$\mathrm{ch}_i(\pi_! \hat{i}^* \mathcal{E}_r)\, =\, \mathrm{ch}_i(
\tilde{i}^* \pi_!\mathcal{E}_r)$.   Note that $\hat{i}^*
\mathcal{E}_r$ is the $r$-th root from the universal  $r$-spin
structure $\mathcal{C}_{\Gamma}$.

Let $F$ be the fiber product
\begin{equation} \label{eq:FDef}
F:= \M_{\tilde{\Gamma}} \times_{\M_{\Gamma}} \M^{1/r}_{\Gamma}.
\end{equation}
We have
\begin{diagram}
\mathcal{C}_{\tilde{\Gamma}} & \lTo & \hat{\mathcal{C}}_F \\ \dTo
& & \dTo^{\theta}\\ & & \mathcal{C}_F & \rTo &
\mathcal{C}_{\Gamma}\\ && j\dTo && \dTo\\
\M^{1/r}_{\tilde{\Gamma}}& \lTo^{p_1} & F & \rTo &
\M^{1/r}_{\Gamma},\\
\end{diagram}
where $\mathcal{C}_{\tilde{\Gamma}}$ is defined to be the
universal curve over $\M^{1/r}_{\tilde{\Gamma}}$, and
$\hat{\mathcal{C}}_F$ and $\mathcal{C}_F$ are the fibered products
$F \times_{\M^{1/r}_{\tilde{\Gamma}}}\mathcal{C}_{\tilde{\Gamma}}$
and $F \times_{\M^{1/r}_{\Gamma}} \mathcal{C}_{\Gamma}$,
respectively.  The morphism $\theta$ is just normalization, and in
fact, if $\hat{\mathcal{C}}_{\Gamma} \rTo^{\theta_{\Gamma}}
\mathcal{C}_{\Gamma}$ is the normalization of
$\mathcal{C}_{\Gamma}$, then we have the following fibered diagram
(all rectangles are Cartesian).\footnote{Elsewhere in the paper
the maps $\tilde{\chi}$ and $\theta$ were called $\tilde{\mu}$ and
$\nu$, respectively, but we have renamed them here in order to
avoid confusion with the cohomology classes $\mu_i$ and $\nu_i$.}

\begin{diagram}
\mathcal{C}_{\tilde{\Gamma}}& \lTo_{\hat{p}_1} &
\hat{\mathcal{C}}_F & \rTo_{\hat{\chi}} &
\hat{\mathcal{C}}_{\Gamma}\\ \dTo^{\pi} & & \dTo^{\theta} & &
\dTo^{\theta_{\Gamma}}\\ & & \mathcal{C}_F & \rTo &
\mathcal{C}_{\Gamma}\\ & & \dTo^{\pi} & & \dTo^{\pi}\\
\M^{1/r}_{\tilde{\Gamma}} & \lTo^{p_1} & F & \rTo^{\tilde{\chi}} &
\M^{1/r}_{\Gamma}\\
\end{diagram}

Moreover, if $\mathcal{E}_{\tilde{\Gamma}}$,
$\mathcal{E}_{\Gamma}$, $\hat{\mathcal{E}}_{\Gamma}$,
$\mathcal{E}_F$, and $\hat{\mathcal{E}}_F$ are the $r$-th roots of
the universal  $r$-spin structures on
$\mathcal{C}_{\tilde{\Gamma}}$, $\mathcal{C}_{\Gamma}$,
$\hat{\mathcal{C}}_{\Gamma}$, $\mathcal{C}_F$, and
$\hat{\mathcal{C}}_F$, respectively, we have $\hat{\mathcal{E}}_F
= \hat{p}_1^* \mathcal{E}_{\tilde{\Gamma}}$. Also, since $\theta$
is finite, $ (\pi \circ \theta)_! \hat{\mathcal{E}}_F = \pi_!
\mathcal{E}_F$.  Since $\tilde{\chi}$ is flat, we have
$\tilde{\chi}^* \pi_! \mathcal{E}_{\Gamma} = \pi_! \tilde{\chi}^*
\mathcal{E}_{\Gamma}= \pi_! \mathcal{E}_F = (\pi \circ \theta)_!
\hat{\mathcal{E}}_F$. Also, since $p_1$ is flat $p_1^* \pi_!
\mathcal{E}_{\tilde{\Gamma}}= (\pi \circ \theta)_! p_1^*
\mathcal{E}_{\tilde{\Gamma}}= (\pi \circ \theta)_!
\hat{\mathcal{E}}_F$. Furthermore,
$p_1^*\,\pi_!\,=\,(\pi\circ\theta)_!\,p_1^*$. If
$\tilde{i}\,:\,\M_\Gamma\,\rTo\,\M^{1/r,\bm}_{g,n}$ is the
inclusion map, then the above implies that
\begin{eqnarray} \label{eq:muformula}
&\tilde{\chi}^*\, \tilde{i}^* \mathrm{ch}_i(\pi_!\,\ce)\,
=\,\tilde{\chi}^*\,\mathrm{ch}_i(\pi_!\,\ce_\Gamma)\,=\,
\mathrm{ch}_i(\pi_!\, \hat{\chi}^*\,\ce_\Gamma)  \,=\, \\
&\mathrm{ch}_i((\pi\circ\theta)_!\,\hat{p}_1^*\,\ce_{\tilde{\Gamma}})
\,=\, \mathrm{ch}_i(p_1^*\,\pi_!\,\ce_{\tilde{\Gamma}}) \,=\,
p_1^*\,\mathrm{ch}_i(\pi_!\,\ce_{\tilde{\Gamma}}).  \nonumber
\end{eqnarray}
The previous equation, the fact that  $\nu_i$ on
$\M_{g,n}^{1/r,\bm}$ is the lift of $\nu_i$ on $\M_{g,n}$, and the
projection formula yield the desired result.

This finishes the proof of Theorem~\ref{th:rcft}.

\subsection{The genus-zero case}

\

As we explained in Remark~\ref{rem:virt}.3, we call the $r$-spin
structure   $(\{\ce_d\},\{ c_{d,d'}\})$ on the universal curve
$\cc^{1/r}_{\Gamma}$ \emph{convex} if $\pi_* \ce_r$ is identically
zero.  This occurs, for example, when $g=0$, as is shown in the
following proposition.
\begin{prop}\label{this2}
  Let $X$ be a prestable curve of genus zero with $n$ punctures and
  markings $(m_1,\ldots,m_n)$, such that $-1 \le m_i \le r-1 $ for all
  $i$ and $m_i \ge 0$ for all $i$ except at most one.  Then, if
  $(\ce,b)$ is an $r$-th root of $\omega_X(-\sum m_i p_i)$, we have $
  H^0 (X, \ce) =0.  $
\end{prop}
\begin{proof}
The degree of $\ce$ is an integer and is equal to $-(2+\sum
m_i)/r$. Thus $\sum m_i \geq r-2$, and the degree of $\ce$ is
strictly negative.  Therefore, when $X$ is irreducible, $\ce$ has
no global sections.  When $X$ is not irreducible, but $\ce$ is
locally free (Ramond) at each node, the same argument holds. If
$\ce$ is Neveu-Schwarz at some nodes, then normalization $\nu:
\tilde{X} \rTo X$ at the nodes of $X$ where $\ce$ is not locally
free gives $\ce = \nu_* \cf$, where $\cf$ is locally free on
$\tilde{X}$. Restricting $\cf$  to $\tilde{X}$ we obtain  an
$r$-th  root of $\omega_{\tilde{X}}(-\sum\tilde{m}_i
\tilde{p}_i)$, where the points $\tilde{p}_i$ are either marked
points or inverse images of nodes, and thus the collection
$\tilde{m}_i$ still meets the hypotheses of the proposition, but
now $\cf$ is locally free on each  component, and hence has no
global sections.  Since $\nu$ is  finite, $H^0 (X, \ce) = H^0
(\tilde{X}, \cf) =0$.
\end{proof}

\medskip
The previous proposition shows that if a class $\cv$ on
$\M^{1/r}_{0,n}$ satisfying Axioms~1--5 exists, then by Axiom~2 it
must be the top Chern class of the bundle with fiber
$H^1(X,\ce_r)^*$ at $[(X,p_1,\dots,p_n,(\{\ce_d\},\{c_{d,d'}\}))]
\in \M^{1/r}_{0,n}$. In this case, it does indeed satisfy the
required properties.
\begin{thm}
\label{thm:gzero}
Define cohomology classes on $\M^{1/r}_{0,n}$ by
\begin{equation}
  \label{eq:genzero}
\cv_{0,n}(\bm)=c_{D} (\pi_! \ce_r) = (-1)^{D} c_{D} (R^1 \pi_* \ce_r),
\end{equation}
where $$ D = \frac{1}{r}(2-r+\sum_{i=1}^n m_i) $$ and $\ce_r$ is
the $r$-th root sheaf of the universal $r$-spin structure. Then
the collection of classes $\cv_{\Gamma}$ defined
by~(\ref{eq:cvgamma}) for decorated stable graphs of genus zero
satisfies Axioms~1--5.
\end{thm}

\begin{proof}
It is clear from the construction of the classes $\cv_{\Gamma}$
that they satisfy Axiom~2 (convexity).

Axiom 1 follows from the fact that since $\ce_{r,\Gamma}= \ce_{r,
\Gamma_1} \oplus \ce_{r, \Gamma_2} $, the top-dimensional Chern
class $c_{D(\Gamma)}=c_{top}$ of $\pi_! \ce_{r,\Gamma}$, is simply
the product of the top-dimensional classes of $\pi_!
\ce_{r,\Gamma_1}$ and $\pi_! \ce_{r,\Gamma_2}$.

Now we will show that Axiom 3 holds.  If $\ce_r$ is the $r$-th
root from the universal $r$-spin structure on
$\cgn^{1/r}\,\to\,\M_{g,n}^{1/r}$, then since $g=0$, $\ce_r$ is
convex.  Repeating the argument in the proof of Axiom~3 in
Section~\ref{sec:generalCohFT} with the Chern character replaced
by the top Chern class yields $$
\cv_\Gamma\,=\,(r^{|E|}/\prod_{e\,\in\,E}
l_e)\,c_D(-R^1\pi_*\ce_r), $$ where $E$ denotes the edge set of
$\Gamma$. Also, we can compute the degree of $p_1:F \rTo
\M^{1/r}_{\tilde{\Gamma}}$ from the diagram
\begin{diagram}
& & F & \rTo^{\tilde{\chi}} & \M^{1/r}_{\Gamma}\\ & \ldTo^{p_1} \\
\M^{1/r}_{\tilde{\Gamma}} & & \dTo^{\tilde{p}} & &
\dTo^{p_{\Gamma}}\\ & \rdTo^{p_2} \\ & & \M_{\tilde{\Gamma}} &
\rTo & \M_{\Gamma}.
\end{diagram}
The morphism $\M_{\tilde{\Gamma}} \rTo \M_{\Gamma}$ has degree
$1$, since $\Gamma$ is a tree, and therefore $\tilde{\chi}$ also
has degree $1$. The morphism $p_{\Gamma}$ has degree $\prod_e l_e
r^{-|V(\Gamma)|}$, as can be seen from the fact that the coarse
moduli space map induced by $p_{\Gamma}$ has degree $1$, and there
are $r^{|V(\Gamma)|}/{\prod_{e} l_e}$ automorphisms of a generic
$r$-spin structure.  Thus the map $\tilde{p}$ also  has degree
$\prod_e l_e r^{-|V(\Gamma)|}$, the map $p_2$ has degree
$r^{-|V(\Gamma)|}$, so $p_1$ has degree exactly $\prod_e l_e$. Now
we can compute $p_{1 *} \tilde{\chi}^*c^{1/r}_{\Gamma}$ using
equation~(\ref{eq:muformula}) (replacing $\mathrm{ch}_i$ with
$c_D$) to obtain
\begin{eqnarray*}
p_{1 *}(r^{|E|}/\prod_e l_e)\tilde{\chi}^* c_D (-R^1 \pi_*
\mathcal{E}_{\Gamma}) & = & (r^{|E|}/{\prod_e l_e}) p_{1 *} c_D(p_1^*
(- R^1 \pi_* \mathcal{E}_{\tilde{\Gamma}}))\\
&=& (r^{|E|}/{\prod_e l_e}) p_{1 *} p_1^* c_D (-R^1 \pi_*
\mathcal{E}_{\tilde{\Gamma}})\\
&=& r^{|E|} c_D (-R^1 \pi_* \mathcal{E}_{\tilde{\Gamma}}) \\ &=&
r^{|E|}c^{1/r}_{\tilde{\Gamma}}
\end{eqnarray*}
as desired. All that remains to check is Axiom 4 (vanishing). Let
$p$ be a point corresponding to a tail marked  by $m=r-1$. Taking
the tensor product of $\ce_r$ with the exact sequence $$ 0 \rTo
\co \rTo \co(p) \rTo \co(p)|_p \rTo 0 $$ gives the exact sequence
$$ 0 \rTo \ce_r \rTo \ce' \rTo \ce'|_p \rTo 0, $$ where
$\ce'=\ce_r \otimes \co(p)$.  Thus $\ce'$ corresponds to a root
with $p_i$ marked by $m'=-1$.  Since $R^1\pi_* (\ce'|_{p}) =0$,
and since the residue isomorphism~(\ref{eq:residue}) $R_p : \pi_*
(\ce'|_p) \irightarrow \co$ shows that $\pi_* (\ce'|_p)$ is a
trivial bundle,  we have $\pi_! (\ce') = \co+\pi_!(\ce_r).$  By
Proposition \ref{this2} the sheaves $\ce_r$ and $\ce'$ are both
convex; thus $\pi_!\ce_r=-R^1 \pi_* \ce_r$ and $\pi_! \ce' = -R^1
\pi_* \ce'$ are both locally free and have the same Chern classes
in all dimensions.  However,  the vector bundle $\pi_! \ce'$ has
dimension $D'=D-1$, and so $\cv=c_D (\pi_! \ce_r)=0$. This gives
Axiom 4.
\end{proof}

\subsection{The case $\M_{1,1}^{1/r,\mathbf{0}}$} \label{sec:M11}

\

In this section we will calculate the virtual class on
$\M_{1,1}^{1/r,\mathbf{0}}$ from the axioms. Let $\Gamma_{m^+}$ be
the decorated graph as in~(\ref{eq:loopgr}) with the underlying
graph $\Gamma =\Gamma_{loop}$ with one tail, one node  of genus
zero, and one loop whose one half-edge is marked with $m^+$ and
the other with $m^-$ given by~(\ref{eq:mneg}).

The stack $\M^{1/r, \mathbf{0}}_{1,1}$ is a disjoint union
$\coprod_{d|r}  \M^{1/r, \mathbf{0},(d)}_{1,1}$, where the
component indexed by $d$ has a generic geometric point
corresponding to a smooth $r$-spin curve $(X, \ce_r, c_{r,1})$
with $\ce_r^{\otimes d}$ isomorphic to $ \mathcal{O}_X$.  That is,
$\ce_r$ is a $d$-torsion point of the Jacobian of $X$. Since
$\ce_r$ has global sections if and only if $d=1$, the case of
$d>1$ is convex. Since the dimension~(\ref{eq:deg}) of the virtual
class $c_{1,1}^{1/r,\mathbf{0}}$  is $0$, we have that
$c_{1,1}^{1/r,\mathbf{0},(d)}=1$ for $d>1$. Moreover, consider the
graph $\Gamma_{r-1}$.  By Axiom 3, $p_{1 *}\tilde{\mu}^*
c_{{\Gamma}_{r-1}} = c_{\widetilde{\Gamma}_{r-1}}$, where
$\widetilde{\Gamma}_{r-1}$ is $\Gamma_{r-1}$ with the loop cut.
Since  both the half-edges of the cut loop will be labeled by
$r-1$, by Axiom 4 the corresponding class must vanish. Therefore,
the (Ramond) case of $\Gamma_{r-1}$ with a trivial gluing (i.e.,
$\ce_r =\co_X$) yields
\begin{equation} \label{eq:M11}
 c^{1/r,\mathbf{0},(1)}_{\Gamma_{r-1}}\,=\,i^* c^{1/r,\mathbf{0},(1)}_{1,1} =-(r-1),
\end{equation}
since  all the remaining Ramond  components have $i^*
c^{1/r,\mathbf{0},(d)}_{1,1}=+1$, and there are $r-1$ of them.
Since $D=0$, this means $c^{1/r,\mathbf{0},(1)}_{1,1}$ is also
equal to $-(r-1)$. \emph{Notice that this differs from the top
Chern class of the bundle
$c_D(\pi_!\,\ce_r)\,=\,c_0(\pi_!\ce_r)\,=\,1$.}

Now, the map $p^{(1)}: \M^{1/r,\mathbf{0},(1)}_{1,1} \rightarrow
\M_{1,1}$ has degree $1/r$, and
$p^{(d)}:\M^{1/r,\mathbf{0},(d)}_{1,1} \rightarrow \M_{1,1}$ has
degree $d^2/r \prod_{p|d} (1-1/p^2)$. The latter is $1/r$ times
the number of points of order precisely $d$ on the Jacobian of the
underlying curve. Thus
\begin{align*}
p_* c^{1/r}_{1,1} &=1/r \sum_{d>1} d^2
\prod_{p|d} (1-1/p^2)-(r-1)/r \\ &= \frac{r^2-1}{r}-\frac{r-1}{r}
=r-1.
\end{align*}
Therefore,
\begin{align*}
\cor{\tau_{1,0}}_1\,&=\,\int_{\M_{1,1}^{1/r,0}}\,\psi_1 c_{1,1}^{1/r,0}\,=\,
\int_{\M_{1,1}}\,\psi_1 p_* c_{1,1}^{1/r,0}\,=\,
(r-1) \int_{\M_{1,1}}\,\psi_1,
\end{align*}
and we conclude that
\begin{equation}
\label{eq:M11Calc}
\cor{\tau_{1,0}}_1 = \frac{r-1}{24}.
\end{equation}

Equation (\ref{eq:M11Calc}) is consistent with the prediction from \cft\
stated in equation (\ref{eq:M11CohFT}).

\subsection{The case $r=2$} \label{sec:rtwo}

\

In this section we will show that in the case of
theta-characteristics (i.e., when $r=2$) there exists a unique
virtual class $\cvtwo$ satisfying the axioms of
Section~\ref{axioms}.

\begin{thm}\label{th:rtwo}
The collection of cohomology classes $$ \cvtwo_{g,n}(\bm) \in
H^\bullet (\mbar^{1/2,\bm}_{g,n}) $$ satisfies Axioms~1--5 of
Section~\ref{axioms} if and only if $\cvtwo_{g,n}(\bm)=0$ for $\bm
\ne \mathbf{0}$, and for $\bm = \mathbf{0}$ the class
$\cvtwo_{g,n}(\mathbf{0})$ belongs to
$H^0(\M^{1/2,\mathbf{0}}_{g,n})$ and is given by
\begin{equation}
  \label{eq:cvrtwo}
\cvtwo_{g,n}(\mathbf{0})=
\begin{cases}
\ \ 1 & \mathrm{\ on\ \ } \mgnbar^{1/2,\mathbf{0},\mathrm{even}}
\\ -1 & \mathrm{\ on\ \ } \mgnbar^{1/2,\mathbf{0},\mathrm{odd}}.
\end{cases}
\end{equation}
\end{thm}

\begin{proof}
Let us show first that the conditions of the theorem are
necessary. If $\bm\ne {\bf 0}$, then by Axiom~4,
$\cvtwo_{g,n}(\bm)$ must vanish; therefore, we can assume that
$\bm={\bf 0}$. In this case the dimension $D$ (given in
equation~(\ref{eq:deg})) of the class $\cvtwo_{g,n}$ is equal to
zero,  and since $\M^{1/2,{\bf 0}}_{g,n}$ has two connected
components ($\mgnbar^{1/2,\mathbf{0},\mathrm{even}}$ and
$\mgnbar^{1/2,\mathbf{0},\mathrm{odd}})$ it will be sufficient to
find $\cvtwo_{\Gamma}$ for two  graphs $\Gamma_0$ and $\Gamma_1$,
such that the intersections $\M^{1/2}_{\Gamma_0} \cap
\mgnbar^{1/2,\mathbf{0},\mathrm{even}}$ and $\M^{1/2}_{\Gamma_1}
\cap \mgnbar^{1/2,\mathbf{0},\mathrm{odd}}$ are non-empty.

Let $\Gamma_0$ be the graph with one genus-zero vertex, $n$ tails,
and  $g$ Neveu-Schwarz (i.e., all half-edges are decorated with
zeroes) loops. In this case,
$\M^{1/2}_{\tilde{\Gamma_0}}\,=\,\M_{\tilde{\Gamma_0}}\,
\times_{\M_{\Gamma_0}}\,\M^{1/2}_{\Gamma_0}$
in~(\ref{eq:Cutting}), so by Axiom~3 (cutting edges) the class
$\cvtwo_{\Gamma_0}$ pulls back to $\cvtwo_{\tilde{\Gamma}_0}$,
where $\tilde{\Gamma}_0$ is the graph with one vertex of genus
zero and $n+2g$ tails. Since the genus is zero, the universal
square root $\ce$ of $\omega$ on the universal curve over
$\M^{1/2,\mathbf{0}}_{0,n+2g}$ is convex by Proposition
\ref{this2}. Therefore, if $E$ is the set of edges of $\Gamma_0$,
we have
\[
c_{g,n}^{1/2,\mathrm{even}}(\mathbf{0})= \frac{1}{2^{|E|}}
c_{\Gamma_0}^{1/2,\mathrm{even}}(\mathbf{0}) =
\cvtwo_{\tilde{\Gamma_0}} = c_0(-R^1 \pi_* \ce)=1,
\]
where the first equality follows from Axiom 1, the second from
Axiom 3, and the third from Axiom 2. To find
$c^{1/2,\mathrm{odd}}$ consider the graph $\Gamma_1$ with a single
vertex of genus one, $n$ tails, and $g-1$ Neveu-Schwarz loops.
Axiom~3 again shows that  $\cvtwo_{\Gamma_1}$ pulls back to
$\cvtwo_{\tilde{\Gamma}_1}$, where $\tilde{\Gamma}_1$ has a single
vertex of genus one and $2g-2+n$ tails. Since $\bm={\bf 0}$,
Axiom~5 (forgetting tails) shows that $c^{1/2,\mathrm{odd}}$ is a
pullback from  $\M_{1,1}^{1/2,\mathbf{0},\mathrm{odd}} \,=\,
\M_{1,1}^{1/2,(1)}$, and $c_{1,1}^{1/2,(1)}\,=\,-1$ by equation
(\ref{eq:M11}).

Now let us show that the classes $\cvtwo_{g,n}(\bm)$ defined above
for $r=2$ indeed satisfy Axioms~1--5.

Axiom~2 (convexity) holds when $\bm=\mathbf{0}$, since in this
case the class has dimension $0$, and if $\ce$ is convex (and,
therefore, even) $R^1\pi_*\ce=0$ and $\cvtwo = c_{top}=1$ as
required.

If $\bm \ne \mathbf{0}$ then $\ce$ is not convex on the universal
curve over $\M^{1/2,\bm}_{g,n}$ for any $g>0$.  In particular,
consider the degenerate curve of genus $g$ which has two
irreducible components $E$ and $C$ joined at a single node, the
component $E$ of genus zero, containing all $n$ marked points, and
the component $C$ of genus $g$. For degree reasons the node must
be Neveu-Schwarz, and so $\ce$ corresponds to $\ce_E \oplus
\ce_C$, for $\ce_E$ a square root of type $\bm$ of $ \omega_E =
\co_E (-2-\sum m_i)$, and for $\ce_C$ a theta-characteristic on
$C$.  In general, $\ce_C$ has non-zero global sections.  Thus
$\ce$ also has non-zero global sections, and the universal  square
root is not convex.

When $g=0$ and $\bm \neq 0$, the sheaf $\ce$ is convex by
Proposition \ref{this2}. By Theorem \ref{thm:gzero} (and Axiom 4)
the class $c_D(-R^1 \pi_*\ce)$ vanishes, and so agrees with our
definition of $\cvtwo_{g,n}(\bm)$. Axiom~1 holds because of the
simple observation that the parity of a root $\ce$ over a curve
with the graph $\Gamma_1 \sqcup \Gamma_2$ is equal to the sum
modulo $2$ of the parities of the restrictions of $\ce$ to the
components corresponding  to $\Gamma_1$ and  $\Gamma_2$.

To prove Axiom~3, we may assume that $\bm\,=\,\mathbf{0}$, and it
is sufficient to check the case that $\Gamma$ has only one edge.
We have by definition $c^{1/2,even}_{\Gamma} =2$ and
$c^{1/2,odd}_{\Gamma} = -2$ for $\Gamma$ a tree.  Let $F$ be
defined as in equation (\ref{eq:FDef}). The canonical morphism
$p_1:F \rTo \M^{1/2}_{\tilde{\Gamma}}$ is actually an  isomorphism
if $\Gamma$ is a tree, and so we get $p_{1 *} \tilde{\mu}^*
c^{1/2}_{\Gamma} =2 \cdot c^{1/2}_{\tilde{\Gamma}}$ since the
parity of $\mathcal{E}_2$ does not change when restricting to the
normalization. In the case of a loop, there are two subcases.
First when $m^+\,=\,1$, $F$ is isomorphic to two copies of
$\M^{1/2}_{\tilde{\Gamma}}$ (because of the  two choices of gluing
data---see Section \ref{sec:NSRnodes}). Consequently, $p_1$ has
degree $2$ if $m^+  =1$. In the second case $m^+\,=\,0$, and  $F$
is isomorphic to $\M^{1/2}_{\tilde{\Gamma}}$, so $p_1$ has degree
$1$.  Also,  $c^{1/2}_{\Gamma}$ is $2   \tilde{i}^* c^{1/2}_{g,n}$
if $m^+ =0$, and $\tilde{i}^* c^{1/2}_{g,n}$ if $m^+=1$. Thus,
when $m^+=0$, $p_{1 *} \tilde{\mu}^*
c^{1/2}_{\Gamma}=2c^{1/2}_{\tilde{\Gamma}}$ as desired.

When $m^+=1$,   the two choices of gluing give different parities.
Since parity is deformation invariant, this can be seen by
degenerating to the special case of the curve $X$, whose partial
normalization $\tilde{X}$ at one node $q$ consists of two
irreducible components joined at a single node $q$.  One component
$C$ is of genus $g-1$ and contains the marked points $p_1, \dots,
p_{n}$.  The other component $\tilde{E}$ is of genus zero and
contains marked points $q^+$ and $q^-$.  Degree reasons force the
node to be Neveu-Schwarz, and so $\tilde{\ce}$ is simply a direct
sum $\ce_C \oplus \ce_{\tilde{E}}$.  Moreover, since
$\tilde{\ce}_{\tilde{E}}$ is a square root of $\omega_{\tilde{E}}
= \co_{\tilde{E}} (-2)$ of type $(0,-1,-1)$, $\ce_{\tilde{E}}$
must be trivial ($\ce_{\tilde{E}} \cong \co_{\tilde{E}}$).  Gluing
$\ce_{\tilde{E}}$ via $+1$ and $-1$ yields the trivial bundle
$\co_E = \ce^+_E$ and another non-trivial bundle $\ce^-_E$ of
degree zero, respectively. Consequently, $h^0 (\ce^+)=h^0(\ce_C
\oplus \co) =1+h^0(\ce_C \oplus \ce^-_E)$. Since the parities of
$\ce^+$ and $\ce^-$ are simply the parities of $h^0(\ce^+)$ and
$h^0(\ce^-)$, respectively, $\ce^+$ and $\ce^-$ have different
parities. The different parities under the two gluings give $p_{1
*} \tilde{\mu}^* c^{1/2}_{\Gamma} =0 =c^{1/2}_{\tilde{\Gamma}}$,
since $\tilde{\Gamma}$ has a tail marked with $+1=m^+$. This
completes the proof of Axiom~3.

Axiom~5 is true because the projection
$\mbar^{1/2,(m_1,\ldots,m_n,0)}_{g,n+1} \rTo
\mbar_{g,n}^{1/2,\mathbf{m}}$ respects the parity of components.
Axiom~4 follows from the definition of the classes $\cvtwo$.
\end{proof}

This theorem together with Theorem~\ref{th:rcft} implies that in the $r=2$
case we obtain a well-defined complete \cft\ of rank one. It turns out to be
the same  as the Witten-Kontsevich rank-one \cft\ of the pure topological
gravity. Namely, we have the following result (cf.\ also~\cite{W}).

\begin{crl}\label{cor:rtwo}
The class $\Lambda_{g,n}(e_{m_1},\ldots,e_{m_n})\in
H^\bullet(\mgnbar)$ of the $2$-spin \cft\ given
by~(\ref{eq:rlambda}) is equal to $1$ if $m_1=m_2=\ldots=m_n=0$
and $0$ otherwise.
\end{crl}
\begin{proof}
We only need to check the case $m_1=m_2=\ldots=m_n=0$. Since the class has
dimension $0$, equation~(\ref{eq:rlambda}) gives
$$
\Lambda_{g,n}(e_0,\ldots,e_0)=2^{1-g}p_*c^{1/2}_{g,n}(\mathbf{0})=
2^{1-g}\left(2^{g-1}(2^g+1)-2^{g-1}(2^g-1)\right)/2 = 1,
$$
where $2^{g-1}(2^g\pm 1)$ are the numbers of even/odd theta characteristics
on a smooth curve of genus $g$ and the last factor of $1/2$ is the local
(orbifold) degree of the map $p$ near a generic point of
$\mbar^{1/2,\mathbf{0}}_{g,n}$.
\end{proof}

In Section~\ref{sec:kdvr} we will see that this corollary together
with Kontsevich's theorem gives the generalized Witten conjecture
for $r=2$.

\section{Intersection numbers and recursion relations} \label{sec:rr}

In this section, we use relations between boundary classes and
tautological classes in order to derive recursion relations
between intersection numbers on the moduli space of stable
$r$-spin curves. Throughout this section, we assume the existence
of    a virtual class $\cv$ satisfying Axioms 1--5 of
Section~\ref{axioms}. This class was shown  in Section~\ref{vc} to
exist in genus zero for arbitrary $r$, and in arbitrary genus for
$r=2$. Let $(\chr,\eta,\Lambda,e_0)$ be the $r$-spin \cft\ with
the standard basis $\{\,e_0\,\ldots, e_{r-2}\,\}$ described in
Corollary \ref{crl:rCFT}.

We will find it convenient to introduce the notion of a very large
phase space potential in this section. For all $v_1,\ldots, v_n$
in $\chr$, define
\[
\ccor{\tau_{a_1}(v_1)\,\cdots\,\tau_{a_n}(v_n)}_g\,:=\,
\cor{\tau_{a_1}(v_1)\,\cdots\,\tau_{a_n}(v_n)\,
\exp({\bt\cdot\btau + \bs \cdot  \bmu  + \bu \cdot \bnu})\,}_g,
\]
where
\[
\bt\,\cdot\,\btau\,=\,\sum_{\substack{0\,\leq\,m\,\leq\,r-2 \\
  a\,\geq\,0}}\,\tau_{a,m}\,t_a^m
\]
and $\tau_{a,m}\,=\,\tau_a(e_m)$.
Here
$$
\bu \cdot \bnu\,=\,\sum_{i\,\geq\,1}\, u_i\,
\nu_i \qquad \bs\cdot\bmu\,=\,\sum_{i\,\geq\,1}\,s_i\mu_i,
$$
where the classes $\nu_i$ and $\mu_i$ are defined by~(\ref{eq:nu})
and~(\ref{eq:mu})  as components of the Chern characters
of the Hodge bundle and its $r$-spin analog.
These expressions should be understood as formal power series in variables
$t_a^m, u_i,  v_i$.
The correlators are defined by
\begin{align*}
&
\cor{\tau_{a_1}(e_{m_1})\,\cdots\,\tau_{a_n}(e_{m_n})\,}_g
\,:=\,\\
&\lambda^{2g-2}\,\int_{\mgnrmbar}\,
  \psi_1^{a_1}\,\cdots\,\psi_n^{a_n}\,\cv_{g,n}(\bm)\, \exp({\bs \cdot \bmu +
\bu   \cdot   \bnu}),
\end{align*}
where $\mathbf{m}\,=\,(m_1,\ldots,m_n)$. In particular, the \emph{very
  large phase space potential} is
\begin{equation}
\label{eq:vlarge}
\Phi(\bt,\bs,\bu)\,=\,\sum_{g\,=\,0}^\infty\,\Phi_g(\bt,
\bs, \bu) ,
\end{equation}
where
\[
\Phi_g(\bt,\bs,\bu)\,=\,\ccor{}_g .
\]

The other two potentials are restrictions of $\Phi(\bt,\bs,\bu)$.
The large phase space potential
$\Phi(\bt)=\Phi(\bt,\mathbf{0},\mathbf{0})$ corresponds to setting
all the $\bs$ and $\bu$ variables to zero. The small phase space
potential is $$ \Phi(\bx)\,:=\,\Phi(x^0,\ldots,x^{r-2}), $$ where
$x^i\,:=\,t_0^i$ and all other variables are set to zero. The
function $\Phi(\bx)$ is the small phase space potential of the
$r$-spin \cft.  It will also be useful to define
\begin{equation}
Z(\bt,\bs,\bu)\,:=\,\exp(\Phi(\bt,\bs,\bu)).
\end{equation}

\subsection{The Euler vector field}

\

We begin with a differential equation arising from the grading.
The dimensions of the moduli spaces and cohomology classes induce
a grading on the potential function.
\begin{df}
The \emph{Euler vector field} $E$ is the differential operator
\begin{align*}
E\,=\,&
\sum_{\substack{a\,\geq\,0 \\ 0\,\leq\,m\,\leq\,r-2}}\,
(a - 1 + \frac{m}{r}) t_a^m\,\frac{\partial}{\partial t_a^m} \,
\\
& +\, \sum_{a\,\geq\,1}\,(a s_a\,\frac{\partial}{\partial
s_a}\,+\, (2 a -1) u_a\,\frac{\partial}{\partial u_a}).\,
\end{align*}
\end{df}

\begin{prop}
The very large phase space potential $\Phi(\bt,\bu,\bs)$ satisfies the
\emph{grading equation}

\begin{equation}\label{eq:grading}
E\,\Phi\,=\,(1 + \frac{1}{r})\lambda\,\frac{\partial}{\partial
\lambda}\,\Phi.
\end{equation}
\end{prop}

\begin{proof}
This follows from the definition of the potential, the dimensions
of the cohomology classes, and the dimensions of the moduli
spaces $\mgnrmbar$. It encodes the fact that intersection numbers
between cohomology classes vanish if the classes do not have
proper the dimension.
\end{proof}

This equation encodes the fact that the potential function is invariant under
the rescaling
$$
t_a^m\,\mapsto\,\epsilon^{a-1 + \frac{m}{r}}\,t_a^m,
u_a\,\mapsto\, \epsilon^{2 a - 1}\,u_a, s_a\,\mapsto\,\epsilon^a\,s_a,
\mathrm{\ and\ }
\lambda\,\mapsto\,\,\epsilon^{-1-\frac{1}{r}}\lambda.
$$

\begin{rem}
This grading shows that our small phase space  potential function
\emph{cannot} arise as the small phase space potential associated to the
Gromov-Witten invariants of a smooth, projective variety. This is because the
elements in $\chr$ have fractional dimension with respect to this
Euler vector field, whereas cohomology classes of a space always have
integral dimension.
\end{rem}

\subsection{The string equation and its cousins}

\

We begin by proving the analog of the string equation, the
dilaton equation, and a new equation arising from our identity on
the $\mu_1$ class.

\begin{thm}\label{thm:string}
Let $(g,n)$ be a pair of nonnegative integers such that $2 g - 2 + n > 0$.
The following identities are satisfied:
\begin{equation}
\cor{\tau_{0,0}\,\tau_{0,m_1}\,\tau_{0,m_2}\,\exp(\bs\,\cdot\,\bmu\,+\,
  \bu\,\cdot\, \bnu)\,}_0\,=\,\eta_{m_1 m_2}
\end{equation}
and
\begin{align*}
&\cor{\tau_{0,0}\tau_{a_1,m_1}\,\cdots\,\tau_{a_n,m_n}\,
  \exp(\bs\,\cdot\,\bmu\,+\,\bu\,\cdot\,\bnu)\,}_g \,= \\
&\sum_{i\,=\,1}^n\,
  \cor{\tau_{a_1,m_1}\,\cdots\,\tau_{a_i-1,m_i}\,\cdots\, \tau_{a_n,m_n}\,
\exp(\bs\,\cdot\,\bmu\,+\,\bu\,\cdot\,\bnu)\,}_g,
\end{align*}
where we assume that the terms in the sum containing $\tau_{a,m}$
with $a\,<\,0$ vanish.

These two equations are equivalent to the \emph{string (or
puncture) equation}
\begin{equation}
\cl_{-1} Z = 0,
\end{equation}
where
\begin{align*}
\cl_{-1}\,:=&\,- \frac{\partial}{\partial t_0^0} +
\sum_{0\,\leq\,m_1,m_2\,\leq\,r-2}\,\frac{1}{2} \,\eta_{m_1 m_2}
\,t_0^{m_1} t_0^{m_2}\,\\ +&\,\sum_{\substack{0\,\leq
\,m\,\leq\,r-2 \\ a\,\geq\,0}}\, t_{a+1}^m\,\frac{\partial
}{\partial  t_a^m} .\\
\end{align*}

Similarly, the following identities are satisfied:
\begin{equation}
\cor{\tau_{1,0}\,\exp({\bs \cdot  \bmu  + \bu \cdot \bnu})\,}_1 \,=
\,\frac{r-1}{24}
\end{equation}
and
\begin{align*}
&\cor{\tau_{1,0}\,\tau_{a_1,m_1}\,\cdots\,\tau_{a_n,m_n}\,
  \exp({ \bs \cdot  \bmu  + \bu \cdot \bnu})\,}_g \,= \\
  \,& (2 g - 2 + n)   \cor{\tau_{a_1,m_1}\,\cdots\,\tau_{a_n,m_n}
\exp({ \bs \cdot  \bmu  + \bu \cdot \bnu})\,}_g.
  \end{align*}
These two equations are equivalent to the \emph{dilaton equation}
\begin{equation}
\dil Z = 0,
\end{equation}
where
\begin{equation}
\dil\,=\, - \frac{\partial}{\partial t_1^0}\, +\, \sum_{\substack{ 0\, \leq\,
    m\,\leq\, r-2 \\ a\,\geq\,0}}\,t_a^m\,\frac{\partial}{\partial t_a^m}\,+\,
\lambda\,\frac{\partial}{\partial \lambda}\,+\,\frac{r-1}{24}.
\end{equation}
Finally, if $\cl_0$ denotes the differential operator
\begin{align*}
\cl_0\,:=&\,-(1 + \frac{1}{r})\,\frac{\partial}{\partial t_1^0}\,+\,
\sum_{\substack{0\,\leq\,m\,\leq\,r-2\\ a\,\geq\,0}} (a +
\frac{m\,+\,1}{r})\,t_a^m\,\frac{\partial}{\partial t_a^m}
\,+\,\frac{r^2-1}{24 r}
\\
\,+&\, \sum_{a\,\geq\,1}\,((2 a - 1)\, u_a\,\frac{\partial}{\partial
  u_a}\,+\,  a\,\frac{\partial}{\partial s_a}),\\
\end{align*}
then the following equation holds:
\begin{equation}
\cl_0\,Z\,=\,0.
\end{equation}

Furthermore, $[\cl_{0},\cl_{-1}]\,=\,\cl_{-1}$. When restricted
to the large phase space ($\bu\,=\,\bs\,=\,\mathbf{0}$), the
operators $\cl_0$ and $\cl_{-1}$ become the usual generators $L_0$ and $L_1$ in
the Virasoro Lie algebra.
\end{thm}

\begin{proof}
Recall that on the moduli space of  stable curves, $\mgnbar$, the
$\psi_i$ classes (where  $i\,=\,1,\ldots,n$) satisfy the equation
$\psi_i\,=\,\pi^*\,\psi_i\,+\,D_{i,n+1}$, where $D_{i,n+1}$ is the
image of the $i$-th canonical section and
$\pi\,:\,\M_{g,n+1}\,\to\,\M_{g,n}$. Let  $p$ be the forgetful
morphism $\mgnrmbar \,\to\, \mgnbar$, then since the  class
$\psi_i$ on $\mgnrmbar$ is the pullback via $p$ of the $\psi_i$
class on $\mgnbar$, one can lift the same formula to $\mgnrmbar$
to obtain $\psi_i\,=\,\pit^*\,\psi_i\,+\,D_{i,n+1}$, where this
equation is now regarded as being on $\mgnrmbar$, $D_{i,n+1}$ is
the pullback via $p$ of the divisors  with the same name on
$\mgnbar$, and $\pit$ is the forgetful morphism
$\M_{g,n+1}^{1/r,\bm\sqcup 0}\,\to\,\mgnrmbar$. Suppose that
$(g,n+1)\,\not=\,(0,3), \, (1,1)$. Using the lifting formula and
canceling trivial terms, we obtain
\[
\psi_j^a\,=\,\pit^*\,\psi_j^a\,+\,D_{j,n+1}\,\pit^*(\psi_j^{a-1}).
\]
Since $\pit^*\cv\,=\,\cv$, we have
\begin{align*}
\cor{\tau_{0,0}\,\tau_{a_1,m_1}\,\cdots\,\tau_{a_n,m_n}}_g \,&=\,
\lambda^{2g-2}\,
\int_{\M_{g,n+1}^{1/r,\bm}}\,\psi_1^{a_1}\,\ldots\,\psi_n^{a_n}\,
\cv \\ \,=\, \sum_{1\,\leq\,i\,\leq\,n}\,\lambda^{2 g-2}\,
\int_{\M_{g,n+1}^{1/r,\bm}}\, &
D_{i,n+1}\,\psi_{1}^{a_1}\,\cdots\,\psi_i^{a_i-1}\,\cdots\,\psi_n^{a_n}\,
\cv ,\\
\end{align*}
where the right hand side is understood to vanish if an exponent
is negative. Integration over the fiber of
$\M_{g,n+1}^{1/r,\bm\sqcup 0}\,\to\, \M_{g,n}^{1/r,\bm}$ yields
the desired result. Inclusion of the additional $\mu$ and $\nu$
classes into the correlators does not change the argument, since
$\pit^*\,\mu_i\,=\,\mu_i$, and similarly for $\nu_i$.

Finally, the exceptional cases follow from dimensional
considerations and the fact that on $\M_{0,3}^{1/r,\bm}$, $\cv$ is
the identity element in cohomology provided that
$m_1\,+\,m_2\,+\,m_3\,=\,r-2$. The dilaton equation is proved by a
similar analysis, where the exceptional case can be computed by
using the explicit presentation for $\psi_1$ on $\M_{1,1}^{1/r,0}$
to obtain
\begin{align}
\label{eq:M11CohFT}
\cor{\tau_{1,0}}_1\,&=\,\frac{1}{24}\,\eta^{m_+ m_-}\,\cor{\tau_{0,0}\,
  \tau_{0,m_+}\, \tau_{0,m_-}}_0\,
=\, \frac{r-1}{24}.
\end{align}

Finally, the equation
 $\cl_{0}\,Z\,=\,0$ is obtained by
combining the dilaton equation and the grading equation
(\ref{eq:grading}).
\end{proof}

The new relation for the $\mu_1$ class yields new equations between
correlators.

\begin{thm}
Let $\xi\,:\,\chr\,\to\,\chr$ be defined by
\[
{\xi_{m_+}}^{m_-}\,:=\,\frac{m_+
m_-}{r^2}\,\delta_{m_+\,+\,m_-,r-2},
\]
relative to the standard basis, and let
$\xi^{m_+,m_-}\,:=\,\eta^{m_+ m}\,{\xi_m}^{m_-}$. The following
differential equation holds:
\begin{align*}
\frac{\partial \Phi}{\partial s_1}\, &=
\frac{r^2-6 r + 6}{r}\,\frac{\partial \Phi}{\partial u_1}\,+\,
\sum_{\substack{{a\,\geq\,0} \\
{0\,\leq\,m_+,m_-\,\leq\,r-2}}}\,
\frac{1}{2} t_a^{m_+}\,{\xi_{m_+}}^{m_-}\,
\frac{\partial \Phi}{\partial t_{a+1}^{m-}}  \\
& - \sum_{m_+ , m_-}\,\frac{1}{2}\,\, \frac{\partial
  \Phi}{\partial t_0^{m_+}}\,\xi^{m_+,m_-}\,\frac{\partial
  \Phi}{\partial   t_0^{m_-}}
- \sum_{m_+, m_-}\,\frac{1}{4}\,\frac{\partial^2 \Phi}{\partial
  t_0^{m_+} \partial t_0^{m_-}}\,\xi^{m_+,m_-},
\end{align*}
where the summation over $m_+$ and $m_-$ runs over $0,\ldots,r-2$. This is
equivalent to the following relations between the correlators:

\begin{align*}
&\ccor{\tau_{a_1,m_1}\,\cdots\,\tau_{a_n,m_n}\,\mu_1}_g\,=\, \frac{r^2 - 6 r
  + 6}{r^2} \,\ccor{\tau_{a_1,m_1}\,\cdots\,\tau_{a_n,m_n}\,\nu_1}_g \\
&+ \sum_{\substack{{1\,\leq\,i\,\leq\,n} \\ m'_i}}\,
  \frac{1}{2}\,{\xi_{m_i}}^{m'_i}\, \ccor{\tau_{a_1,m_1}\, \cdots\,
  \tau_{a_i+1,m'_i}\,\cdots\,\tau_{a_n,m_n}}_g \\
&- \sum_{\substack{I_+\,\sqcup\, I_-\,=\,[n]\\ g_+ + g_- = g \\ m_+ , m_-
  }}\, \frac{1}{2}\,
\ccor{(\,\prod_{i\,\in\,I_+}{\tau_{a_i,m_i}})\tau_{0,m_+}\,}_{g_+}
\,\xi^{m_+,m_-}\, \ccor{\tau_{0,m_-}\,
  (\prod_{j\,\in\,I_-}\,\tau_{a_j,m_j})\,}_{g_-}\\
& - \sum_{m_+ , m_- }\,\frac{1}{4}\,
\ccor{\tau_{a_1,m_1}\,\cdots\, \tau_{a_n,m_n}\,\tau_{0,m_+} \tau_{0,
    m_-}}_{g-1}\, \xi^{m_+,m_-}, \\
\end{align*}
where we use the notation
$$ [n] = \{1,2,\ldots,n\} .
$$

The class $\nu_1$ which appears above is precisely $\lambda_1$. This class
vanishes on $\mgnrmbar$ for $g\,=\,0$.
\end{thm}

\begin{proof}
The proof follows from the facts that $p_*(\cv\,\exp( \bs \bmu ))$
forms a \cft, and that $\nu$ and $\psi$ are lifts of the analogous
classes on the moduli space of stable curves, and from
Proposition~\ref{prop:murelation}.
\end{proof}

\subsection{Topological recursion relations}

\

Topological recursion relations are relations between correlators which arise
from presentations of tautological classes in terms of boundary classes.

\begin{thm}\label{thm:trr0}
The following topological recursion relations hold in genus zero:

\begin{align*}
&\ccor{\tau_{a_1+1,m_1}\tau_{a_2,m_2}\tau_{a_3,m_3}}_0\,=\, \\
&\sum_{m_+, m_-}\, \ccor{\tau_{a_1,m_1}\tau_{0,m_+}\,}_0\,
\eta^{m_+ m_-}\,
\ccor{\tau_{0,m_-}\tau_{a_2,m_2}\tau_{a_3,m_3}}_0,\\
\end{align*}

which is equivalent to the differential equation
\[
\frac{\partial^3 \Phi_0}{\partial t_{a_1+1}^{m_1} \partial t_{a_2}^{m_2}
  t_{a_3}^{m_3}}\,=\, \sum_{m_+,m_-}\,\frac{\partial^2 \Phi_0}{\partial
  t_{a_1}^{m_1} \partial t_{0}^{m_+}}\, \eta^{m_+ m_-}\,\frac{\partial^3
  \Phi_0}{\partial t_0^{m_-}\,\partial t_{a_2}^{m_2}\,\partial
  t_{a_3}^{m_3}}.
\]

\end{thm}

\begin{proof}
On $\M_{0,n}^{1/r,\bm}$, the class $\psi_1$ can be written in terms of
boundary classes as
\[
\psi_1\,=\,\sum_{\substack{I_+\,\sqcup\,I_-\,=\,[n] \\ n-1, n \in I_+ \\
    1\,\in\,I_-}}\,\delta_{0;I_+}.
\]
This equation is obtained from lifting the analogous relation on $\M_{0,n}$.
The classes $\psi_i$ can be written similarly by applying an element of
the permutation group $S_n$.  The recursion relation follows from this
presentation and the restriction properties of the $\psi_i$ to the boundary
strata.
\end{proof}

\begin{thm}
The following topological recursion relation holds in genus one:
\begin{align*}
&\ccor{\tau_{a_1+1,m_1}}_1\,=\,
\frac{1}{24}\ccor{\tau_{a_1,m_1}\tau_{0,m_+}\tau_{0,m_-}}_0
\eta^{m_+ m_-} \\ &+\,\sum_{m_+,m_-}\, \ccor{\tau_{a_1,m_1}
\tau_{0,m_+} }_0 \,\eta^{m_+ m_-}\, \ccor{\tau_{0,m_-}}_1.
\end{align*}
This is equivalent to
\begin{align*}
\frac{\partial \Phi_1}{\partial t_{a_1+1}^{m_1}}\,=\,&\frac{1}{24}
  \sum_{m_+\,,\,m_-} \,\frac{\partial^3 \Phi_0}{\partial
  t_{a_1}^{m_1} \partial t_0^{m_+} \partial t_0^{m_-}} \,\eta^{m_+ m_-}\,\,\\
&+ \sum_{m_+, m_-}\,\frac{\partial^2 \Phi_0}{\partial t_{a_1}^{m_1}
  \partial t_0^{m_+}}\,\eta^{m_+ m_-}\,\frac{\partial \Phi_1}{\partial
  t_0^{m_-}}.\\
\end{align*}

The topological recursion relation for $\nu_1\,=\,\lambda_1$ from~\cite{KK}
\begin{align*}
&\ccor{\nu_1}_1\,=\\ \,
&\frac{1}{24}\,\sum_{m_+ , m_-}\,\eta^{m_+ m_-}
\ccor{\tau_{0,m_+} \tau_{0,m_-}}_0\\
\end{align*}
can be written as
\[
\frac{\partial \Phi_1}{\partial u_1}\,=\,\frac{1}{24}\,\sum_{m_+,m_-}\,
\eta^{m_+ m_-}\,\frac{\partial^2 \Phi_0}{\partial t_0^{m_+} \partial
  t_0^{m_-}}.
\]
\end{thm}
\begin{proof}
The proof of the first topological recursion relation arises from the
relation on $\M_{1,n}^{1/r,\bm}$
\[
\psi_1\,=\,\frac{1}{12}\,\delta_{irr}\,+\,
\sum_{\substack{I_+\sqcup I_-=[n] \\ n-1,n \in I_+ \\ 1 \in
I_-}}\, \delta_{0;I_+},
\]
which is obtained from lifting the analogous relation from
$\M_{1,n}$. The action of $S_n$ yields $\psi_i$. This, combined
with the restriction properties of the $\psi$ classes, yields the
desired result.  The second comes from the presentation of
$\lambda_1$ on $\M_{1,n}^{1/r,\bm}$
\[
\lambda_1\,=\,\frac{1}{12}\,\delta_{irr}
\]
and the restriction properties of $\lambda_1\,=\,\nu_1$.
\end{proof}
These two relations allow one to completely reduce the large phase
space potential to the small phase space potential in genus zero
and one. Combined with the previous equations, we can compute
$\Phi(\bt,\bs,\bu)$ in genus zero and one when we set
$s_i\,=\,u_i\,=\,0$ for all $i\,\geq\,2$.

The genus one potential satisfies an analog of the WDVV equation
due to Getzler \cite{Ge}, which arises from relations between
codimension-two boundary classes on $\M_{1,4}$. Using this
equation, Dubrovin and Zhang \cite{DuZh} showed that if the
Frobenius manifold is semisimple, then the genus-one potential is
determined by the Frobenius structure. Since the Frobenius
manifold structure on $(\chr,\eta)$ associated to $\Phi_0(\bx)$ is
known to be semisimple \cite{Du}, $\Phi_1(\bx)$ is determined. On
the other hand, the latter must vanish due to dimensional
considerations. Together with the topological recursion relations
in genus zero and one, we obtain the following corollary.

\begin{crl}
Let $\bs=\bu = \mathbf{0}$ so that we are on the large phase. Let
$v^m := \sum_{m=0}^{r-2} \ccor{\tau_{0,0} \tau_{0,l}}_0 \eta^{l
m}$. Let $\Delta(\bt)$  denote the matrix with entries
$\frac{\partial v^m}{\partial t_0^l}$ where $m,l = 0,\ldots,r-2$,
then
\begin{equation}
\label{eq:DiWiFormula}
\Phi_1(\bt) = \frac{1}{24}\ln \det \Delta(\bt).
\end{equation}
\end{crl}

\begin{proof}
Dijkgraaf and Witten \cite{DiWi} write down a formula (see
Theorem~15 of \cite{Ge2} for an explicit proof) for the large
phase space potential $\Phi_1(\bt)$ which in our case is
\[
\Phi_1(\bt) = \cor{\exp( \tau_0(v) )}_1 + \frac{1}{24}\ln \det
\Delta(\bt),
\]
where $v := \sum_{m=0}^{r-2} v^m e_m$. The term $\cor{\exp (
\tau_0(v) )}_1$ is equal to the small phase space potential
$\Phi_1(\bx)$, evaluated at $x^m = v^m$ for all $m =
0,\ldots,r-2$, but the small phase potential $\Phi_1(\bx)$
vanishes.
\end{proof}

In genus 2, there exist presentations of products of $\psi$
classes in terms of boundary classes \cite{Ge2, BPa} which give
rise to topological recursion relations, but, in general, they do
not allow one to reduce the computation of the large phase space
potential to the computation of the small phase space potential.

\section{The genus-zero large phase space potential} \label{sec:LPSCalc}

In this section, we compute the genus-zero, three- and four-point
correlators and show that they completely determine the
genus-zero, large phase space potential function $\Phi_0(\bt)$. We
roughly follow the outline provided by Witten \cite{W} and are
able to rigorously prove the validity of his computations, now
that the relevant moduli spaces and classes have been constructed.

\emph{Throughout the rest of this paper, we will consider only
 the large phase space potential $\Phi(\bt)$ (setting
 the other variables $u_a$ and $s_a$ of the very large phase space potential
 to zero).  We will also fix the coupling constant $\lambda$
 as $\lambda\,=\,\frac{1}{\sqrt{r}}$, since this is the value
 which is relevant to the generalized Witten conjecture.}

The following proposition rigorously demonstrates the formulas
from \cite{W}, but the idea of our proof is quite different, as it
uses our new relation for the $\mu_1$ class.
\begin{prop} \label{prop:potential}
The three-point and four-point correlators of the $r$-spin \cft\
are given by the following formulas:
\[
\cor{\tau_{0,m_1}\,
  \tau_{0,m_2}\, \tau_{0,m_3}}_0\,=\,\delta_{m_1\,+\,m_2\,+\,m_3,r-2} \]
and
\[
\cor{\tau_{0,m_1}\,\ldots\,\tau_{0,m_4}}_0\,=\,\frac{1}{r}\,
\mathrm{Min}_{1\leq i\leq 4}(m_i,r-1-m_i),
\]
where $\mathrm{Min}$ is minimum value.
\end{prop}

\begin{proof}
$\M_{0,n}^{1/r,\bm}$ is nonempty if and only if
$(2\,+\,\sum_i\,m_i)/r\,\in\,\nz$, where
$\bm\,=\,(m_1,\ldots,m_n)$ and $m_i\,=\,0\,\ldots\,r-1$ for all
$i$.  The genus zero correlators are given by
\[
\cor{\tau_{0,m_1}\,\ldots\,\tau_{0,m_n}}_0\,=\,r\,\int_{\M_{0,n}^{1/r,\bm}}\,\cv,
\]
where $\cv\,=\,c_D(-R^1\,\pi_*\,\ce_r)$ is the (top) Chern class
of degree
\[
D\,=\, -1\,+\,\frac{2}{r}\,+\,\frac{1}{r}\,\sum_i\,m_i.
\]
The class $\cv$ vanishes unless $m_i\,=\,0,\ldots\,r-2$ for all
$i$ by Theorem \ref{thm:gzero}. Furthermore, the correlator can
only be nonzero if $D\,=\,n-3$.

If $n\,=\,3$ then the dimensionality condition becomes
$m_1\,+\,m_2\,+\,m_3\,=\,r-2$, in which case $\cv$ is the
identity. This proves the first part of the proposition.

If $n\,=\,4$ then the dimensionality condition becomes
\[
m_1\,+\,\ldots\,+\,m_4\,=\, 2 r - 2.
\]
If this condition is satisfied then $\cv\,=\,\mu_1$. The correlator is
\[
\cor{\tau_{0,m_1}\,\ldots\,\tau_{0,m_4}}_0\,=\,
r\,\int_{\M_{0,4}^{1/r,\bm}}\,\mu_1.
\]

The right hand side can be computed using the relation for the
class $\mu_1$ in Proposition~\ref{prop:murelation}, which
becomes, in genus zero,

\begin{align*}
\mu_1\,&=\,\sum_{1\leq i\leq n}\,\frac{m_i(r-2-m_i)}{2 r^2} \psi_i
\\ & + \sum_{I_+ \subset [n]}\,\frac{r-1-(m_+ + 1) (m_-+1)}{2
r^2}\, \,\delta_{0;I_+},
\\
\end{align*}
where $m_+$ and $m_-$ are uniquely determined by the divisor
$\delta_{0:I_+}$. Let $\delta_{ij,kl}$ denote the divisor
$\delta_{0;\{\,i,j\,\}}$ on $\M_{0,4}^{1/r,\bm}$.

Plugging in this formula, one obtains (after doing a little case by case
analysis to write $m_+$ and $m_-$ in terms of $m_1,\ldots,m_4$)
\begin{align*}
\int_{\M_{0,4}^{1/r,\bm}}\,\mu_1\,=&\,\sum_{i\,=\,1}^4\, \frac{m_i
(r-2-m_i)}{2 r^2}\,\int_{\M_{0,4}^{1/r,\bm}}\,\psi_i\,\\ & +\,
\frac{r-1-(\chi_{12,34}\,+\,1) (r-1-\chi_{12,34})}{2 r^2}
\,\int_{\M_{0,4}^{1/r,\bm}}\, \delta_{12,34}\\ & +\,
\frac{r-1-(\chi_{13,24}\,+\,1) (r-1-\chi_{13,24})}{2 r^2}
\,\int_{\M_{0,4}^{1/r,\bm}}\, \delta_{13,24}\\ & +\,
\frac{r-1-(\chi_{14,23}\,+\,1) (r-1-\chi_{14,23})}{2 r^2}
\,\int_{\M_{0,4}^{1/r,\bm}}\, \delta_{14,23},\\
\end{align*}
where $\chi_{ij,kl}\,:=\,Min(m_i+m_j,m_k+m_l)$. Since each
$\delta_{ij,kl}$ is Poincar\'e dual to the (topological) homology
class represented by a point, one has
\[
r\,\int_{\M_{0,4}^{1/r,\bm}}\, \delta_{ij,kl}\,=\,1.
\]
Similarly, each class $\psi_i$ can be represented by
$\delta_{ij,kl}$ for some $i,j,k,l$.
Therefore, one obtains

\begin{align*}
r\,\int_{\M_{0,4}^{1/r,\bm}}\,\mu_1\,=
& \,\sum_{i\,=\,1}^4\,
\frac{m_i (r-2-m_i)}{2 r^2}  \\
& +\, \frac{r-1-(\chi_{12,34}\,+\,1) (r-1-\chi_{12,34})}{2 r^2}\\
& +\, \frac{r-1-(\chi_{13,24}\,+\,1) (r-1-\chi_{13,24})}{2 r^2} \\
& +\, \frac{r-1-(\chi_{14,23}\,+\,1) (r-1-\chi_{14,23})}{2 r^2}.
\end{align*}

The right hand side of this equation can be shown to be equal to
$\frac{1}{r}\,\mathrm{Min}_{1\leq i\leq 4}(m_i,r-1-m_i)$, an
elementary but not obvious identity.

\end{proof}

\begin{prop} \label{prop:unique} \cite{W}
The genus zero potential $\Phi_0(\bt)$ is completely determined by
$\cor{\tau_{a_1,m_1}\, \tau_{a_2,m_2}\,\tau_{a_3,m_3}\,\tau_{a_4,m_4}}_0$
and the fact that
\[ \cor{\tau_{a_1,m_1}\,
  \tau_{a_2,m_2}\,\tau_{a_3,m_3}}_0\,=\,\delta_{m_1+m_2+m_3,r-2}.
\]
\end{prop}

\begin{proof}
The large phase space, genus zero potential $\Phi_0(\bt)$ is
completely determined by its values on the small phase space by
the topological recursion relations. Let $\Phi_0(\bx)$ denote the
small phase space potential, which must satisfy the WDVV equation
since $\cv$ yields a \cft.  Furthermore, the grading equation
(\ref{eq:grading}) shows that the small phase space potential is a
polynomial in the variables $\{\,x^0,\,\ldots,\,x^{r-2}\,\} $ of
degree of at most $r+1$. One then performs an induction on the
degree of the polynomial to show that $\Phi_0(\bx)$ is uniquely
determined by the above data and the WDVV equation. The proof is
straightforward.
\end{proof}

\section{Gelfand-Dickey hierarchies and the generalized Witten conjecture}

In this section we present a mathematical formulation of the
generalized Witten conjecture, relating intersection theory on
$\mgnrmbar$ with Gelfand-Dickey integrable hierarchies and prove
this conjecture in various cases.

\subsection{Gelfand-Dickey hierarchies and their potentials}

In order to fix notation and normalization constants for the
generalized Witten conjecture, we recall the definition of the
Gelfand-Dickey hierarchies $KdV_r$ and   their special solutions.
A more detailed review can be found, for example, in~\cite{Lo,W}.

Fix an integer $r\ge 2$ and consider the space
\begin{equation}
  \label{eq:dops}
\cd=\{D^r-\sum_{m=0}^{r-2}u_m(x)D^m \}
\end{equation}
 of  differential operators
in $D=\frac{i}{\sqrt{r}} \frac{\partial}{\partial x}$ (the factor
$\frac{i}{\sqrt{r}}$ is added for convenience), where $u_m$ are
formal functional variables. For every operator $L\in \cd$ there
exists a unique pseudodifferential operator $$
L^{1/r}=D+\sum_{m>0}w_m D^{-m}, $$ such that $(L^{1/r})^r=L$. All
coefficients $w_m$ of $L^{1/r}$ are differential polynomials in
$u_0,u_1,\ldots,u_{r-2}$.

For a  pseudodifferential operator $Q=\displaystyle\sum_{m\ge
-n}v_mD^{-m}$, denote by $ Q_+=\displaystyle\sum_{m=-n}^0 v_m
D^{-m}$  its differential part, and consider the following
infinite family of differential equations on  $\cd$:
\begin{equation}
  \label{eq:kdv}
i\frac{\partial   L}{\partial   t^m_n} =
\frac{k_{n,m}}{\sqrt{r}}\left[(L^{n+\frac{m+1}{r}})_+,L\right],
\end{equation}
where the constants $$
k_{n,m}=\frac{(-1)^nr^{n+1}}{(m+1)(r+m+1)\ldots (nr+m+1)} $$ have
been introduced for convenience. It can be shown that the
corresponding flows on $\cd$ commute, and thus the following
definition makes sense.

\begin{df}
The infinite system~(\ref{eq:kdv}) of partial differential
equations with $r-1$ unknown functions $u_i(x,t^m_n)$,
$i=0,\ldots,r-2$, $m=0,\ldots,r-1$, $n\ge 0$ is called the
\emph{$r$-th Gelfand-Dickey hierarchy} or $KdV_r$.
\end{df}
The $KdV_2$ hierarchy is the usual Korteweg-de Vries hierarchy.

For $L=D^r-\displaystyle\sum_{m=0}^{r-2}u_m(x)D^m$, consider the functions
\begin{equation}
  \label{eq:v}
  v_n = -\frac{r}{n+1} \mathrm{res}(L^{1/r})^{n+1},
\end{equation}
where the residue of a pseudodifferential operator is defined as
the coefficient of $D^{-1}$. The functions $v_k$ can be expressed
in terms of $u_j$ by a triangular system of differential
polynomials. This means that $u_j$ can be expressed in terms of
$v_n$ in a similar way, and we may regard $v_0,v_1,\ldots,
v_{r-2}$ as a new system of coordinates for $\cd$.

\begin{df}
A formal power series $\Psi(\bt)$, in variables $t_n^m$,
$m=0,\ldots,r-2$, $n\ge 0$, is called a {\em   potential\/} of the
$\KdV_r$ hierarchy if it satisfies the following conditions:
\begin{enumerate}
\item
 $\Psi({\bf 0})=0$,
\item
 the functions
$$ v_m(\bt) = \frac{\partial^2 \Psi(\bt)}{\partial  t^0_0 \partial
t_0^m} $$ satisfy the equations~(\ref{eq:kdv}) with $x=t_0^0$ and
$u_j$ related to $v_m$ via~(\ref{eq:v}),
\item
 $\Psi(\bt)$ satisfies the
{\em   string equation\/}
\begin{equation}
  \label{eq:string}
\frac{\partial \Psi(\bt)}{\partial  t^0_0 } =
\frac{1}{2}\sum_{m,n=0}^{r-2} \eta_{mn}t_0^mt_0^n +
\sum_{k=0}^\infty \sum_{m=0}^{r-2}t_{k+1}^m \frac{\partial
\Psi(\bt)}{\partial  t^m_k},
\end{equation}
where $\eta_{mn}=\delta_{m+n,r-2}$.
\end{enumerate}
\end{df}
It can be shown that the potential $\Psi(\bt)$ is uniquely determined
by these  conditions (cf.~\cite{W3}).
\medskip

Finally, we introduce the semiclassical  limit  of the hierarchy
$\KdV_r$~(\ref{eq:kdv}) and its potential.

For a differential operator
$L=D^r-\displaystyle\sum_{m=0}^{r-2}u_m(x)D^m \in \cd$, denote by
$\widetilde{L}=p^r-\displaystyle\sum_{m=0}^{r-2}u_m(x)p^m$ the
polynomial in a formal variable $p$ obtained by replacing $D$ with
$p$. The commutator $[L,Q]$ of differential operators will be
replaced in~(\ref{eq:kdv}) by the Poisson bracket $$
\{\widetilde{L},\widetilde{Q}\} = \frac{\partial
\widetilde{L}}{\partial   p} \frac{\partial
\widetilde{Q}}{\partial x} -\frac{\partial \widetilde{Q}}{\partial
p} \frac{\partial \widetilde{L}}{\partial x} . $$

\begin{df}
{\em The semiclassical limit $\KdV^s_r$ of the $\KdV_r$
hierarchy\/}, is the system  of equations
\begin{equation}
\label{eq:kdv0}
\frac{\partial \widetilde{L}}{\partial  t^m_n}
=\frac{k_{m,n}}{r}
\left\{\widetilde{L^{n+\frac{m+1}{r}}},\widetilde{L}\right\}
\end{equation}
in unknown functions $u_0,\ldots,u_{r-2}$.
\end{df}

The corresponding potential function $\Psi_0(\bt)$ is defined as
the unique function satisfying the string
equation~(\ref{eq:string}) and the condition $\Psi_0({\bf 0})=0$,
and such that the functions $u_0,\ldots,u_{r-2}$ given
by~(\ref{eq:v}) and
\begin{equation}\label{eq:vzero}
v_m(\bt) = \frac{\partial^2 \Psi_0(\bt)}{\partial  t^0_0 \partial
t_0^m}
\end{equation}
satisfy the equations of the hierarchy~(\ref{eq:kdv0}).

\subsection{The Generalized Witten conjecture}\label{sec:kdvr}

Even before the moduli space $\mgnrmbar$ of $r$-spin curves was
constructed, Witten \cite{W} conjectured that these moduli spaces
would exist, and that intersection numbers on them would assemble
into the potential $\Psi(\bt)$ of the $KdV_r$ hierarchy. Now we
can give this conjecture the following mathematical formulation.

\begin{conj} \label{conj}
There exists an $r$-spin virtual (cohomology) class $\cv$  on
$\mgnrmbar$ satisfying Axioms 1---5 of Definition~\ref{axioms},
such that the large phase space potential $\Phi(\bt)$  of the
$r$-spin \cft~(\ref{eq:rlambda}) coincides with the potential
function $\Psi(\bt)$ of the $KdV_r$ hierarchy.
\end{conj}

Using results from Sections~\ref{vc}, ~\ref{sec:rr}, and
\ref{sec:LPSCalc}, we prove this conjecture in two special cases.

\begin{thm}
Conjecture~\ref{conj} holds for $r=2$ and arbitrary $g$.
\end{thm}

\begin{proof}
Theorem~\ref{th:rtwo} shows that when $r=2$ the class given by
(\ref{eq:cvrtwo}) satisfies the axioms of a virtual class, and
Corollary~\ref{cor:rtwo} implies that the large phase space
potential of the corresponding $2$-spin \cft\ is equal to the
generating function of tautological intersection numbers on
$\mgnbar$ (the large space potential of pure topological gravity).
By Kontsevich's theorem~\cite{Ko}, this generating function
coincides with the potential function of the Korteweg-de Vries
hierarchy, which is the same as the $KdV_2$ hierarchy.
\end{proof}

\begin{thm}\label{thm:conj0}
Conjecture~\ref{conj} holds for $g=0$ and arbitrary $r$.
\end{thm}

\begin{proof}
In this case, the conjecture means that the genus zero part
$\Phi_0(\bt)$ of the large phase space
potential~(\ref{eq:LPSPotential}) of the $r$-spin
\cft~(\ref{eq:rlambda}) coincides with the potential $\Psi_0(\bt)$
of the semiclassical limit of the $KdV_r$ hierarchy.

In genus zero the virtual class $\cv$ exists by
Theorem~\ref{thm:gzero}. From Theorem~\ref{thm:string} it follows
that the corresponding potential function $\Phi_0$ satisfies the
string equation~(\ref{eq:string}).

Because of the uniqueness of the potential function of the $KdV_r$
hierarchy (and its semiclassical approximation) all that remains
is the proof of the following proposition.
\end{proof}

\begin{prop}\label{prop:phi0}
The functions $u_m(\bt)$, $m=0,\ldots,r-2$, given
by~(\ref{eq:v}) and
\begin{equation}
v_m(\bt) = \frac{\partial^2 \Phi_0(\bt)}{\partial  t^0_0 \partial
t_0^m},
\end{equation}
satisfy the equations of the semiclassical limit of the $KdV_r$
hierarchy~(\ref{eq:kdv0}).
\end{prop}

\begin{proof}
By Proposition~\ref{prop:unique} there is a unique formal power
series $\Phi_0(\bt)$, of the proper grading, satisfying the
equations of Proposition~\ref{prop:potential}, WDVV, and the
genus-zero topological recursion relations.  Witten \cite{W} shows
by a straightforward computation that any such power series yields
a solution of the semiclassical limit of the $KdV_r$ hierarchy.
\end{proof}

\begin{crl}
\label{cor:frob} The Frobenius manifold structure on
$(\chr,\eta)$, defined by Theorems~\ref{th:rcft}
and~\ref{thm:gzero}, is isomorphic to the Frobenius structure on
the base of the versal deformation of the  $A_{r-1}$ singularity.
\end{crl}
\begin{proof}
The proof follows from Theorem~\ref{thm:conj0} and the fact that
the potential of the Frobenius structure on the base of the versal
deformation of the  $A_{r-1}$ singularity is equal to the
potential $\Psi_0$ of the semiclassical limit of the $KdV_r$
hierarchy (cf.~\cite{Du}).
\end{proof}

\medskip

\begin{rems}
\

\begin{enumerate}
\item
The generalized Witten conjecture, as it is stated here, should be
viewed as a refinement of Witten's original formulation of his
conjecture \cite{W}, since it is not clear that his construction
yields a class with the desired factorization properties.

\item
The coincidence of Frobenius structures given by
Corollary~\ref{cor:frob} appears to be a genus zero manifestation
of some mirror phenomenon \cite{Ma2}, relating the moduli space of
$r$-spin curves and singularities of type $A_{r-1}$.

\item
There is additional evidence for       Conjecture~\ref{conj} in
genus one for arbitrary $r$. Witten \cite{W} states that the
formula (\ref{eq:M11Calc}) for the intersection numbers when $g=1$
can be derived from the conjecture for all $r \geq  2$.
Furthermore, when $r\leq 4$, it can be shown that
equation~(\ref{eq:DiWiFormula}) holds for the genus-one part of
the potential of the $\KdV_r$ hierarchy (cf.~\cite{DiWi,DuZh}).

\item The fact that $\Psi(\bt)$ (in all genera) is independent of the
variables $t_n^{r-1}$ for all $n\,\geq\,0$ is consistent with Axiom 4
(Vanishing) of the virtual class.
\end{enumerate}
\end{rems}
\medskip

The exponential of the $\KdV$ potential function is called a
$\tau$-function and can be defined as the unique function $Z(\bt)$
annihilated by certain differential operators $L_i$, $i\ge -1$,
generating (a part of) the Virasoro Lie algebra. This gives an
alternate formulation of the original Witten conjecture.
Similarly, the exponential $Z(\bt)$ of the $\KdV_r$ potential is
annihilated by a series of differential operators which forms a
so-called $W_r^+$-algebra~\cite{AvM,Krich} (part of which forms a
subalgebra isomorphic to (half of) the Virasoro algebra). Thus we
obtain an alternate formulation of the generalized Witten
conjecture.

\begin{conj} (W-algebra conjecture)
There exist a collection of differential operators forming a
$W_r^+$ algebra (in which the generators $\{ L_n \}_{n \geq -1}$
of the Virasoro algebra form a subset) which annihilates and
completely determines $Z(\bt)$.
\end{conj}

This conjecture can be regarded as the $\KdV_r$ analog of a
refinement of the Virasoro conjecture \cite{EHX}. When $r=2$, this
conjecture reduces to the usual Virasoro highest weight condition.

\bibliographystyle{amsplain}

\providecommand{\bysame}{\leavevmode\hbox to3em{\hrulefill}\thinspace}

\end{document}